\documentclass[submission,copyright,creativecommons]{eptcs}
\providecommand{\event}{Applied Category Theory, 2021} 
\usepackage{breakurl}             
\usepackage{etex}

\usepackage[pdftex,dvipsnames]{xcolor}  
\usepackage{amsmath,amsthm,amssymb,amsfonts,graphicx}
\usepackage{stmaryrd}
\usepackage{tikz}
\usepackage{proof}
\usepackage{enumerate}
\usepackage{mathtools}
\usepackage{xypic}
\usepackage{lscape}
\usepackage{multicol} 
\usepackage{leftidx}
\usepackage{bbm}
\usepackage{rotating}
\usepackage{extarrows}
\usepackage{cmll}

\usepackage{array}   
\usepackage{float}

\newcolumntype{L}{>{$}l<{$}} 

\tikzstyle{strings}=[baseline={([yshift=-.5ex]current bounding box.center)}]

\tikzset{every picture/.append style={scale=.5}, transform shape, strings}

\tikzset{%
symbol/.style={%
draw=none,
every to/.append style={%
edge node={node [sloped, allow upside down, auto=false]{$#1$}}}
}
}  

\usetikzlibrary{shapes.geometric}
\usetikzlibrary{patterns}
\usetikzlibrary{fit}
\usetikzlibrary{positioning}
\usetikzlibrary{calc}
\usetikzlibrary{arrows}
\usetikzlibrary{decorations.markings}
\usetikzlibrary{decorations.pathreplacing}
\usetikzlibrary{shapes}

\pgfdeclarelayer{nodelayer}
\pgfdeclarelayer{edgelayer}
\pgfsetlayers{edgelayer,nodelayer,main}

\tikzset{simple/.style={}}
\tikzset{nothing/.style={outer sep=-3.4pt}}
\tikzset{map/.style={draw,fill=white, rectangle}}

\tikzset{dot/.style={thick, fill=black, circle, scale=1, inner sep = .05cm}}

\tikzset{oa/.style={draw, scale=0.9,minimum height=.1cm,circle,append after command={
[shorten >=\pgflinewidth, shorten <=\pgflinewidth,]
(\tikzlastnode.north) edge (\tikzlastnode.south)
(\tikzlastnode.east) edge (\tikzlastnode.west)
} } }

\tikzset{ox/.style={draw, scale=0.9,minimum height=.1cm,circle,append after command={
[shorten >=\pgflinewidth, shorten <=\pgflinewidth,]
(\tikzlastnode.north west) edge (\tikzlastnode.south east)
(\tikzlastnode.north east) edge (\tikzlastnode.south west) } } }

\tikzstyle{none}=[inner sep=-1pt]
\tikzstyle{circle}=[shape=circle,draw]
\tikzstyle{black}=[fill=black, draw=black, shape=circle]
\tikzstyle{onehalfcircle}=[fill=white, draw=black, shape=circle, scale=1.5]
\tikzstyle{twocircle}=[fill=white, draw=black, shape=circle, scale=2]
\tikzstyle{blackmap}=[fill=black, draw=black, shape=rectangle]
\tikzstyle{oa}=[draw = black, scale=0.9, minimum height=.1cm, shape= circle, append after command={
[shorten >=\pgflinewidth, shorten <=\pgflinewidth,]
(\tikzlastnode.north) edge (\tikzlastnode.south)
(\tikzlastnode.east) edge (\tikzlastnode.west)
} ]
\tikzstyle{ox}=[draw = black, scale=0.9, minimum height=.1cm, shape= circle, append after command={
    [shorten >=\pgflinewidth, shorten <=\pgflinewidth,]
(\tikzlastnode.north west) edge (\tikzlastnode.south east)
(\tikzlastnode.north east) edge (\tikzlastnode.south west)
} ]

\tikzstyle{filled}=[-, fill=black]

\tikzset{wires/.style={}}

\tikzset{box/.style={inner sep=0pt, thick, draw=black, text height=1.5ex, text depth=.25ex, 
text centered, minimum height=3em, anchor=center}}

\usepackage{pgfplots}
\usepgfplotslibrary{external} 
\tikzset{external/up to date check={md5}}
\tikzset{external/mode=convert with system call} 
\tikzsetexternalprefix{build/}

\tikzset{external/export=false}

\renewcommand{\u}{{\sf{u}}}
\renewcommand{\v}{{\sf{v}}}

\newcommand{\X}{\mathbb{X}}

\newcommand{\C}{\mathbb{C}}

\newcommand{\U}{\mathbb{U}}

\newcommand{\m}{{\sf m}}

\newcommand{\ox}{\otimes}

\newcommand{\oa}{\oplus}
\newcommand{\op}{\mathsf{op}}

\newcommand{\mx}{\mathsf{mx}}

\newcommand{\FMat}{\mathsf{FMat}}

\newcommand{\Core}{\mathsf{Core}}
\newcommand{\Unitary}{\mathsf{Unitary}}

\newcommand{\dashvv}{\dashv \!\!\!\!\! \dashv}

\renewcommand{\phi}{\varphi}

\newcounter{dummy} 
\numberwithin{dummy}{section}

\newtheorem{lemma}[dummy]{Lemma}  
\newtheorem{theorem}[dummy]{Theorem}
\newtheorem{definition}[dummy]{Definition}

\newtheorem{proposition}[dummy]{Proposition} 
\newtheorem{corollary}[dummy]{Corollary} 
\numberwithin{equation}{dummy}

\newlength{\llcfoo}

\makeatletter

\newdimen\w@dth

\def\setw@dth#1#2{\setbox\z@\hbox{\scriptsize $#1$}\w@dth=\wd\z@
\setbox\@ne\hbox{\scriptsize $#2$}\ifnum\w@dth<\wd\@ne \w@dth=\wd\@ne \fi
\advance\w@dth by 1.2em}

\def\t@^#1_#2{\allowbreak\def\n@one{#1}\def\n@two{#2}\mathrel
{\setw@dth{#1}{#2}
\mathop{\hbox to \w@dth{\rightarrowfill}}\limits
\ifx\n@one\empty\else ^{\box\z@}\fi
\ifx\n@two\empty\else _{\box\@ne}\fi}}
\def\t@@^#1{\@ifnextchar_ {\t@^{#1}}{\t@^{#1}_{}}}

\def\t@left^#1_#2{\def\n@one{#1}\def\n@two{#2}\mathrel{\setw@dth{#1}{#2}
\mathop{\hbox to \w@dth{\leftarrowfill}}\limits
\ifx\n@one\empty\else ^{\box\z@}\fi
\ifx\n@two\empty\else _{\box\@ne}\fi}}
\def\t@@left^#1{\@ifnextchar_ {\t@left^{#1}}{\t@left^{#1}_{}}}

\def\two@^#1_#2{\def\n@one{#1}\def\n@two{#2}\mathrel{\setw@dth{#1}{#2}
\mathop{\vcenter{\hbox to \w@dth{\rightarrowfill}\kern-1.7ex
                 \hbox to \w@dth{\rightarrowfill}}%
       }\limits
\ifx\n@one\empty\else ^{\box\z@}\fi
\ifx\n@two\empty\else _{\box\@ne}\fi}}
\def\tw@@^#1{\@ifnextchar_ {\two@^{#1}}{\two@^{#1}_{}}}

\def\tofr@^#1_#2{\def\n@one{#1}\def\n@two{#2}\mathrel{\setw@dth{#1}{#2}
\mathop{\vcenter{\hbox to \w@dth{\rightarrowfill}\kern-1.7ex
                 \hbox to \w@dth{\leftarrowfill}}%
       }\limits
\ifx\n@one\empty\else ^{\box\z@}\fi
\ifx\n@two\empty\else _{\box\@ne}\fi}}
\def\t@fr@^#1{\@ifnextchar_ {\tofr@^{#1}}{\tofr@^{#1}_{}}}

\newdimen\W@dth
\def\setW@dth#1#2{\setbox\z@\hbox{$#1$}\W@dth=\wd\z@
\setbox\@ne\hbox{$#2$}\ifnum\W@dth<\wd\@ne \W@dth=\wd\@ne \fi
\advance\W@dth by 1.2em}

\def\T@^#1_#2{\allowbreak\def\N@one{#1}\def\N@two{#2}\mathrel
{\setW@dth{#1}{#2}
\mathop{\hbox to \W@dth{\rightarrowfill}}\limits
\ifx\N@one\empty\else ^{\box\z@}\fi
\ifx\N@two\empty\else _{\box\@ne}\fi}}
\def\T@@^#1{\@ifnextchar_ {\T@^{#1}}{\T@^{#1}_{}}}

\def\T@left^#1_#2{\def\N@one{#1}\def\N@two{#2}\mathrel{\setW@dth{#1}{#2}
\mathop{\hbox to \W@dth{\leftarrowfill}}\limits
\ifx\N@one\empty\else ^{\box\z@}\fi
\ifx\N@two\empty\else _{\box\@ne}\fi}}
\def\T@@left^#1{\@ifnextchar_ {\T@left^{#1}}{\T@left^{#1}_{}}}

\def\Tofr@^#1_#2{\def\N@one{#1}\def\N@two{#2}\mathrel{\setW@dth{#1}{#2}
\mathop{\vcenter{\hbox to \W@dth{\rightarrowfill}\kern-1.7ex
                 \hbox to \W@dth{\leftarrowfill}}%
       }\limits
\ifx\N@one\empty\else ^{\box\z@}\fi
\ifx\N@two\empty\else _{\box\@ne}\fi}}
\def\T@fr@^#1{\@ifnextchar_ {\Tofr@^{#1}}{\Tofr@^{#1}_{}}}

\def\Two@^#1_#2{\def\N@one{#1}\def\N@two{#2}\mathrel{\setW@dth{#1}{#2}
\mathop{\vcenter{\hbox to \W@dth{\rightarrowfill}\kern-1.7ex
                 \hbox to \W@dth{\rightarrowfill}}%
       }\limits
\ifx\N@one\empty\else ^{\box\z@}\fi
\ifx\N@two\empty\else _{\box\@ne}\fi}}
\def\Tw@@^#1{\@ifnextchar_ {\Two@^{#1}}{\Two@^{#1}_{}}}

\def\to{\@ifnextchar^ {\t@@}{\t@@^{}}}
\def\from{\@ifnextchar^ {\t@@left}{\t@@left^{}}}
\def\tofro{\@ifnextchar^ {\t@fr@}{\t@fr@^{}}}
\def\To{\@ifnextchar^ {\T@@}{\T@@^{}}}
\def\From{\@ifnextchar^ {\T@@left}{\T@@left^{}}}
\def\Two{\@ifnextchar^ {\Tw@@}{\Tw@@^{}}}
\def\Tofro{\@ifnextchar^ {\T@fr@}{\T@fr@^{}}}

\makeatother

\newcommand{\linmonwtik} {\begin{tikzpicture}
	\begin{pgfonlayer}{nodelayer}
		\node [style=none] (0) at (-2.7, 1.17) {};
		\node [style=none] (1) at (-1.85, 1.17) {};
		\node [style=none] (2) at (-2, 1.35) {};
		\node [style=none] (3) at (-2, 1) {};
		\node [style=none] (4) at (-1.85, 1) {};
		\node [style=none] (5) at (-1.85, 1.35) {};
		\node [style=circle, scale=0.6] (6) at (-2.35, 1.45) {};
		\node [style=none] (7) at (-1.6, 1.17) {};
		\node [style=none] (8) at (-2.95, 1.17) {};
	\end{pgfonlayer}
	\begin{pgfonlayer}{edgelayer}
		\draw (2.center) to (3.center);
		\draw (5.center) to (4.center);
		\draw (0.center) to (1.center);
	\end{pgfonlayer}
\end{tikzpicture}}

\newcommand{\expmonwtik} {\begin{tikzpicture}
	\begin{pgfonlayer}{nodelayer}
		\node [style=none] (0) at (-2.8, 1.17) {};
		\node [style=none] (1) at (-1.85, 1.17) {};
		\node [style=none] (2) at (-2, 1.35) {};
		\node [style=none] (3) at (-2, 1) {};
		\node [style=none] (4) at (-1.85, 1) {};
		\node [style=none] (5) at (-1.85, 1.35) {};
		\node [style=circle, scale=0.6] (6) at (-2.35, 1.55) {};
		\node [style=map, scale=1.7, fill opacity=0] (9) at (-2.35, 1.55) {};
		\node [style=none] (7) at (-1.6, 1.17) {};
		\node [style=none] (8) at (-3.05, 1.17) {};
	\end{pgfonlayer}
	\begin{pgfonlayer}{edgelayer}
		\draw (2.center) to (3.center);
		\draw (5.center) to (4.center);
		\draw (0.center) to (1.center);
	\end{pgfonlayer}
\end{tikzpicture}}

\newcommand{\dagmonwtik} {\begin{tikzpicture}
	\begin{pgfonlayer}{nodelayer}
		\node [style=none] (0) at (-2.7, 1.17) {};
		\node [style=none] (1) at (-1.85, 1.17) {};
		\node [style=none] (2) at (-2, 1.35) {};
		\node [style=none] (3) at (-2, 1) {};
		\node [style=none] (4) at (-1.85, 1) {};
		\node [style=none] (5) at (-1.85, 1.35) {};
		\node [style=circle, scale=0.6] (6) at (-2.25, 1.45) {};
		\node [style=none] (7) at (-1.6, 1.17) {};
		\node [style=none] (8) at (-2.95, 1.17) {};
		\node [style=none, scale=1.5] (9) at (-2.55, 1.45) {$\dag$};
	\end{pgfonlayer}
	\begin{pgfonlayer}{edgelayer}
		\draw (2.center) to (3.center);
		\draw (5.center) to (4.center);
		\draw (0.center) to (1.center);
	\end{pgfonlayer}
\end{tikzpicture}}

\newcommand{\lincomonwtik} {\begin{tikzpicture} 
	\begin{pgfonlayer}{nodelayer}
		\node [style=none] (0) at (-2.7, 1.17) {};
		\node [style=none] (1) at (-1.85, 1.17) {};
		\node [style=none] (2) at (-2, 1.35) {};
		\node [style=none] (3) at (-2, 1) {};
		\node [style=none] (4) at (-1.85, 1) {};
		\node [style=none] (5) at (-1.85, 1.35) {};
		\node [style=circle, scale=0.6] (6) at (-2.35, 0.9) {};
		\node [style=none] (7) at (-1.6, 1.17) {};
		\node [style=none] (8) at (-2.95, 1.17) {};
	\end{pgfonlayer}
	\begin{pgfonlayer}{edgelayer}
		\draw (2.center) to (3.center);
		\draw (5.center) to (4.center);
		\draw (0.center) to (1.center);
	\end{pgfonlayer}
\end{tikzpicture}}

\newcommand{\dagcomonwtik} {\begin{tikzpicture} 
	\begin{pgfonlayer}{nodelayer}
		\node [style=none] (0) at (-2.7, 1.17) {};
		\node [style=none] (1) at (-1.85, 1.17) {};
		\node [style=none] (2) at (-2, 1.35) {};
		\node [style=none] (3) at (-2, 1) {};
		\node [style=none] (4) at (-1.85, 1) {};
		\node [style=none] (5) at (-1.85, 1.35) {};
		\node [style=circle, scale=0.5] (6) at (-2.2, 0.9) {};
		\node [style=none, scale=1.4] (9) at (-2.5, 0.8) {$\dag$};
		\node [style=none] (7) at (-1.6, 1.17) {};
		\node [style=none] (8) at (-2.95, 1.17) {};
	\end{pgfonlayer}
	\begin{pgfonlayer}{edgelayer}
		\draw (2.center) to (3.center);
		\draw (5.center) to (4.center);
		\draw (0.center) to (1.center);
	\end{pgfonlayer}
\end{tikzpicture}}

\newcommand{\lincomonbtik} {\begin{tikzpicture} 
	\begin{pgfonlayer}{nodelayer}
		\node [style=none] (0) at (-2.7, 1.17) {};
		\node [style=none] (1) at (-1.85, 1.17) {};
		\node [style=none] (2) at (-2, 1.35) {};
		\node [style=none] (3) at (-2, 1) {};
		\node [style=none] (4) at (-1.85, 1) {};
		\node [style=none] (5) at (-1.85, 1.35) {};
		\node [style=circle, scale=0.6, fill=black] (6) at (-2.35, 0.9) {};
		\node [style=none] (7) at (-1.6, 1.17) {};
		\node [style=none] (8) at (-2.95, 1.17) {};
	\end{pgfonlayer}
	\begin{pgfonlayer}{edgelayer}
		\draw (2.center) to (3.center);
		\draw (5.center) to (4.center);
		\draw (0.center) to (1.center);
	\end{pgfonlayer}
\end{tikzpicture}}

\newcommand{\lincomonwtritik} {\begin{tikzpicture} 
	\begin{pgfonlayer}{nodelayer}
		\node [style=none] (0) at (-2.7, 1.17) {};
		\node [style=none] (1) at (-1.85, 1.17) {};
		\node [style=none] (2) at (-2, 1.35) {};
		\node [style=none] (3) at (-2, 1) {};
		\node [style=none] (4) at (-1.85, 1) {};
		\node [style=none] (5) at (-1.85, 1.35) {};
		\node [style=none] (6) at (-2.2, 1) {};
		\node [style=none] (7) at (-2.5, 1) {};
		\node [style=none] (8) at (-2.35, 0.82) {};
		\node [style=none] (9) at (-1.6, 1.17) {};
		\node [style=none] (10) at (-2.95, 1.17) {};
	\end{pgfonlayer}
	\begin{pgfonlayer}{edgelayer}
		\draw (2.center) to (3.center);
		\draw (5.center) to (4.center);
		\draw (0.center) to (1.center);
		\draw (6.center) -- (7.center) -- (8.center) -- (6.center);
	\end{pgfonlayer}
\end{tikzpicture}}

\newcommand{\linbialgwtik} {\begin{tikzpicture}
	\begin{pgfonlayer}{nodelayer}
		\node [style=none] (0) at (-2.8, 1.17) {};
		\node [style=none] (1) at (-1.85, 1.17) {};
		\node [style=none] (2) at (-2, 1.35) {};
		\node [style=none] (3) at (-2, 1) {};
		\node [style=none] (4) at (-1.85, 1) {};
		\node [style=none] (5) at (-1.85, 1.35) {};
		\node [style=none] (6) at (-2.2, 1) {};
		\node [style=none] (7) at (-2.5, 1) {};
		\node [style=none] (8) at (-2.35, 0.82) {};
		\node [style=none] (9) at (-1.6, 1.17) {};
		\node [style=none] (10) at (-3.05, 1.17) {};
		\node [style=circle, scale=0.6] (11) at (-2.35, 1.45) {};
	\end{pgfonlayer}
	\begin{pgfonlayer}{edgelayer}
		\draw (2.center) to (3.center);
		\draw (5.center) to (4.center);
		\draw (0.center) to (1.center);
		\draw (6.center) -- (7.center) -- (8.center) -- (6.center);
	\end{pgfonlayer}
\end{tikzpicture}}

\newcommand{\expbialgwtik} {\begin{tikzpicture}
	\begin{pgfonlayer}{nodelayer}
		\node [style=none] (0) at (-2.7, 1.17) {};
		\node [style=none] (1) at (-1.85, 1.17) {};
		\node [style=none] (2) at (-2, 1.35) {};
		\node [style=none] (3) at (-2, 1) {};
		\node [style=none] (4) at (-1.85, 1) {};
		\node [style=none] (5) at (-1.85, 1.35) {};
		\node [style=none] (6) at (-2.2, 1) {};
		\node [style=none] (7) at (-2.5, 1) {};
		\node [style=none] (8) at (-2.35, 0.82) {};
		\node [style=none] (9) at (-1.6, 1.17) {};
		\node [style=none] (10) at (-2.95, 1.17) {};
		\node [style=circle, scale=0.6 ] (12) at (-2.35, 1.55) {};
		\node [style=map, scale=1.7, fill opacity=0] (11) at (-2.35, 1.55) {};
	\end{pgfonlayer}
	\begin{pgfonlayer}{edgelayer}
		\draw (2.center) to (3.center);
		\draw (5.center) to (4.center);
		\draw (0.center) to (1.center);
		\draw (6.center) -- (7.center) -- (8.center) -- (6.center);
	\end{pgfonlayer}
\end{tikzpicture}}

\newcommand{\dagbialgwtik} {\begin{tikzpicture}
	\begin{pgfonlayer}{nodelayer}
		\node [style=none] (0) at (-2.7, 1.17) {};
		\node [style=none] (1) at (-1.85, 1.17) {};
		\node [style=none] (2) at (-2, 1.35) {};
		\node [style=none] (3) at (-2, 1) {};
		\node [style=none] (4) at (-1.85, 1) {};
		\node [style=none] (5) at (-1.85, 1.35) {};
		\node [style=none] (6) at (-2.2, 1) {};
		\node [style=none] (7) at (-2.5, 1) {};
		\node [style=none] (8) at (-2.35, 0.82) {};
		\node [style=none] (9) at (-1.6, 1.17) {};
		\node [style=none] (10) at (-2.95, 1.17) {};
		\node [style=circle, scale=0.5] (11) at (-2.25, 1.45) {};		
		\node [style=none, scale=1.5] (12) at (-2.55, 1.45) {$\dagger$};
		\node [style=none] (13) at (-2.25, 0.5) {}; 
	\end{pgfonlayer}
	\begin{pgfonlayer}{edgelayer}
		\draw (2.center) to (3.center);
		\draw (5.center) to (4.center);
		\draw (0.center) to (1.center);
		\draw (6.center) -- (7.center) -- (8.center) -- (6.center);
	\end{pgfonlayer}
\end{tikzpicture}}

\newcommand{\linbialgbtik} {\begin{tikzpicture}
	\begin{pgfonlayer}{nodelayer}
		\node [style=none] (0) at (-2.7, 1.17) {};
		\node [style=none] (1) at (-1.85, 1.17) {};
		\node [style=none] (2) at (-2, 1.35) {};
		\node [style=none] (3) at (-2, 1) {};
		\node [style=none] (4) at (-1.85, 1) {};
		\node [style=none] (5) at (-1.85, 1.35) {};
		\node [style=none] (6) at (-2.2, 1) {};
		\node [style=none] (7) at (-2.5, 1) {};
		\node [style=none] (8) at (-2.35, 0.82) {};
		\node [style=none] (9) at (-1.6, 1.17) {};
		\node [style=none] (10) at (-2.95, 1.17) {};
		\node [style=circle, scale=0.6, fill=black] (11) at (-2.35, 1.45) {};
	\end{pgfonlayer}
	\begin{pgfonlayer}{edgelayer}
		\draw (2.center) to (3.center);
		\draw (5.center) to (4.center);
		\draw (0.center) to (1.center);
		\draw[fill=black] (6.center) -- (7.center) -- (8.center) -- (6.center);
	\end{pgfonlayer}
\end{tikzpicture}}

\newcommand{\expbialgbtik} {\begin{tikzpicture}
	\begin{pgfonlayer}{nodelayer}
		\node [style=none] (0) at (-2.7, 1.17) {};
		\node [style=none] (1) at (-1.85, 1.17) {};
		\node [style=none] (2) at (-2, 1.35) {};
		\node [style=none] (3) at (-2, 1) {};
		\node [style=none] (4) at (-1.85, 1) {};
		\node [style=none] (5) at (-1.85, 1.35) {};
		\node [style=none] (6) at (-2.2, 1) {};
		\node [style=none] (7) at (-2.5, 1) {};
		\node [style=none] (8) at (-2.35, 0.82) {};
		\node [style=none] (9) at (-1.6, 1.17) {};
		\node [style=none] (10) at (-2.95, 1.17) {};
		\node [style=circle, scale=0.6, fill=black ] (12) at (-2.35, 1.55) {};
		\node [style=map, scale=1.7, fill opacity=0] (11) at (-2.35, 1.55) {};
	\end{pgfonlayer}
	\begin{pgfonlayer}{edgelayer}
		\draw (2.center) to (3.center);
		\draw (5.center) to (4.center);
		\draw (0.center) to (1.center);
		\draw (6.center) -- (7.center) -- (8.center) -- (6.center);
	\end{pgfonlayer}
\end{tikzpicture}}

\newcommand{\mulmap}[2]{
\begin{tikzpicture}[scale=#1, rotate=180]
	\begin{pgfonlayer}{nodelayer}
		\node [style=circle,scale=0.4, fill=#2] (0) at (-1, 3) {};
		\node [style=none] (1) at (-1.25, 2.75) {};
		\node [style=none] (2) at (-0.75, 2.75) {};
		\node [style=none] (3) at (-1, 3.25) {};
	\end{pgfonlayer}
	\begin{pgfonlayer}{edgelayer}
		\draw [style=none] (3.center) to (0);
		\draw [style=none, bend right] (0.center) to (1.center);
		\draw [style=none, bend left] (0.center) to (2.center);
	\end{pgfonlayer}
\end{tikzpicture}
}

\newcommand{\leftaction}[2]{
	\begin{tikzpicture}[scale=#1]
		\begin{pgfonlayer}{nodelayer}
			\node [style=none] (0) at (0.5, 0.25) {};
			\node [style=none] (1) at (0.5, 0.75) {};
			\node [style=none] (2) at (0.25, 0.5) {};
			\node [style=none] (3) at (0.5, -0.5) {};
			\node [style=none] (4) at (0.5, 1.5) {};
			\node [style=none] (5) at (-0.5, 1.5) {};
			\node [style=none] (6) at (0.5, 0.75) {};
			\node [style=none] (7) at (0.5, 0.25) {};
		\end{pgfonlayer}
		\begin{pgfonlayer}{edgelayer}
			\draw (3.center) to (0.center);
			\draw (1.center) to (4.center);
			\draw [in=-90, out=180, looseness=1.25] (2.center) to (5.center);
			\draw [bend right=90, looseness=2.00] (6.center) to (7.center);
			\draw (7.center) to (6.center);
		\end{pgfonlayer}
	\end{tikzpicture}}

\newcommand{\rightaction}[2]{
	\begin{tikzpicture}[scale=#1, xscale=-1]
		\begin{pgfonlayer}{nodelayer}
			\node [style=none] (0) at (0.5, 0.25) {};
			\node [style=none] (1) at (0.5, 0.75) {};
			\node [style=none] (2) at (0.25, 0.5) {};
			\node [style=none] (3) at (0.5, -0.5) {};
			\node [style=none] (4) at (0.5, 1.5) {};
			\node [style=none] (5) at (-0.5, 1.5) {};
			\node [style=none] (6) at (0.5, 0.75) {};
			\node [style=none] (7) at (0.5, 0.25) {};
		\end{pgfonlayer}
		\begin{pgfonlayer}{edgelayer}
			\draw (3.center) to (0.center);
			\draw (1.center) to (4.center);
			\draw [in=-90, out=180, looseness=1.25] (2.center) to (5.center);
			\draw [bend right=90, looseness=2.00] (6.center) to (7.center);
			\draw (7.center) to (6.center);
		\end{pgfonlayer}
   \end{tikzpicture} }

\newcommand{\leftcoaction}[2]{
	\begin{tikzpicture}[scale=#1, yscale=-1]
		\begin{pgfonlayer}{nodelayer}
			\node [style=none] (0) at (0.5, 0.25) {};
			\node [style=none] (1) at (0.5, 0.75) {};
			\node [style=none] (2) at (0.25, 0.5) {};
			\node [style=none] (3) at (0.5, -0.5) {};
			\node [style=none] (4) at (0.5, 1.5) {};
			\node [style=none] (5) at (-0.5, 1.5) {};
			\node [style=none] (6) at (0.5, 0.75) {};
			\node [style=none] (7) at (0.5, 0.25) {};
		\end{pgfonlayer}
		\begin{pgfonlayer}{edgelayer}
			\draw (3.center) to (0.center);
			\draw (1.center) to (4.center);
			\draw [in=-90, out=180, looseness=1.25] (2.center) to (5.center);
			\draw [bend right=90, looseness=2.00] (6.center) to (7.center);
			\draw (7.center) to (6.center);
		\end{pgfonlayer}
	\end{tikzpicture} }

\newcommand{\rightcoaction}[2]{
	\begin{tikzpicture}[scale=#1, yscale=-1, xscale=-1]
		\begin{pgfonlayer}{nodelayer}
			\node [style=none] (0) at (0.5, 0.25) {};
			\node [style=none] (1) at (0.5, 0.75) {};
			\node [style=none] (2) at (0.25, 0.5) {};
			\node [style=none] (3) at (0.5, -0.5) {};
			\node [style=none] (4) at (0.5, 1.5) {};
			\node [style=none] (5) at (-0.5, 1.5) {};
			\node [style=none] (6) at (0.5, 0.75) {};
			\node [style=none] (7) at (0.5, 0.25) {};
		\end{pgfonlayer}
		\begin{pgfonlayer}{edgelayer}
			\draw (3.center) to (0.center);
			\draw (1.center) to (4.center);
			\draw [in=-90, out=180, looseness=1.25] (2.center) to (5.center);
			\draw [bend right=90, looseness=2.00] (6.center) to (7.center);
			\draw (7.center) to (6.center);
		\end{pgfonlayer}
	\end{tikzpicture}
}

\newcommand{\unitmap}[2]{
\begin{tikzpicture}[scale=#1]
	\begin{pgfonlayer}{nodelayer}
		\node [style=circle, scale=0.4, fill=#2] (0) at (0, 0) {};
		\node [style=none] (1) at (0, -0.4) {};
		\node [style=none] (4) at (0.13, 0) {};
		\node [style=none] (4) at (-0.13, 0) {};
	\end{pgfonlayer}
	\begin{pgfonlayer}{edgelayer}
		\draw [style=none] (0) to (1.center);
	\end{pgfonlayer}
\end{tikzpicture} }

\newcommand{\counitmap}[2]{
\begin{tikzpicture}[scale=#1, rotate=180]
	\begin{pgfonlayer}{nodelayer}
		\node [style=circle, scale=0.4, fill=#2] (0) at (0, 0) {};
		\node [style=none] (1) at (0, -0.4) {};
		\node [style=none] (4) at (0.13, 0) {};
	\end{pgfonlayer}
	\begin{pgfonlayer}{edgelayer}
		\draw [style=none] (0) to (1.center);
	\end{pgfonlayer}
\end{tikzpicture}
}

\newcommand{\comulmap}[2]{
\begin{tikzpicture}[scale=#1]
	\begin{pgfonlayer}{nodelayer}
		\node [style=circle,scale=0.4, fill=#2] (0) at (-1, 3) {};
		\node [style=none] (1) at (-1.25, 2.75) {};
		\node [style=none] (2) at (-0.75, 2.75) {};
		\node [style=none] (3) at (-1, 3.25) {};
		\node [style=none] (4) at (-0.68, 2.75) {};
	\end{pgfonlayer}
	\begin{pgfonlayer}{edgelayer}
		\draw [style=none] (3.center) to (0);
		\draw [style=none, bend right] (0.center) to (1.center);
		\draw [style=none, bend left] (0.center) to (2.center);
	\end{pgfonlayer}
\end{tikzpicture}}

\newcommand{\trianglemult}[1]{
\begin{tikzpicture}[scale=#1]
	\begin{pgfonlayer}{nodelayer}
		\node [style=none] (0) at (-0.25, 3.5) {};
		\node [style=none] (1) at (-0.5, 3.75) {};
		\node [style=none] (2) at (0, 3.75) {};
		\node [style=none] (3) at (-0.25, 3) {};
		\node [style=none] (4) at (0.25, 4.25) {};
		\node [style=none] (5) at (-0.75, 4.25) {};
	\end{pgfonlayer}
	\begin{pgfonlayer}{edgelayer}
		\draw [bend right, looseness=1.00] (2.center) to (4.center);
		\draw (0.center) to (1.center);
		\draw (0.center) to (2.center);
		\draw (2.center) to (1.center);
		\draw [in=-90, out=165, looseness=0.75] (1.center) to (5.center);
		\draw (0.center) to (3.center);
	\end{pgfonlayer}
\end{tikzpicture} }

\newcommand{\trianglecomult}[1]{
\begin{tikzpicture}[scale=#1]
	\begin{pgfonlayer}{nodelayer}
		\node [style=none] (0) at (-0.25, 3.75) {};
		\node [style=none] (1) at (-0.5, 3.5) {};
		\node [style=none] (2) at (0, 3.5) {};
		\node [style=none] (3) at (-0.25, 4.25) {};
		\node [style=none] (4) at (0.25, 3) {};
		\node [style=none] (5) at (-0.75, 3) {};
	\end{pgfonlayer}
	\begin{pgfonlayer}{edgelayer}
		\draw [bend left, looseness=1.00] (2.center) to (4.center);
		\draw (0.center) to (1.center);
		\draw (0.center) to (2.center);
		\draw (2.center) to (1.center);
		\draw [in=90, out=-165, looseness=0.75] (1.center) to (5.center);
		\draw (0.center) to (3.center);
	\end{pgfonlayer}
\end{tikzpicture}
}

\newcommand{\trianglecounit}[1]{
\begin{tikzpicture}[scale=#1]
	\begin{pgfonlayer}{nodelayer}
		\node [style=none] (0) at (-0.25, 3.5) {};
		\node [style=none] (1) at (-0.5, 3.25) {};
		\node [style=none] (2) at (0, 3.25) {};
		\node [style=none] (3) at (-0.25, 4.15) {};
		\node [style=none] (4) at (-0.25, 3) {};
	\end{pgfonlayer}
	\begin{pgfonlayer}{edgelayer}
		\draw (0.center) to (1.center);
		\draw (0.center) to (2.center);
		\draw (2.center) to (1.center);
		\draw (0.center) to (3.center);
	\end{pgfonlayer}
\end{tikzpicture} 
}

\newcommand{\triangleunit}[1]{
\begin{tikzpicture}[scale=#1]
	\begin{pgfonlayer}{nodelayer}
		\node [style=none] (0) at (-0.25, 4) {};
		\node [style=none] (1) at (-0.5, 4.25) {};
		\node [style=none] (2) at (0, 4.25) {};
		\node [style=none] (3) at (-0.25, 3.25) {};
	\end{pgfonlayer}
	\begin{pgfonlayer}{edgelayer}
		\draw (0.center) to (1.center);
		\draw (0.center) to (2.center);
		\draw (2.center) to (1.center);
		\draw (0.center) to (3.center);
	\end{pgfonlayer}
\end{tikzpicture}
}

\newcommand{\tricounit}[1]{\trianglecounit{#1}}
\newcommand{\triunit}[1]{\triangleunit{#1}}

\newcommand{\linduala}[1]{\xymatrixcolsep{6mm} \xymatrix{ \ar@{-||}[r]^{(#1)} & }}
\newcommand{\dagdual}{\xymatrixcolsep{4mm} \xymatrix{ \ar@{-||}[r]^{\dagger} & }}

\newcommand{\whitelin}{\xymatrixcolsep{3mm} \xymatrix{  \ar@{-||}[r]^{\circ} &  }}
\newcommand{\blacklin}{\xymatrixcolsep{3mm} \xymatrix{ \ar@{-||}[r]^{\bullet} & }}

\newcommand{\whitedag}[1]\dagmonw 
\newcommand{\blackdag}[1]\dagmonb 

\newcommand{\linmonw} {\xymatrixcolsep{4mm} \xymatrix{ \ar@{-||}[r]^{\circ} & }}
\newcommand{\linmonb} {\xymatrixcolsep{4mm} \xymatrix{  \ar@{-||}[r]^{\bullet} & }}
\newcommand{\dagmonw} {\xymatrixcolsep{4mm} \xymatrix{ \ar@{-||}[r]_{}^{\dagger\circ} & }}
\newcommand{\dagmonb} {\xymatrixcolsep{4mm} \xymatrix{  \ar@{-||}[r]_{}^{\dagger\bullet} & }}

\newcommand{\lincomonw} {\xymatrixcolsep{4mm} \xymatrix{ \ar@{-||}[r]_{\circ} & }}
\newcommand{\lincomonb} {\xymatrixcolsep{4mm} \xymatrix{  \ar@{-||}[r]_{\bullet} & }}
\newcommand{\dagcomonw} {\xymatrixcolsep{4mm} \xymatrix{ \ar@{-||}[r]_{\dagger\circ} & }}
\newcommand{\dagcomonb} {\xymatrixcolsep{4mm} \xymatrix{  \ar@{-||}[r]_{\dagger\bullet} & }}

\xymatrixrowsep{5mm}
\newdir{(>}{{}*!/-5pt/\dir{>}}  

\newcommand{\tri}[1]{
\begin{tikzpicture}[scale=#1]
	\begin{pgfonlayer}{nodelayer}
		\node [style=none] (0) at (2, 3.5) {};
		\node [style=none] (1) at (1.75, 3.25) {};
		\node [style=none] (2) at (2.25, 3.25) {};
	\end{pgfonlayer}
	\begin{pgfonlayer}{edgelayer}
		\draw (0.center) to (1.center);
		\draw (0.center) to (2.center);
		\draw (2.center) to (1.center);
	\end{pgfonlayer}
\end{tikzpicture} } 

\newcommand{\btri}[1]{
\begin{tikzpicture}[scale=#1]
	\begin{pgfonlayer}{nodelayer}
		\node [style=none] (0) at (2, 3.5) {};
		\node [style=none] (1) at (1.75, 3.25) {};
		\node [style=none] (2) at (2.25, 3.25) {};
	\end{pgfonlayer}
	\begin{pgfonlayer}{edgelayer}
		\draw [fill=black] (0.center) -- (1.center) -- (2.center) -- (0.center);
	\end{pgfonlayer}
\end{tikzpicture} } 

\newcommand{\linbialg}{\xymatrixcolsep{5mm} \xymatrix{ \ar@{-||}[r]^{\tri{0.75}} & }}
\newcommand{\linbialgblack}{\xymatrixcolsep{5mm} \xymatrix{ \ar@{-||}[r]^{\btri{0.75}} & }}
\newcommand{\dagbialg}{\xymatrixcolsep{6mm} \xymatrix{ \ar@{-||}[r]^{\dagger \tri{0.75}} & }}
\newcommand{\dagbialgb}{\xymatrixcolsep{6mm} \xymatrix{ \ar@{-||}[r]^{\dagger \btri{0.75}} & }}

\title{Exponential Modalities and Complementarity \\ (extended abstract)}
\author{Robin Cockett\thanks{Partially supported by NSERC, Canada} 
\qquad \qquad Priyaa Varshinee Srinivasan
\email{\quad \quad robin@ucalgary.ca\qquad \qquad priyaavarshinee.srin@ucalgary.ca}
\institute{Department of Computer Science\\
University of Calgary\\
Alberta, Canada}
}
\def\titlerunning{Exponential Modalities and Complementarity}
\def\authorrunning{Robin Cockett, Priyaa Varshinee Srinivasan}

\begin{document}
	
 \maketitle
	
\begin{abstract} 
The exponential modalities of linear logic have been used by various authors to model infinite-dimensional quantum systems.  
This paper explains how these modalities can also give rise to the complementarity principle of quantum mechanics.  

The paper uses a formulation of quantum systems based on $\dagger$-linear logic, 
whose categorical semantics lies in  mixed unitary categories, and a formulation of 
measurement therein.   The main result exhibits a complementary system as the result of 
measurements on free exponential modalities.  Recalling that, in linear logic, 
exponential modalities have two distinct but dual components, $!$ and $?$, 
this shows how these components under measurement become  ``compacted'' 
into the usual notion of a complementary Frobenius algebras from categorical 
quantum mechanics. 
\end{abstract}



 \section{Introduction}  

 Linear logic introduced by Girard in his seminal paper \cite{Gir87} treats logical statements 
 as resources, which cannot be duplicated or destroyed. The word ``linear" refers to this 
 resource sensitivity of the logic: a proof of a statement in linear logic may thus be regarded as a series of 
 resource transformations. In full linear logic the classical ability to duplicate and destroy resources is 
 recaptured by the exponential (or storage) modality written $!$ (pronounced the ``bang'').  The type $!A$ 
 may be interpreted as an unbounded ``store'' from which resources of type $A$ can be extracted an arbitrary (including 0)
 number of times. 
  
The exponential modality $!$ has been proposed as a structure for modelling infinite dimensional systems: \cite{Vic08} used the 
exponential modality to model the quantum harmonic oscillator, and \cite{BPS94} used it to model the bosonic Fock space. 
However, these uses did not explain what exponential modalities  $!$ and its dual $?$ (pronounced the ``whimper'')
have to do with the complementarity principle of quantum mechanics  \cite{CoD11}.  A pair of 
quantum observables (physical properties of a system) is said to be complementary if 
measuring one observable increases the uncertainty regarding the value of the other. 
The classic example is that the more one knows about position of a particle the less one knows about its momentum.
The purpose of this article is to exhibit a relationship between the exponential modalities $!$ and its dual $?$,
and complementary observables --- a relationship which suggests a possibly new perspective on measurement in quantum systems.

Linearly distributive categories (LDCs) \cite{CS97} provide a categorical semantics for the multiplicative fragment of linear logic (MLL). 
Thus, LDCs are equipped with two distinct tensors called the ``tensor'', $\ox$, and the ``par'', $\oplus$.\footnote{In the linear 
logic community the par is often denoted by $\parr$ but we follow the convention in \cite{CS97} and use $\oplus$.}  
These are related by a linear distributor.   It is not assumed that the tensor is dual to the par  --- which would be normal in linear logic.  
In an LDC, having a dual is a property which an object may or may not possess.  When every object possesses a dual then the category is $*$-autonomous.

In this development, LDCs which satisfy the so-called ``mix law" are particularly important.  The mix law provides a natural transformation from 
the tensor to the par called the {\em mixor}.  When the mixor is a natural isomorphism the LDC becomes equivalent to a monoidal category.  
Conversely, monoidal categories can be viewed as being degenerate or {\em compact} LDCs in which the mixor is the identity map.  
Thus, from this perspective, a compact closed category is a compact LDC with duals, that is, a compact $*$-autonomous category.  
 
In \cite{CCS18}, we introduced $\dagger$-linearly distributive categories ($\dagger$-LDCs) with mix for modelling possibly infinite-dimensional 
quantum processes.  In a mix LDC, there is always a set of objects which cannot distinguish between the tensor and the par in the sense 
that $A \ox \_ \simeq A \oa \_$: these objects form a compact subLDC called the {\em core\/}.  In a mix $\dagger$-LDC, it is possible to go one 
step further, and identify a {\em unitary core} in which every object is not only in the core but isomorphic to its $\dagger$-dual in a coherent way.  
A unitary core is equivalent to a $\dagger$-monoidal category and when this category has duals it is equivalently a $\dagger$-compact closed 
category \cite{AbC04,Sel07}. This is the main structure underlying categorical quantum mechanics (CQM) \cite{CoA17,HeV19}: 
finite-dimensional Hilbert spaces provide the paradigmatic example.

The general notion of a {\em mixed unitary category\/} (MUC) is essentially a mix $\dagger$-LDC with a specified unitary core.  In particular, the unitary core 
may be viewed as comprising the finite-dimensional processes while the larger category extends this to include infinite-dimensional processes.  
An example of a MUC is given by the embedding of complex finite matrices into the category of finiteness matrices \cite{Ehr05}.  Another example is given by 
the embedding of the finite-dimensional Hilbert spaces within the category of 
Chu spaces \cite{Bar06} of vector spaces over complex numbers 
with the field $\C$ as the dualizing object.  For further details see \cite{CCS18,CS19}, where  completely positive maps, and environment structures 
for MUCs are described.  In this article, we explore the notions of measurement and complementarity in MUCs.  

In CQM, Coecke and Pavlovic \cite{CoP07} described a ``demolition'' measurement in a 
$\dagger$-monoidal category as a map, $m: A \to X$, with 
$m^\dagger m = 1_X$, to a special commutative $\dagger$-Frobenius algebra, $X$.   Interpreted 
in the category of finite-dimensional Hilbert spaces, the notion of the demolition measurement 
models the Projection-Valued Measures (PVMs) of quantum mechanics. 
Generalizing this idea to MUCs to model measurements here is complicated by the fact that, in a MUC,
generally, $A \neq A^\dagger$, except in the unitary core.  Thus, in a MUC, a measurement can be viewed as a two-step process 
in which one first ``compacts'' an object into the unitary core by a retraction and then one performs Coecke and Pavlovic's 
demolition measurement.  The compaction process is discussed in Section \ref{Sec: compaction} and is already quite interesting:  it gives 
rise to a $\dagger$-binary idempotent. Conversely, a $\dagger$-binary idempotent, which is ``coring'' and splits, gives rise to a compaction 
into the ``canonical'' unitary core.  

In CQM, quantum observables are characterized by certain $\dagger$-Frobenius algebras \cite{CPV13} in 
$\dagger$-monoidal categories. Two such $\dagger$-Frobenius algebras, $(A, \mulmap{1.5}{white}, \unitmap{1.5}{white}, 
\comulmap{1.5}{white}, \counitmap{1.5}{white})$ and  $(A, \mulmap{1.5}{black}, \unitmap{1.5}{black}, 
\comulmap{1.5}{black}, \counitmap{1.5}{black})$ are said to be complementary \cite{CoD11} if  
$(A, \mulmap{1.5}{white}, \unitmap{1.5}{white}, \comulmap{1.5}{black}, \counitmap{1.5}{black})$
and $(A, \mulmap{1.5}{black}, \unitmap{1.5}{black}, \comulmap{1.5}{white}, \counitmap{1.5}{white})$ 
are Hopf algebras. An object which is a Frobenius algebra is always self-dual. 
In an LDC, a linear monoid, 
$A \linmonw B$, is a $\ox$-monoid $A$ together with a dual $B$.  Because $B$ is dual to $A$ ---
and $A$ is a $\ox$-monoid --- it follows that $B$ is a $\oa$-comonoid. 
In contrast a linear {\em comonoid}\footnote{Note that this is {\em not\/} the dual notion of a linear monoid as a linear monoid is a self-dual notion in an LDC.}, 
$A \lincomonw B$, is a $\ox$-comonoid $A$ together with a dual $B$: this means that $B$ is a $\oa$-monoid. 
A linear monoid and a linear comonoid interact to produce a linear bialgebra: this has a $\ox$-bialgebra  on $A$ and a $\oa$-bialgebra on $B$. 
In a MUC, the linear bialgebras in the unitary core are the base for defining complimentary systems. 
These structures are presented in Section \ref{Sec: complementarity}.

Section \ref{Sec: exponential modalities} describes the connection between the free exponential modalities 
and complimentary systems in a $\dagger$-isomix setting. An LDC is said to have exponential modalities, if 
it has a monoidal comonad $(!, \delta, \epsilon)$, a comonoidal monad $(?, \mu, \eta)$, and for all objects $A$ 
$(!A, \Delta_A, \tricounit{0.65}_A)$ is a natural commutative $\ox$-monoid and $(?A, \nabla_A, \triunit{0.65}_A)$ 
is a natural commutative $\oa$-monoid. The modalities are said to be free if $(!A, \Delta_A, \tricounit{0.65}_A)$ is 
cofree and $(?A, \nabla_A, \triunit{0.65}_A)$ is free. The main result of this paper is that in a MUC,
every $\dagger$-complementary system in the unitary core arises as the splitting of a $\dagger$-binary 
idempotent on the $\dagger$-linear bialgebra induced on the free exponentials.  This is an interesting result 
since it shows that complementary observables arise from compacting dual but distinct systems of arbitrary dimensions. 

\medskip

\noindent
 {\bf Notation:} Diagrammatic order of composition is used: so $fg$ should be read as $f$ followed by $g$.  Circuit diagrams should be 
be read top to bottom: that is following the direction of gravity! 

\vspace{0.5em}

\noindent
A full version of this article containing all proofs  is available in arXiv \cite{CoS21}. 



\section{Preliminaries}
\label{Sec: preliminaries}

In this section, we recall the definitions of dagger isomix categories, unitary categories, and 
mixed unitary categories from \cite{CCS18}. To achieve this we start by recalling the definitions of linearly distributive 
categories and isomix categories.

A linearly distributive category (LDC) \cite{CS97}, $(\X, \ox, \oa)$, is a category with two tensor products --- $\ox$ called the tensor with 
unit $\top$, and the $\oa$ called the par with unit $\bot$. The tensor and the par interact by means of  
linear distributors which are natural transformations (which, in general, are not isomorphisms): 
\[ \partial^L: A \ox (B \oa C) \rightarrow  (A \ox B) \oa C ~~~~~~~~~~~~  
\partial^R: (B \oa C) \ox A \rightarrow B \oa (C \ox A) \]

A symmetric LDC is an LDC in which both monoidal structures are symmetric,  with symmetry maps $c_\ox: A \ox B \to B \ox A$ and $c_\oa: A \oa B \to B \oa A$, such that
$\partial^R = c_\ox (1 \ox c_\oa) \partial^L (c_\ox \oa 1) c_\oa$.  
LDCs provide a categorical semantics for linear logic, and are equipped with a graphical calculus; see \cite[Section 2]{CCS18} 
and \cite{BCST96}.

A mix category is an LDC with a mix map, $\m: \bot \to \top$; when $\m$ is an isomorphism it is an isomix category.  The mix map gives a natural mixor map, 
$\mx: A \ox B \to A \oa B$, which, even if the mix map is an isomorphism, is usually not an isomorphism.  An isomix category in which every mixor map 
is an isomorphism is a compact LDC. A compact LDC with $\m = 1$ and $\mx = 1$ is just a monoidal category.

The {\bf core}, $\Core(\X)$, of an isomix category $\X$ is the full subcategory given by the objects, $U$, 
such that for all $A \in \X$, the maps $\mx_{U,A}: U \ox A \to U \oa A$ and $\mx_{A,U}: A \ox U \to A \oa U$ are isomorphisms.
The units, $\top$ and $\bot$, are always in the core.  The core $\Core(\X)$ of an isomix category $\X$ is always a compact LDC.

A {\bf $\dagger$-linearly distributive category} \cite{CCS18} is an LDC $\X$ with a functor $(\_)^\dag:
\X^\op\to \X$ and the following natural isomorphisms satisfying the coherence conditions which are 
described in \cite{CCS18}. 
\begin{align*}
&\text{ \bf tensor laxors: }  A^\dag \ox B^\dag \to^{\lambda_\ox} (A\oa B)^\dag ~~~~~~~~~~
 A^\dag \oa B^\dag \to^{ \lambda_\oa} (A\ox B)^\dag \\
&\text{ \bf unit laxors: } \top \to^{\lambda_\top} \bot^\dag ~~~~~~~~~~~ 
\bot \to^{\lambda_\bot} \top^\dag \\ 
&\text{ \bf involutor: }  A \xrightarrow{\iota} (A^\dag)^\dag 
\end{align*}

In a $\dagger$-LDC, it is generally the case that $A \neq A^\dagger$ because $\dagger$ swaps the tensor and the par. A {\bf $\dagger$-mix category} is a $\dagger$-LDC which has a mix map which  satisfies, in addition the following commuting diagram:

{\centering $\mbox{\bf [$\dagger$-\text{mix}]}  ~~~~\begin{array}[c]{c} 
\xymatrix{
\bot                 \ar@{->}[r]^{{\sf m}} \ar@{->}[d]_{\lambda_\bot}  
  & \top             \ar@{->}[d]^{\lambda_\top}\\
\top^\dag            \ar@{->}[r]_{{\sf m}^\dag}
  & \bot^\dag
} 
\end{array}$ \par}

If ${\sf m}$ is an isomorphism, then $\X$ is an {\bf $\dagger$-isomix category}.  A {\bf compact $\dagger$-LDC} 
is a compact LDC which is also a $\dagger$-isomix category. 

In a $\dagger$-monoidal category a unitary isomorphism is an isomorphism $f$ with $f^\dagger = f$.  
In a $\dagger$-LDC, an object, $A$, does not necessarily coincide with its dagger, $A^\dagger$: this 
means that describing unitary isomorphism for $\dagger$-LDCs is more complicated.  To accomplish this 
the notion of {\bf unitary structure} which is described in \cite{CCS18} is used.  
A {\bf pre-unitary} object in a $\dagger$-isomix category is an object $A$ in the core with an isomorphism 
$\varphi: A \to A^\dagger$ such that $\varphi (\varphi^{-1})^\dagger = \iota$ (where $\iota$ is the involutor).   Unitary structure for a 
$\dagger$-isomix category is given by a family of pre-unitary objects satisfying certain closure and 
coherences requirements.  

A {\bf unitary category} is a compact $\dagger$-LDC equipped with unitary structure which 
makes every object a (pre)unitary object.  A $\dagger$-monoidal category \cite{HeV19} 
is a unitary category with $\ox = \oa$ and the unitary structure given by the 
identity map.  Conversely, every unitary category is $\dagger$-linearly equivalent to a 
$\dagger$-monoidal category via the $\dagger$-linear functor ${\sf Mx_\downarrow}: (\X, \ox, \oa) \to (\X, \ox, \ox)$, 
see \cite[Prop. 5.11]{CCS18}.    

A {\bf mixed unitary category} (MUC) \cite{CCS18} is a $\dagger$-isomix category, $\C$, 
equipped with a strong $\dagger$-isomix functor $M: \U \to \C$ from a unitary category 
$\U$ such that there are natural transformations:

\[ \mx': M(U) \ox X \to M(U) \oa X \text{ with } \mx~ \mx' = 1, ~ \mx' \mx = 1 \]

Thus, a mixed unitary category can be visualized schematically as:
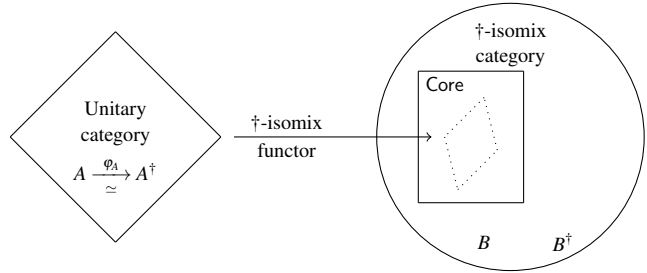
\[  \begin{tikzpicture} [scale=1.4]
	\begin{pgfonlayer}{nodelayer}
		\node [style=circle, scale=14] (0) at (4.5, 0) {};
		\node [style=none] (1) at (-3, 2) {};
		\node [style=none] (2) at (-5, 0) {};
		\node [style=none] (3) at (-3, -2) {};
		\node [style=none] (4) at (-1, 0) {};
		\node [style=none] (5) at (-3, -0.75) {$A \to^{\varphi_A}_{\simeq} A^\dagger$};
		\node [style=none] (6) at (-3, 0.5) {Unitary};
		\node [style=none] (7) at (-3, 0) {category};
		\node [style=none] (8) at (-0.75, 0) {};
		\node [style=none] (9) at (3, 0) {};
		\node [style=none] (10) at (0.25, 0.25) {$\dagger$-isomix};
		\node [style=none] (11) at (0.25, -0.25) {functor};
		\node [style=none] (12) at (4.5, 2) {$\dagger$-isomix};
		\node [style=none] (13) at (4.5, 1.5) {category};
		\node [style=none] (14) at (4, -2) {$B$};
		\node [style=none] (15) at (5.5, -2) {$B^\dagger$};
		\node [style=none] (16) at (4, 0.75) {};
		\node [style=none] (17) at (3.25, 0) {};
		\node [style=none] (18) at (3.5, -1) {};
		\node [style=none] (19) at (4.25, -0.25) {};
		\node [style=none] (20) at (2.75, 1.25) {};
		\node [style=none] (21) at (4.75, 1.25) {};
		\node [style=none] (22) at (4.75, -1.25) {};
		\node [style=none] (23) at (2.75, -1.25) {};
		\node [style=none] (24) at (3.25, 1) {$\Core$};
	\end{pgfonlayer}
	\begin{pgfonlayer}{edgelayer}
		\draw (1.center) to (2.center);
		\draw (2.center) to (3.center);
		\draw (3.center) to (4.center);
		\draw (4.center) to (1.center);
		\draw [->] (8.center) to (9.center);
		\draw [dotted] (16.center) to (17.center);
		\draw [dotted] (17.center) to (18.center);
		\draw [dotted] (18.center) to (19.center);
		\draw [dotted] (19.center) to (16.center);
		\draw (20.center) to (21.center);
		\draw (21.center) to (22.center);
		\draw (22.center) to (23.center);
		\draw (23.center) to (20.center);
	\end{pgfonlayer}
\end{tikzpicture}
  \]
  Within the unitary category, $A \simeq A^\dagger$ by the means of the unitary structure map. 
However, outside the unitary core,  an object is not in general isomorphic to its dagger. 

Given any $\dagger$-isomix category $\X$, the preunitary objects always 
form a unitary category, $\Unitary(\X)$ with a forgetful  $\dagger$-isomix functor $U: \Unitary(\X) \to \X$
which produces a MUC.  $\Unitary(\X)$ satisfies a couniversal property, see \cite[Section 5.2]{CCS18}, 
and is the ``largest'' possible unitary core for the $\dagger$-isomix category $\X$.  We shall call 
$\Unitary(\X)$ the {\bf canonical} unitary core of $\X$.



\section{Measurement}
\label{Sec: compaction}

A measurement in a MUC can be broken into two steps: a compaction step into an object in the unitary core followed by 
a demolition measurement within the unitary core.

\begin{definition}
Let $M: \U \to \C$ be a MUC. A {\bf compaction } of an object $A \in \C$ to $U \in \U$ is a retraction, 
$r: A \to M(U)$.  This means that  there is a section $s: M(U) \to A$ such that $sr = 1_{M(U)}$.
A compaction is said to be {\bf canonical} when $\U = \Unitary(\X)$ (so $U$ is a preunitary object).
\end{definition}

The compact object, $M(U)$, has a unitary structure map which is an isomorphism between $M(U)$ 
and $M(U)^\dagger$ given by composing the unitary structure map of $U$ with the preservator $\rho$ (see \cite[Definition 3.17]{CCS18} for the complete definiton of a preservator): 
\[ \psi := M(U) \to^{M(\varphi)} M(U^\dagger) \to^{\rho} M(U)^\dagger \] 
Once one has reached $M(U)$ by a compaction, one can follow with a classical demolition measurement $U \to^w X$ to obtain an overall 
compaction $A \to^{r M(w)} M(X)$, which gives a (demolition) measurement in a MUC. 

We start by showing how a compaction gives rise to a binary idempotent:

\begin{definition}
A {\bf binary idempotent} in any category is a pair of maps $(\u,\v)$ with $\u: A \to B$, and $\v: B \to A$ such that $\u\v\u = \u$, and $\v \u \v = \v$.
\end{definition}

A binary idempotent, $(\u, \v): A \to B$ gives a pair of idempotents: $e_A := \u \v : A \to A$, and $e_B := \v \u : B \to B$. 
We say the binary idempotent $(\u, \v)$ {\bf splits} in case the idempotents $e_A$ and $e_B$ split.

\begin{lemma}
	\label{Lemma: binary idempotent equivalence}
	In any category the following are equivalent: 
\begin{enumerate}[(i)]
\item $(\u, \v) : A \to B$ is a binary idempotent which splits.
\item  $e: A \to A$, and $d:B \to B$ are a pair of idempotents which split through isomorphic objects.
\end{enumerate}
\end{lemma}

Observe that a compaction of an object, say $A$, in any MUC, gives the following system of maps:

{ \centering 
$\xymatrixcolsep{12mm}
 \xymatrix{
    A  \ar@<3pt>[r]^r &  
    M(U) \ar@<3pt>[l]^s \ar@<3pt>[r]^{\psi := M(\varphi)\rho} & 
    M(U)^\dagger \ar@<3pt>[r]^{r^\dagger} \ar@<3pt>[l]^{\psi^{-1}} &
    A^\dagger \ar@<3pt>[l]^{s^\dagger}    } $
\par }

Thus the compaction gives rise to a binary idempotent  $(\u,\v): A \to A^\dagger$ 
where $\u:= r\psi r^\dagger$ and $\v:= s^\dagger \psi^{-1} s$.

Because $U$ is a unitary object, we have that $\varphi (\varphi^{-1 \dagger}) = \iota$.   
The preservator, on the other hand, satisfies $\iota \rho^\dagger = M(\iota) \rho$  
(see after Definition 3.17 in \cite{CCS18}).  Thus, $\iota \rho^\dagger = M(\iota) \rho = 
M(\varphi\varphi^{-1 \dagger})\rho = M(\varphi)\rho M(\varphi^{-1})^\dagger$ and hence 
$\psi = M(\varphi) \rho= \iota \rho^\dagger M(\varphi)^\dagger = 
\iota(M(\varphi)\rho)^\dagger = \iota \psi^\dagger$.  
This allows us to observe:
\begin{align*} 
\iota\u^\dagger & = \iota (r \psi r^\dagger)^\dagger = \iota r^{\dagger\dagger} \psi^\dagger r^\dagger = r \iota \psi^\dagger r^\dagger = r \psi r^\dagger = \u \\
\v^\dagger & = (s^\dagger \psi^{-1} s)^\dagger = s^\dagger (\psi^{\dagger})^{-1} s^{\dagger\dagger} 
= s^\dagger (\iota^{-1}\psi)^{-1} s^{\dagger\dagger} 
 = s^\dagger \psi^{-1} \iota s^{\dagger\dagger} = s^\dagger \psi^{-1} s \iota  = \v \iota
\end{align*}

This leads to  the following definition:

\begin{definition}
A binary idempotent, $(\u, \v): A \to A^\dagger$ in a $\dagger$-LDC is a {\bf $\dagger$-binary idempotent}, 
written $\dagger(\u, \v)$, if $\u  = \iota \u^\dagger$
and $\v^\dagger = \v \iota$.
\end{definition}

In a $\dagger$-monoidal category, where $A=A^\dagger$ and $\iota=1_A$ 
this makes $\u=\u^\dagger$ and $\v=\v^\dagger$; thus  $\u\v= (\v\u)^\dagger$.  This means 
that if we require $\u\v = \v\u$ we obtain a dagger idempotent in the sense of \cite{Sel08}.  

Splitting a $\dagger$-binary idempotent almost produces a preunitary object.    In a $\dagger$-LDC, we shall call an object $A$
 with an isomorphism $\varphi: A \to A^\dagger$ such that $\varphi \varphi^{\dagger -1} = \iota$ a {\bf weak preunitary object}.  
 Clearly, in a $\dagger$-isomix category, a weak preunitary object $(A,\varphi)$ is a preunitary object when, in addition, 
 $A$ is in the core.  We next observe that dagger binary idempotent splits through weak preunitary objects:

\begin{lemma} In a $\dagger$-LDC with a $\dagger$-binary idempotent $\dagger(\u,v): A \to A^\dagger$:
\label{Lemma: dagger splitting}
\begin{enumerate}[(i)]
\item $e_{A^\dagger} := \v \u = (\u \v)^\dagger =: (e_A)^\dagger$;
\item if $\dagger(\u,\v)$ splits with $e_A = A \to^r E \to^s A$ then $E$ is a weak preunitary object.
\end{enumerate}
\end{lemma}

Thus, in a $\dagger$-isomix category, an object which splits a $\dagger$-binary idempotent is always weakly preunitary.  
In order to ensure that the splitting of a $\dagger$-binary idempotent is a preunitary object --- and so a canonical compaction --- it 
remains to ensure that the splitting is in the core. This leads to the following definition:

\begin{definition}
	An idempotent $A \to^{e} A$ in an isomix category, $\X$, is a {\bf coring idempotent} if it is equipped with natural
	$\kappa^L_X: X \oa A \to X \ox A$ and $\kappa^R_X: A \oa X \to A \ox X$ such that the following diagrams 
	commute:
    \[ \xymatrix{
    X \ox A \ar[r]^{1 \ox e}  \ar[d]_{1 \ox e} \ar@{}[dr]|{\bf [KL.1]}  & X \ox A \ar[d]^{\mx}\\
    X \ox A  & X \oa A \ar[l]^{\kappa^L_X}
    } ~~~~~~~~ \xymatrix{
    X \oa A \ar[r]^{1 \oa e}  \ar[d]_{1 \oa e} \ar@{}[dr]|{\bf [KL.2]}  & X \oa A \ar[d]^{\kappa^L_X}\\
    X \oa A  & X \ox A \ar[l]^{\mx}
    }
    ~~~~~~~~ \xymatrix{
    A \ox X \ar[r]^{e \ox 1}  \ar[d]_{e \ox 1} \ar@{}[dr]|{\bf [KR.1]}  & A \ox X \ar[d]^{\mx}\\
    A \ox X  & A \oa X \ar[l]^{\kappa^R_X}
    } ~~~~~~~~ \xymatrix{
    A \oa X \ar[r]^{e \oa 1}  \ar[d]_{e \oa 1} \ar@{}[dr]|{\bf [KR.2]}  & A \oa X \ar[d]^{\kappa^R_X}\\
    A \oa X  & A \ox X \ar[l]^{\mx}
    } \]
\end{definition}

For a coring idempotent $A \to^{e} A$, the transformations $\kappa^{\_}_X$ act on a splitting as the inverse of the mixor, $\mx$.  Thus, a coring idempotent  
always splits through the core:

\begin{lemma} In a mix category: 
\label{Lemma: pseudocore}
\begin{enumerate}[(i)]
\item An idempotent splits through the core if and only if it is coring;
\item If $(u,v)$ is a binary idempotent then $\u\v$ is coring if and only if $\v\u$ is coring.
\end{enumerate}
\end{lemma}

This allows:

\begin{definition}
A {\bf coring binary idempotent} in a mix category is a binary idempotent, $(\u,\v)$, for which either $\u\v$ or $\v \u$ 
is a coring idempotent.
\end{definition}

These observations can be summarized by the following:

\begin{theorem}
In the MUC $M\!:\!\Unitary(\X) \to \X$, with $\dagger$-isomix category $\X$, an object $U$ is a compaction of $A$ if and only if $U$ is the
splitting of a coring $\dagger$-binary idempotent $\dagger(\u,\v): A \to A^\dagger$.
\end{theorem}

Using this characterization of canonical compaction, we will show that, in the presence of free $\dagger$-exponential 
modalities, complementarity always arises as a canonical compaction of a $\dagger$-linear bialgebra on 
the free exponential modalities. 


\section{Complementarity}
\label{Sec: complementarity}

The objective of this section, is to describe strong complementarity within a $\dagger$-isomix category.  Strong complementarity 
classically is, in a $\dagger$-monoidal setting,  between two special commutative $\dagger$-Frobenius algebras. 
In a linear setting with two distinct tensor products, Frobenius Algebras are generalized by linear monoids \cite{Egg10, CKS00} 
which consist of a $\ox$-monoid and a dual $\oa$-comonoid.  The directionality of the linear distributor makes a bialgebraic interaction between two 
$\dagger$-linear monoids impossible.  However, such an interaction is possible between a 
$\dagger$-linear monoid and a $\dagger$-linear comonoid, and this gives a $\dagger$-linear bialgebra. These  $\dagger$-linear bialgebras provide the basis for complementarity in a $\dagger$-isomix category.   
In a MUC, one can, furthermore, consider the effect of a compaction which preserves these structures to arrive back at the classical CQM notion 
of a complementary system.

\subsection{Duals}

\begin{definition} A {\bf dual} in an LDC, $(\eta, \epsilon): A \dashvv B$, consists of maps, 
$\eta: \top \to A \oa B$, and $\epsilon: B \ox A \to \bot$ such that the 
snake diagrams hold.  A {\bf morphism of duals}, $(f,g): (\eta,  \epsilon): A \dashvv B \to (\tau, \gamma): A' \dashvv B'$, is given by a pair of maps
$f: A \to A'$ and $g: B' \to B$ such that:  
\begin{center}
(a) ~~~  \begin{tikzpicture}
	\begin{pgfonlayer}{nodelayer}
		\node [style=none] (0) at (-1, 2) {};
		\node [style=none] (1) at (0.5, 3) {};
		\node [style=none] (2) at (-1, 3) {};
		\node [style=none] (3) at (-0.25, 4) {$\tau$};
		\node [style=none] (4) at (0.5, 2) {};
		\node [style=none] (5) at (-1.25, 2.25) {$A'$};
		\node [style=none] (6) at (1, 3.5) {$B'$};
		\node [style=circle, scale=1.5] (7) at (0.5, 2.75) {};
		\node [style=none] (8) at (0.5, 2.75) {$g$};
		\node [style=none] (9) at (1, 2.25) {$B$};
	\end{pgfonlayer}
	\begin{pgfonlayer}{edgelayer}
		\draw (4.center) to (7);
		\draw (1.center) to (7);
		\draw (2.center) to (0.center);
		\draw [bend left=90, looseness=1.75] (2.center) to (1.center);
	\end{pgfonlayer}
\end{tikzpicture}
= \begin{tikzpicture}
	\begin{pgfonlayer}{nodelayer}
		\node [style=none] (0) at (0.5, 2) {};
		\node [style=none] (1) at (-1, 3) {};
		\node [style=none] (2) at (0.5, 3) {};
		\node [style=none] (3) at (-0.25, 4) {$\eta$};
		\node [style=none] (4) at (-1, 2) {};
		\node [style=none] (5) at (0.75, 2.25) {$B$};
		\node [style=none] (6) at (-1.5, 3.5) {$A$};
		\node [style=circle, scale=1.5] (7) at (-1, 2.75) {};
		\node [style=none] (8) at (-1, 2.75) {$f$};
		\node [style=none] (9) at (-1.5, 2.25) {$A'$};
	\end{pgfonlayer}
	\begin{pgfonlayer}{edgelayer}
		\draw (4.center) to (7);
		\draw (1.center) to (7);
		\draw (2.center) to (0.center);
		\draw [bend right=90, looseness=1.75] (2.center) to (1.center);
	\end{pgfonlayer}
\end{tikzpicture}
  ~~~~~~(b)~~~ \begin{tikzpicture}
	\begin{pgfonlayer}{nodelayer}
		\node [style=none] (0) at (0.5, 3.75) {};
		\node [style=none] (1) at (-1, 2.75) {};
		\node [style=none] (2) at (0.5, 2.75) {};
		\node [style=none] (3) at (-0.25, 1.75) {$\gamma$};
		\node [style=none] (4) at (-1, 3.75) {};
		\node [style=none] (5) at (0.75, 3.5) {$B'$};
		\node [style=none] (6) at (-1.5, 2.25) {$A'$};
		\node [style=circle, scale=1.5] (7) at (-1, 3) {};
		\node [style=none] (8) at (-1, 3) {$f$};
		\node [style=none] (9) at (-1.5, 3.5) {$A$};
	\end{pgfonlayer}
	\begin{pgfonlayer}{edgelayer}
		\draw (4.center) to (7);
		\draw (1.center) to (7);
		\draw (2.center) to (0.center);
		\draw [bend right=90, looseness=1.75] (1.center) to (2.center);
	\end{pgfonlayer}
\end{tikzpicture}
 =  \begin{tikzpicture}
	\begin{pgfonlayer}{nodelayer}
		\node [style=none] (0) at (-1, 3.75) {};
		\node [style=none] (1) at (0.5, 2.75) {};
		\node [style=none] (2) at (-1, 2.75) {};
		\node [style=none] (3) at (-0.25, 1.75) {$\epsilon$};
		\node [style=none] (4) at (0.5, 3.75) {};
		\node [style=none] (5) at (-1.25, 3.5) {$A$};
		\node [style=none] (6) at (1, 2.25) {$B$};
		\node [style=circle, scale=1.5] (7) at (0.5, 3) {};
		\node [style=none] (8) at (0.5, 3) {$g$};
		\node [style=none] (9) at (1, 3.5) {$B'$};
	\end{pgfonlayer}
	\begin{pgfonlayer}{edgelayer}
		\draw (4.center) to (7);
		\draw (1.center) to (7);
		\draw (2.center) to (0.center);
		\draw [bend right=90, looseness=1.75] (2.center) to (1.center);
	\end{pgfonlayer}
\end{tikzpicture}
\end{center}
A {\bf self-duality} is a dual $(\eta, \epsilon): A \dashvv B$ in which $A$ is isomorphic to $B$ (or indeed $A=B$). 
\end{definition}

A morphism $(f,g)$ of duals is determined by either of the maps, as $f$ is 
dual to $g$: they are Australian mates; see \cite{CKS00}. 
In a $\dagger$-LDC, if $A$ is dual to $B$, then $B^\dagger$ is dual to $A^\dagger$:

A binary idempotent can implicitly express a morphism of duals, which becomes explicit when the idempotent splits.

\begin{definition}
A binary idempotent $(\u,\v)$ is {\bf retractional} on a dual $(\eta, \epsilon): A \dashvv B$ if equations $(a)$ and $(b)$, below, hold.
On the other hand $(\u,\v)$ is {\bf sectional}, if equations $(c)$ and $(d)$ hold: 
\begin{equation*}
 (a)~
	\begin{tikzpicture}
		\begin{pgfonlayer}{nodelayer}
			\node [style=none] (0) at (-5, 0.25) {};
			\node [style=none] (1) at (-5, 1.5) {};
			\node [style=none] (2) at (-3.5, 1.5) {};
			\node [style=none] (3) at (-3.5, 0.25) {};
			\node [style=circle, scale=1.8] (4) at (-5, 1) {};
			\node [style=none] (5) at (-5, 1) {$e_A$};
			\node [style=none] (6) at (-5.25, 0.5) {$A$};
			\node [style=none] (7) at (-3.25, 0.5) {$B$};
		\end{pgfonlayer}
		\begin{pgfonlayer}{edgelayer}
			\draw (0.center) to (4);
			\draw (1.center) to (4);
			\draw [bend left=90, looseness=1.75] (1.center) to (2.center);
			\draw (3.center) to (2.center);
		\end{pgfonlayer}
	\end{tikzpicture} =\begin{tikzpicture}
		\begin{pgfonlayer}{nodelayer}
			\node [style=none] (0) at (-5, 0.25) {};
			\node [style=none] (1) at (-5, 1.5) {};
			\node [style=none] (2) at (-3.5, 1.5) {};
			\node [style=none] (3) at (-3.5, 0.25) {};
			\node [style=circle, scale=1.8] (4) at (-5, 1) {};
			\node [style=none] (5) at (-5, 1) {$e_A$};
			\node [style=none] (6) at (-5.25, 0.5) {$A$};
			\node [style=none] (7) at (-3.25, 0.5) {$B$};
			\node [style=circle, scale=1.8] (8) at (-3.5, 1) {};
			\node [style=none] (9) at (-3.5, 1) {$e_B$};
		\end{pgfonlayer}
		\begin{pgfonlayer}{edgelayer}
			\draw (0.center) to (4);
			\draw (1.center) to (4);
			\draw [bend left=90, looseness=1.75] (1.center) to (2.center);
			\draw (3.center) to (8);
			\draw (8) to (2.center);
		\end{pgfonlayer}
	\end{tikzpicture}
~~~~
(b) ~ \begin{tikzpicture}
	\begin{pgfonlayer}{nodelayer}
		\node [style=none] (0) at (0.75, 1.5) {};
		\node [style=none] (1) at (0.75, 0.25) {};
		\node [style=none] (2) at (2.25, 0.25) {};
		\node [style=none] (3) at (2.25, 1.5) {};
		\node [style=circle, scale=1.8] (4) at (0.75, 0.75) {};
		\node [style=none] (5) at (0.75, 0.75) {$e_B$};
		\node [style=none] (6) at (2.5, 1.25) {$A$};
		\node [style=none] (7) at (0.5, 1.25) {$B$};
		\node [style=none] (8) at (1.5, 2) {};
	\end{pgfonlayer}
	\begin{pgfonlayer}{edgelayer}
		\draw (0.center) to (4);
		\draw (1.center) to (4);
		\draw [bend right=90, looseness=1.75] (1.center) to (2.center);
		\draw (3.center) to (2.center);
	\end{pgfonlayer}
\end{tikzpicture}
 =\begin{tikzpicture}
	\begin{pgfonlayer}{nodelayer}
		\node [style=none] (0) at (-3.25, 1.5) {};
		\node [style=none] (1) at (-3.25, 0.25) {};
		\node [style=none] (2) at (-4.75, 0.25) {};
		\node [style=none] (3) at (-4.75, 1.5) {};
		\node [style=circle, scale=1.8] (4) at (-3.25, 0.75) {};
		\node [style=none] (5) at (-3.25, 0.75) {$e_A$};
		\node [style=none] (6) at (-3, 1.25) {$A$};
		\node [style=none] (7) at (-5, 1.25) {$B$};
		\node [style=circle, scale=1.8] (8) at (-4.75, 0.75) {};
		\node [style=none] (9) at (-4.75, 0.75) {$e_B$};
		\node [style=none] (10) at (-4, 2) {};
	\end{pgfonlayer}
	\begin{pgfonlayer}{edgelayer}
		\draw (0.center) to (4);
		\draw (1.center) to (4);
		\draw [bend left=90, looseness=1.75] (1.center) to (2.center);
		\draw (3.center) to (8);
		\draw (8) to (2.center);
	\end{pgfonlayer}
\end{tikzpicture}
~~~~~~
(c) ~\begin{tikzpicture}
	\begin{pgfonlayer}{nodelayer}
		\node [style=none] (0) at (-3.25, 0.25) {};
		\node [style=none] (1) at (-3.25, 1.5) {};
		\node [style=none] (2) at (-4.75, 1.5) {};
		\node [style=none] (3) at (-4.75, 0.25) {};
		\node [style=circle, scale=1.8] (4) at (-3.25, 1) {};
		\node [style=none] (5) at (-3.25, 1) {$e_B$};
		\node [style=none] (6) at (-3, 0.5) {$B$};
		\node [style=none] (7) at (-5, 0.5) {$A$};
	\end{pgfonlayer}
	\begin{pgfonlayer}{edgelayer}
		\draw (0.center) to (4);
		\draw (1.center) to (4);
		\draw [bend right=90, looseness=1.75] (1.center) to (2.center);
		\draw (3.center) to (2.center);
	\end{pgfonlayer}
\end{tikzpicture}  =\begin{tikzpicture}
	\begin{pgfonlayer}{nodelayer}
		\node [style=none] (0) at (-5, 0.25) {};
		\node [style=none] (1) at (-5, 1.5) {};
		\node [style=none] (2) at (-3.5, 1.5) {};
		\node [style=none] (3) at (-3.5, 0.25) {};
		\node [style=circle, scale=1.8] (4) at (-5, 1) {};
		\node [style=none] (5) at (-5, 1) {$e_A$};
		\node [style=none] (6) at (-5.25, 0.5) {$A$};
		\node [style=none] (7) at (-3.25, 0.5) {$B$};
		\node [style=circle, scale=1.8] (8) at (-3.5, 1) {};
		\node [style=none] (9) at (-3.5, 1) {$e_B$};
	\end{pgfonlayer}
	\begin{pgfonlayer}{edgelayer}
		\draw (0.center) to (4);
		\draw (1.center) to (4);
		\draw [bend left=90, looseness=1.75] (1.center) to (2.center);
		\draw (3.center) to (8);
		\draw (8) to (2.center);
	\end{pgfonlayer}
\end{tikzpicture}
~~~
(d) ~ \begin{tikzpicture}
	\begin{pgfonlayer}{nodelayer}
		\node [style=none] (0) at (-3.5, 1.5) {};
		\node [style=none] (1) at (-3.5, 0.25) {};
		\node [style=none] (2) at (-5, 0.25) {};
		\node [style=none] (3) at (-5, 1.5) {};
		\node [style=circle, scale=1.8] (4) at (-3.5, 0.75) {};
		\node [style=none] (5) at (-3.5, 0.75) {$e_A$};
		\node [style=none] (6) at (-5.25, 1.25) {$B$};
		\node [style=none] (7) at (-3.25, 1.25) {$A$};
		\node [style=none] (8) at (-4.25, 2) {};
	\end{pgfonlayer}
	\begin{pgfonlayer}{edgelayer}
		\draw (0.center) to (4);
		\draw (1.center) to (4);
		\draw [bend left=90, looseness=1.75] (1.center) to (2.center);
		\draw (3.center) to (2.center);
	\end{pgfonlayer}
\end{tikzpicture}
 =\begin{tikzpicture}
	\begin{pgfonlayer}{nodelayer}
		\node [style=none] (0) at (-3.25, 1.5) {};
		\node [style=none] (1) at (-3.25, 0.25) {};
		\node [style=none] (2) at (-4.75, 0.25) {};
		\node [style=none] (3) at (-4.75, 1.5) {};
		\node [style=circle, scale=1.8] (4) at (-3.25, 0.75) {};
		\node [style=none] (5) at (-3.25, 0.75) {$e_A$};
		\node [style=none] (6) at (-3, 1.25) {$A$};
		\node [style=none] (7) at (-5, 1.25) {$B$};
		\node [style=circle, scale=1.8] (8) at (-4.75, 0.75) {};
		\node [style=none] (9) at (-4.75, 0.75) {$e_B$};
		\node [style=none] (10) at (-4, 2) {};
	\end{pgfonlayer}
	\begin{pgfonlayer}{edgelayer}
		\draw (0.center) to (4);
		\draw (1.center) to (4);
		\draw [bend left=90, looseness=1.75] (1.center) to (2.center);
		\draw (3.center) to (8);
		\draw (8) to (2.center);
	\end{pgfonlayer}
\end{tikzpicture}
\end{equation*}
where $e_A:=\u\v$ and $e_B:= \v\u$.
\end{definition}
The idempotent pair $(e_A, e_B)$ is a morphism of duals only when the binary idempotent is both sectional and retractional. 

\begin{lemma} 
	\label{Lemma: sectional dual morphism}
	In an LDC, a binary idempotent $(\u,\v)$ on a dual $(\eta, \epsilon): A \dashvv B$, 
with splitting $A \to^r E \to^s$ A and $B \to^{r'} E' \to^{s'} B$ is sectional (respectively retractional) if and only if the 
section $(s, r')$ (respectively the retraction $(r,s')$) is a morphism for $(\eta(r \oa r'), (s' \ox s) \epsilon): E \dashvv E'$.
\end{lemma}


Splitting binary idempotents which are either sectional or retractional on a dual produces a self-duality.

\medskip
We next observe that the dagger of a dual is itself a dual:
\begin{lemma} \label{daggering-a-dual}
	Suppose $\X$ is a $\dagger$-LDC, and $(\eta, \epsilon): A \dashvv B$ is a dual in $\X$. 
	Then, $(\epsilon\dagger, \eta\dagger): B^\dagger \dashvv A^\dagger$ is a dual where:

	\vspace{-3em}
	
	\begin{align*}
	\epsilon\dagger &:= \top \to^{\lambda_\top} \bot^\dagger \to^{\epsilon^\dagger} 
	(B \ox A)^\dagger \to^{\lambda_\oa^{-1}} B^\dagger \oa A^\dagger \\
	\eta\dagger &:= A^\dagger \ox B^\dagger \to^{\lambda_\ox} (A \oa B)^\dagger 
	\to^{\eta^\dagger} \top^\dagger \to^{\lambda_\bot^{-1}} \bot
	\end{align*}
\end{lemma}

\begin{definition}
	\label{defn: right dagger dual}
	In a $\dagger$-LDC, a {\bf $\dagger$-dual}, $A \dagdual A^\dagger$ is a dual $(\eta, \epsilon):A \dashvv A^\dagger$ such that 
	\[(\iota_A, 1_{A^\dagger}): (\eta, \epsilon):A \dashvv A^\dagger \to (\eta^\dagger,\epsilon^\dagger):A^{\dagger\dagger} \dashvv A^\dagger\]
        is an isomorphism of duals (see \ref{Eqn: left dagger dual}  (a), (b)). A {\bf self $\dagger$-dual} is a right $\dagger$-dual with 
		an isomorphism $\alpha: A \to A^\dagger$ such that $\alpha \alpha^{-1 \dagger}= \iota$.
		A {\bf morphism of $\dagger$-duals} consists of a pair of maps $(f, f^\dagger): ((\eta, \epsilon): A \dagdual A^\dagger) \to ((\eta', \epsilon'): B \dagdual B^\dagger)$ 
		which are morphism of duals. 
\end{definition}
$(\iota_A, 1_{A^\dagger})$ being an isomorphism of the duals means that the following equations hold: 
	\begin{align}
		\label{Eqn: left dagger dual}
		(a)~~~ \begin{tikzpicture}
			\begin{pgfonlayer}{nodelayer}
				\node [style=none] (0) at (0.5, 0.25) {};
				\node [style=none] (1) at (-1, 1.5) {};
				\node [style=none] (2) at (0.5, 1.5) {};
				\node [style=none] (3) at (-0.25, 2.5) {$\eta$};
				\node [style=none] (4) at (-1, 0.25) {};
				\node [style=none] (5) at (-1.5, 0.25) {$A^{\dagger  \dagger}$};
				\node [style=none] (6) at (1, 0.25) {$A^{\dagger}$};
				\node [style=circle, scale=1.5] (7) at (-1, 1) {};
				\node [style=none] (8) at (-1, 1) {$\iota$};
			\end{pgfonlayer}
			\begin{pgfonlayer}{edgelayer}
				\draw (0.center) to (2.center);
				\draw (1.center) to (7);
				\draw (4.center) to (7);
				\draw [bend left=90, looseness=1.75] (1.center) to (2.center);
			\end{pgfonlayer}
		\end{tikzpicture}  =
		\begin{tikzpicture}
			\begin{pgfonlayer}{nodelayer}
				\node [style=none] (0) at (0.25, 3.25) {};
				\node [style=none] (1) at (-1, 5.5) {};
				\node [style=none] (2) at (0.4999997, 5.5) {};
				\node [style=none] (3) at (-0.25, 4.5) {$\epsilon$};
				\node [style=none] (4) at (-0.75, 3.25) {};
				\node [style=none] (5) at (-0.65, 5.25) {$A^\dagger$};
				\node [style=none] (6) at (-1.25, 3.25) {$A^{\dagger \dagger}$};
				\node [style=none] (7) at (0.5, 3.25) {$A^\dagger$};
				\node [style=none] (8) at (0.75, 4) {};
				\node [style=none] (9) at (-1.25, 4) {};
				\node [style=none] (10) at (-1.25, 5.5) {};
				\node [style=none] (11) at (0.75, 5.5) {};
				\node [style=none] (12) at (0.25, 4) {};
				\node [style=none] (13) at (-0.75, 4) {};
				\node [style=none] (14) at (0.25, 5.25) {$A$};
			\end{pgfonlayer}
			\begin{pgfonlayer}{edgelayer}
				\draw (8.center) to (9.center);
				\draw (10.center) to (9.center);
				\draw (10.center) to (11.center);
				\draw (11.center) to (8.center);
				\draw (13.center) to (4.center);
				\draw (0.center) to (12.center);
				\draw [bend right=90, looseness=1.75] (1.center) to (2.center);
			\end{pgfonlayer}
		\end{tikzpicture}
		~~~~~~~ \text{(or equivalently)} ~~~~~~~ (b)~~~
		\begin{tikzpicture}
			\begin{pgfonlayer}{nodelayer}
				\node [style=none] (0) at (0.25, 2.25) {};
				\node [style=none] (1) at (-1.25, 0.75) {};
				\node [style=none] (2) at (0.25, 0.75) {};
				\node [style=none] (3) at (-0.5, -0.25) {$\epsilon$};
				\node [style=none] (4) at (-1.25, 2.25) {};
				\node [style=none] (5) at (-1.25, 1.5) {};
				\node [style=none] (6) at (-1.5, 2.25) {$A^{\dagger}$};
				\node [style=none] (7) at (0.5, 2.25) {$A$};
			\end{pgfonlayer}
			\begin{pgfonlayer}{edgelayer}
				\draw (1.center) to (4.center);
				\draw (2.center) to (0.center);
				\draw [bend right=90, looseness=1.75] (1.center) to (2.center);
			\end{pgfonlayer}
		\end{tikzpicture}
		= \begin{tikzpicture}
			\begin{pgfonlayer}{nodelayer}
				\node [style=none] (0) at (-0.75, 2.25) {};
				\node [style=none] (1) at (0.5, -0) {};
				\node [style=none] (2) at (-1, -0) {};
				\node [style=none] (3) at (-0.25, 1) {$\eta$};
				\node [style=none] (4) at (0.2499997, 2.25) {};
				\node [style=none] (5) at (0.15, 0.25) {$A^\dagger$};
				\node [style=none] (6) at (0.75, 2.25) {$A$};
				\node [style=none] (7) at (-1.25, 2.25) {$A^{\dagger}$};
				\node [style=none] (8) at (-1.25, 1.25) {};
				\node [style=none] (9) at (0.75, 1.25) {};
				\node [style=none] (10) at (0.75, 0) {};
				\node [style=none] (11) at (-1.25, 0) {};
				\node [style=none] (12) at (-0.75, 1.25) {};
				\node [style=none] (13) at (0.2499997, 1.25) {};
				\node [style=none] (14) at (-0.75, 0.25) {$A$};
				\node [style=circle, scale=1.5] (15) at (0.25, 1.75) {};
				\node [style=none] (16) at (0.25, 1.75) {$\iota$};
			\end{pgfonlayer}
			\begin{pgfonlayer}{edgelayer}
				\draw (8.center) to (9.center);
				\draw (10.center) to (9.center);
				\draw (10.center) to (11.center);
				\draw (11.center) to (8.center);
				\draw (12.center) to (0.center);
				\draw (13.center) to (15);
				\draw (4.center) to (15);
				\draw [bend left=90, looseness=1.75] (2.center) to (1.center);
			\end{pgfonlayer}
		\end{tikzpicture}
		\end{align} 

Lemma \ref{Lemma: sectional dual morphism} can be lifted to 
$\dagger$-idempotents on $\dagger$-duals which $\dagger$-splits 
to produce a self-$\dagger$-dual \cite[Lemma 4.11]{CoS21}.

\subsection{Linear monoid}

The simplest way to describe a linear monoid is as a  $\ox$-monoid  on an object together 
with a dual for that object.  Their similarity to Frobenius algebras becomes more apparent when one regards 
a linear monoid as a $\ox$-monoid and a $\oa$-comonoid with actions and coactions. 

\begin{definition}
	A  {\bf linear monoid} \cite{CKS00,Egg10}, $A \linmonw B$, in an LDC consists of a monoid $(A,e: \top \to A,m: A \ox A \to A)$, 
	a left dual  $(\eta_L, \epsilon_L): A \dashvv B$, and  a right dual $(\eta_R, \epsilon_R): B \dashvv A$ such that: 
	\begin{equation}
		\label{eqn: linear monoid} 
		\begin{tikzpicture}
			\begin{pgfonlayer}{nodelayer}
				\node [style=none] (12) at (1.75, 0.25) {$B$};
				\node [style=none] (13) at (0.25, 0.25) {$B$};
				\node [style=circle] (17) at (1, 1.25) {};
				\node [style=none] (18) at (1, 2.5) {};
				\node [style=none] (19) at (0.5, 0) {};
				\node [style=none] (20) at (1.5, 0) {};
				\node [style=none] (21) at (1.25, 2.25) {$B$};
			\end{pgfonlayer}
			\begin{pgfonlayer}{edgelayer}
				\draw [in=-150, out=90, looseness=1.25] (19.center) to (17);
				\draw [in=90, out=-30, looseness=1.25] (17) to (20.center);
				\draw (17) to (18.center);
			\end{pgfonlayer}
		\end{tikzpicture} := 		
		 \begin{tikzpicture}
		\begin{pgfonlayer}{nodelayer}
			\node [style=circle] (0) at (-2.75, 0.75) {};
			\node [style=none] (1) at (-3.25, 1.25) {};
			\node [style=none] (2) at (-2.75, 0.5) {};
			\node [style=none] (3) at (-2.25, 1.25) {};
			\node [style=none] (4) at (-1.5, 1.25) {};
			\node [style=none] (5) at (-1.5, -0) {};
			\node [style=none] (6) at (-3.25, 1.5) {};
			\node [style=none] (7) at (-1, 1.5) {};
			\node [style=none] (8) at (-1, -0) {};
			\node [style=none] (9) at (-3.75, 0.5) {};
			\node [style=none] (10) at (-3.75, 2.25) {};
			\node [style=none] (11) at (-4, 2) {$B$};
			\node [style=none] (12) at (-1.75, 0.25) {$B$};
			\node [style=none] (13) at (-0.75, 0.25) {$B$};
			\node [style=none] (14) at (-2.25, 2.5) {$\eta_L$};
			\node [style=none] (15) at (-2, 1.75) {$\eta_L$};
			\node [style=none] (16) at (-3.25, -0.25) {$\epsilon_L$};
		\end{pgfonlayer}
		\begin{pgfonlayer}{edgelayer}
			\draw [in=150, out=-90, looseness=1.00] (1.center) to (0);
			\draw [in=-90, out=30, looseness=1.00] (0) to (3.center);
			\draw (0) to (2.center);
			\draw [bend left=90, looseness=1.50] (2.center) to (9.center);
			\draw (9.center) to (10.center);
			\draw (6.center) to (1.center);
			\draw [bend left=90, looseness=1.25] (6.center) to (7.center);
			\draw (7.center) to (8.center);
			\draw (4.center) to (5.center);
			\draw [bend right=90, looseness=1.25] (4.center) to (3.center);
		\end{pgfonlayer}
	\end{tikzpicture} = \begin{tikzpicture}
		\begin{pgfonlayer}{nodelayer}
			\node [style=circle] (0) at (-2, 0.75) {};
			\node [style=none] (1) at (-1.5, 1.25) {};
			\node [style=none] (2) at (-2, 0.5) {};
			\node [style=none] (3) at (-2.5, 1.25) {};
			\node [style=none] (4) at (-3.25, 1.25) {};
			\node [style=none] (5) at (-3.25, -0) {};
			\node [style=none] (6) at (-1.5, 1.5) {};
			\node [style=none] (7) at (-3.75, 1.5) {};
			\node [style=none] (8) at (-3.75, -0) {};
			\node [style=none] (9) at (-1, 0.5) {};
			\node [style=none] (10) at (-1, 2.25) {};
			\node [style=none] (11) at (-0.75, 2) {$B$};
			\node [style=none] (12) at (-3, 0.25) {$B$};
			\node [style=none] (13) at (-4, 0.25) {$B$};
			\node [style=none] (14) at (-2.5, 2.5) {$\eta_R$};
			\node [style=none] (15) at (-2.75, 1.75) {$\eta_R$};
			\node [style=none] (16) at (-1.5, -0.25) {$\epsilon_R$};
		\end{pgfonlayer}
		\begin{pgfonlayer}{edgelayer}
			\draw [in=30, out=-90, looseness=1.00] (1.center) to (0);
			\draw [in=-90, out=150, looseness=1.00] (0) to (3.center);
			\draw (0) to (2.center);
			\draw [bend right=90, looseness=1.50] (2.center) to (9.center);
			\draw (9.center) to (10.center);
			\draw (6.center) to (1.center);
			\draw [bend right=90, looseness=1.25] (6.center) to (7.center);
			\draw (7.center) to (8.center);
			\draw (4.center) to (5.center);
			\draw [bend left=90, looseness=1.25] (4.center) to (3.center);
		\end{pgfonlayer}
	\end{tikzpicture} 
	~~~~~~~~~~~~ 
	\begin{tikzpicture}
		\begin{pgfonlayer}{nodelayer}
			\node [style=circle] (17) at (1, 0.25) {};
			\node [style=none] (18) at (1, 2.5) {};
			\node [style=none] (21) at (1.25, 2.25) {$B$};
		\end{pgfonlayer}
		\begin{pgfonlayer}{edgelayer}
			\draw (17) to (18.center);
		\end{pgfonlayer}
	\end{tikzpicture} := 	
	\begin{tikzpicture}
		\begin{pgfonlayer}{nodelayer}
			\node [style=circle] (0) at (-0.75, 1.5) {};
			\node [style=none] (1) at (-0.75, 0.5) {};
			\node [style=none] (2) at (-1.75, 0.5) {};
			\node [style=none] (3) at (-1.75, 2.25) {};
			\node [style=none] (4) at (-2, 2) {$B$};
			\node [style=none] (5) at (-1.25, -0.25) {$\epsilon_L$};
		\end{pgfonlayer}
		\begin{pgfonlayer}{edgelayer}
			\draw (0) to (1.center);
			\draw [bend left=90, looseness=1.50] (1.center) to (2.center);
			\draw (2.center) to (3.center);
		\end{pgfonlayer}
	\end{tikzpicture} = 
	\begin{tikzpicture}
		\begin{pgfonlayer}{nodelayer}
			\node [style=circle] (0) at (-2, 1.5) {};
			\node [style=none] (1) at (-2, 0.5) {};
			\node [style=none] (2) at (-1, 0.5) {};
			\node [style=none] (3) at (-1, 2.25) {};
			\node [style=none] (4) at (-0.75, 2) {$B$};
			\node [style=none] (5) at (-1.5, -0.25) {$\epsilon_R$};
		\end{pgfonlayer}
		\begin{pgfonlayer}{edgelayer}
			\draw (0) to (1.center);
			\draw [bend right=90, looseness=1.50] (1.center) to (2.center);
			\draw (2.center) to (3.center);
		\end{pgfonlayer}
	\end{tikzpicture} 
\end{equation}
\end{definition}
In a symmetric LDC, a linear monoid is {\bf symmetric} when its duals are symmetric, i.e, $\eta_R = \eta_L c_\oa$ and 
$\epsilon_R = c_\ox \epsilon_L$.  A symmetric linear monoid is determined by a monoid $(A, m, u)$ and a dual $(\eta, \epsilon): A \dashvv B$. 
There is a more useful form for linear monoids in which their similarity to the usual description of Frobenius algebras in CQM is evident:

\begin{proposition}
	\label{Lemma: alternate presentation of linear monoids}
	A linear monoid, $A \linmonw B$, in an LDC is equivalent to the following data:
	\begin{itemize}
	\item a monoid $(A, \mulmap{1.5}{white}: A \ox A \to A, \unitmap{1.5}{white}: \top \to A) $
	\item a comonoid $(B, \comulmap{1.5}{white}: B \to B \oa B, \counitmap{1.5}{white}: B \to \bot) $
	\item actions, $\leftaction{0.5}{white}: A \ox B \to B$, $\rightaction{0.5}{white}: B \ox A \to B$,
	and coactions $\leftcoaction{0.55}{white}: A \to B \oa A$, $\rightcoaction{0.5}{white}: A \to A \oa B$,
    \end{itemize}
	such that the following axioms (and their `op' and `co' symmetric forms) hold:  
 
	$(a) ~~\begin{tikzpicture}
		\begin{pgfonlayer}{nodelayer}
			\node [style=none] (0) at (2, 0.25) {};
			\node [style=none] (1) at (1.75, 0.25) {};
			\node [style=none] (2) at (2, 0.5) {};
			\node [style=none] (3) at (2, -1) {};
			\node [style=none] (4) at (2, 1.5) {};
			\node [style=none] (5) at (2, 1.75) {$B$};
			\node [style=circle, scale=0.75] (6) at (1, 1.5) {};
			\node [style=none] (7) at (2, 0.5) {};
			\node [style=none] (8) at (2, 0) {};
		\end{pgfonlayer}
		\begin{pgfonlayer}{edgelayer}
			\draw (3.center) to (0.center);
			\draw (2.center) to (0.center);
			\draw (4.center) to (2.center);
			\draw [in=143, out=-90] (6) to (1.center);
			\draw [bend right=90, looseness=1.50] (7.center) to (8.center);
		\end{pgfonlayer}
	\end{tikzpicture} = \begin{tikzpicture}
	\begin{pgfonlayer}{nodelayer}
	\node [style=none] (0) at (1, -1) {};
	\node [style=none] (1) at (1, 1.5) {};
	\node [style=none] (2) at (1, 1.75) {$B$};
	\end{pgfonlayer}
	\begin{pgfonlayer}{edgelayer}
	\draw (1.center) to (0.center);
	\end{pgfonlayer}
	\end{tikzpicture} 	
	~~~~~~(b)~~
	\begin{tikzpicture}
		\begin{pgfonlayer}{nodelayer}
			\node [style=none] (0) at (2.5, -0.5) {};
			\node [style=none] (1) at (2.25, -0.25) {};
			\node [style=none] (2) at (2.5, -0.25) {};
			\node [style=none] (3) at (2.5, -1) {};
			\node [style=none] (4) at (1, 1.5) {};
			\node [style=none] (5) at (2.5, 1.5) {};
			\node [style=none] (6) at (2, 1.5) {};
			\node [style=none] (7) at (1, 1.75) {$A$};
			\node [style=none] (8) at (2.5, 1.75) {$B$};
			\node [style=none] (9) at (2, 1.75) {$A$};
			\node [style=circle, scale=0.75] (10) at (1.5, 0.75) {};
			\node [style=none] (11) at (2.5, 0) {};
			\node [style=none] (12) at (2.5, -0.5) {};
			\node [style=none] (13) at (2.75, -0.75) {$B$};
		\end{pgfonlayer}
		\begin{pgfonlayer}{edgelayer}
			\draw (3.center) to (0.center);
			\draw (2.center) to (0.center);
			\draw (5.center) to (2.center);
			\draw [bend left] (6.center) to (10);
			\draw [bend left] (10) to (4.center);
			\draw [bend right] (10) to (1.center);
			\draw [bend right=90, looseness=1.75] (11.center) to (12.center);
		\end{pgfonlayer}
	\end{tikzpicture} = \begin{tikzpicture}
		\begin{pgfonlayer}{nodelayer}
			\node [style=none] (14) at (5.25, -0.5) {};
			\node [style=none] (15) at (5, -0.25) {};
			\node [style=none] (17) at (5.25, -1) {};
			\node [style=none] (18) at (3.75, 1.5) {};
			\node [style=none] (19) at (5, 0.75) {};
			\node [style=none] (20) at (5.25, 1.5) {};
			\node [style=none] (21) at (5.25, 1) {};
			\node [style=none] (22) at (4.5, 1.5) {};
			\node [style=none] (23) at (5.25, 0.5) {};
			\node [style=none] (24) at (3.75, 1.75) {$A$};
			\node [style=none] (25) at (5.25, 1.75) {$B$};
			\node [style=none] (26) at (4.5, 1.75) {$A$};
			\node [style=none] (27) at (5.25, 0) {};
			\node [style=none] (28) at (5.25, -0.5) {};
			\node [style=none] (29) at (5.25, 1) {};
			\node [style=none] (30) at (5.25, 0.5) {};
		\end{pgfonlayer}
		\begin{pgfonlayer}{edgelayer}
			\draw (17.center) to (14.center);
			\draw (21.center) to (20.center);
			\draw [bend left] (19.center) to (22.center);
			\draw (21.center) to (23.center);
			\draw (23.center) to (14.center);
			\draw [bend left] (15.center) to (18.center);
			\draw [bend right=90, looseness=1.75] (27.center) to (28.center);
			\draw [bend right=90, looseness=1.75] (29.center) to (30.center);
		\end{pgfonlayer}
	\end{tikzpicture}	
	~~~~~~(c)~~
	\begin{tikzpicture}
		\begin{pgfonlayer}{nodelayer}
			\node [style=none] (31) at (7.25, -0.5) {};
			\node [style=none] (32) at (7.25, -0.25) {};
			\node [style=none] (33) at (7, -0.25) {};
			\node [style=none] (34) at (7.25, -1) {};
			\node [style=none] (35) at (6.25, 1.5) {};
			\node [style=none] (36) at (7.5, 0.5) {};
			\node [style=none] (37) at (7.25, 1.5) {};
			\node [style=none] (38) at (7.25, 0.75) {};
			\node [style=none] (39) at (8, 1.5) {};
			\node [style=none] (40) at (7.25, 0.5) {};
			\node [style=none] (41) at (6.25, 1.75) {$A$};
			\node [style=none] (42) at (7.25, 1.75) {$B$};
			\node [style=none] (43) at (8, 1.75) {$A$};
			\node [style=none] (44) at (7.25, 0.75) {};
			\node [style=none] (45) at (7.25, 0.25) {};
			\node [style=none] (46) at (7.25, 0.75) {};
			\node [style=none] (47) at (7.25, 0.25) {};
			\node [style=none] (48) at (7.25, 0) {};
			\node [style=none] (49) at (7.25, -0.5) {};
			\node [style=none] (50) at (7.25, 0) {};
			\node [style=none] (51) at (7.25, -0.5) {};
		\end{pgfonlayer}
		\begin{pgfonlayer}{edgelayer}
			\draw (34.center) to (31.center);
			\draw [in=270, out=150] (33.center) to (35.center);
			\draw (32.center) to (31.center);
			\draw (38.center) to (37.center);
			\draw [bend right] (36.center) to (39.center);
			\draw (38.center) to (40.center);
			\draw (40.center) to (31.center);
			\draw [bend left=90, looseness=1.75] (46.center) to (47.center);
			\draw [bend right=90, looseness=1.75] (50.center) to (51.center);
		\end{pgfonlayer}
	\end{tikzpicture} = \begin{tikzpicture}
		\begin{pgfonlayer}{nodelayer}
			\node [style=none] (31) at (7, -0.5) {};
			\node [style=none] (32) at (7, -0.25) {};
			\node [style=none] (33) at (7.25, -0.25) {};
			\node [style=none] (34) at (7, -1) {};
			\node [style=none] (35) at (8, 1.5) {};
			\node [style=none] (36) at (6.75, 0.5) {};
			\node [style=none] (37) at (7, 1.5) {};
			\node [style=none] (38) at (7, 0.75) {};
			\node [style=none] (39) at (6.25, 1.5) {};
			\node [style=none] (40) at (7, 0.5) {};
			\node [style=none] (41) at (8, 1.75) {$A$};
			\node [style=none] (42) at (7, 1.75) {$B$};
			\node [style=none] (43) at (6.25, 1.75) {$A$};
			\node [style=none] (44) at (7, 0.75) {};
			\node [style=none] (45) at (7, 0.25) {};
			\node [style=none] (46) at (7, 0.75) {};
			\node [style=none] (47) at (7, 0.25) {};
			\node [style=none] (48) at (7, 0) {};
			\node [style=none] (49) at (7, -0.5) {};
			\node [style=none] (50) at (7, 0) {};
			\node [style=none] (51) at (7, -0.5) {};
		\end{pgfonlayer}
		\begin{pgfonlayer}{edgelayer}
			\draw (34.center) to (31.center);
			\draw [in=-90, out=30] (33.center) to (35.center);
			\draw (32.center) to (31.center);
			\draw (38.center) to (37.center);
			\draw [bend left] (36.center) to (39.center);
			\draw (38.center) to (40.center);
			\draw (40.center) to (31.center);
			\draw [bend right=90, looseness=1.75] (46.center) to (47.center);
			\draw [bend left=90, looseness=1.75] (50.center) to (51.center);
		\end{pgfonlayer}
	\end{tikzpicture}
	 ~~~~~ (d) ~~ 
	 \begin{tikzpicture}
		\begin{pgfonlayer}{nodelayer}
			\node [style=none] (38) at (10.25, 1.25) {};
			\node [style=none] (44) at (10.25, 1.25) {};
			\node [style=none] (45) at (10.25, 0.75) {};
			\node [style=none] (46) at (10.25, 1.25) {};
			\node [style=none] (47) at (10.25, 0.75) {};
			\node [style=none] (52) at (10.25, 1.25) {};
			\node [style=none] (53) at (10.25, 1) {};
			\node [style=none] (54) at (10, 1) {};
			\node [style=none] (55) at (10.25, 1.75) {};
			\node [style=none] (56) at (9.5, -0.75) {};
			\node [style=none] (57) at (10.25, 2) {$A$};
			\node [style=none] (58) at (9.5, -1) {$B$};
			\node [style=none] (59) at (9.5, -0.25) {};
			\node [style=none] (60) at (8.75, 1.75) {};
			\node [style=none] (61) at (9.5, 0) {};
			\node [style=none] (62) at (9.25, 0) {};
			\node [style=none] (63) at (10.25, -0.75) {};
			\node [style=none] (64) at (8.75, 2) {$A$};
			\node [style=none] (65) at (10.25, -1) {$A$};
			\node [style=none] (66) at (9.5, 1) {$B$};
			\node [style=none] (67) at (9.5, 0.25) {};
			\node [style=none] (68) at (9.5, 0.25) {};
			\node [style=none] (69) at (9.5, -0.25) {};
			\node [style=none] (70) at (9.5, 0.25) {};
			\node [style=none] (71) at (9.5, -0.25) {};
		\end{pgfonlayer}
		\begin{pgfonlayer}{edgelayer}
			\draw [bend right=90, looseness=1.75] (46.center) to (47.center);
			\draw (55.center) to (52.center);
			\draw [in=90, out=-150] (54.center) to (56.center);
			\draw (53.center) to (52.center);
			\draw [in=-90, out=165] (62.center) to (60.center);
			\draw (61.center) to (59.center);
			\draw (53.center) to (63.center);
			\draw [bend right=90, looseness=1.75] (70.center) to (71.center);
		\end{pgfonlayer}
	\end{tikzpicture} = \begin{tikzpicture}
		\begin{pgfonlayer}{nodelayer}
			\node [style=none] (72) at (12, 0.25) {};
			\node [style=none] (73) at (12, 0) {};
			\node [style=none] (74) at (11.75, 0.25) {};
			\node [style=none] (75) at (11.25, -0.75) {};
			\node [style=none] (76) at (12.5, 2) {$A$};
			\node [style=none] (77) at (11.25, -1) {$B$};
			\node [style=none] (78) at (12, -0.75) {};
			\node [style=none] (79) at (12, -1) {$A$};
			\node [style=circle] (80) at (12, 1) {};
			\node [style=none] (81) at (12.5, 1.75) {};
			\node [style=none] (82) at (11.5, 1.75) {};
			\node [style=none] (83) at (11.5, 2) {$A$};
			\node [style=none] (84) at (12, 0.5) {};
			\node [style=none] (85) at (12, 0) {};
		\end{pgfonlayer}
		\begin{pgfonlayer}{edgelayer}
			\draw [in=90, out=-135] (74.center) to (75.center);
			\draw (73.center) to (72.center);
			\draw (73.center) to (78.center);
			\draw [bend right] (80) to (81.center);
			\draw [bend left] (80) to (82.center);
			\draw (80) to (72.center);
			\draw [bend right=90, looseness=1.75] (84.center) to (85.center);
		\end{pgfonlayer}
	\end{tikzpicture}= \begin{tikzpicture}
		\begin{pgfonlayer}{nodelayer}
			\node [style=none] (74) at (11.25, 1.25) {};
			\node [style=none] (75) at (10.75, -0.75) {};
			\node [style=none] (76) at (11.5, 2.25) {$A$};
			\node [style=none] (77) at (10.75, -1) {$B$};
			\node [style=none] (78) at (12, -0.75) {};
			\node [style=none] (79) at (12, -1) {$A$};
			\node [style=circle] (80) at (12, 0) {};
			\node [style=none] (81) at (12.5, 0.5) {};
			\node [style=none] (82) at (11.5, 0.5) {};
			\node [style=none] (83) at (12.5, 2.25) {$A$};
			\node [style=none] (84) at (11.5, 1.5) {};
			\node [style=none] (85) at (11.5, 1) {};
			\node [style=none] (86) at (11.5, 2) {};
			\node [style=none] (87) at (12.5, 2) {};
		\end{pgfonlayer}
		\begin{pgfonlayer}{edgelayer}
			\draw [in=90, out=-135] (74.center) to (75.center);
			\draw [bend right] (80) to (81.center);
			\draw [bend left] (80) to (82.center);
			\draw [bend right=90, looseness=1.75] (84.center) to (85.center);
			\draw (80) to (78.center);
			\draw (82.center) to (86.center);
			\draw (87.center) to (81.center);
		\end{pgfonlayer}
	\end{tikzpicture} $

	\noindent If the linear monoid is symmetric, then: 
$ \begin{tikzpicture}[scale=1.1]
	\begin{pgfonlayer}{nodelayer}
		\node [style=none] (0) at (-0.75, 2.75) {};
		\node [style=none] (1) at (-0.75, 3) {};
		\node [style=none] (2) at (-0.5, 3) {};
		\node [style=none] (4) at (0, 3.25) {};
		\node [style=none] (7) at (-0.75, 3.25) {};
		\node [style=circle] (8) at (-0.75, 2) {};
		\node [style=none] (9) at (0, 4.25) {};
		\node [style=none] (10) at (-0.75, 4.25) {};
		\node [style=none] (11) at (-1, 4) {$A$};
		\node [style=none] (12) at (0.25, 4) {$B$};
		\node [style=none] (13) at (-0.75, 3.25) {};
		\node [style=none] (14) at (-0.75, 2.75) {};
	\end{pgfonlayer}
	\begin{pgfonlayer}{edgelayer}
		\draw [bend right=45] (2.center) to (4.center);
		\draw (0.center) to (1.center);
		\draw (1.center) to (7.center);
		\draw (8) to (0.center);
		\draw [in=-90, out=90] (4.center) to (10.center);
		\draw [in=-90, out=90, looseness=1.25] (7.center) to (9.center);
		\draw [bend left=90, looseness=1.75] (13.center) to (14.center);
	\end{pgfonlayer}
\end{tikzpicture} = \begin{tikzpicture}[scale=1.1]
	\begin{pgfonlayer}{nodelayer}
		\node [style=none] (0) at (0, 2.75) {};
		\node [style=none] (1) at (0, 3) {};
		\node [style=none] (2) at (-0.25, 3) {};
		\node [style=none] (4) at (-0.75, 3.25) {};
		\node [style=none] (7) at (0, 3.25) {};
		\node [style=circle] (8) at (0, 2) {};
		\node [style=none] (9) at (0, 4.25) {};
		\node [style=none] (10) at (-0.75, 4.25) {};
		\node [style=none] (11) at (-1, 4) {$A$};
		\node [style=none] (12) at (0.25, 4) {$B$};
		\node [style=none] (13) at (0, 3.25) {};
		\node [style=none] (14) at (0, 2.75) {};
	\end{pgfonlayer}
	\begin{pgfonlayer}{edgelayer}
		\draw [bend left=45] (2.center) to (4.center);
		\draw (0.center) to (1.center);
		\draw (1.center) to (7.center);
		\draw (8) to (0.center);
		\draw [in=-90, out=90] (4.center) to (10.center);
		\draw [in=-90, out=90, looseness=1.25] (7.center) to (9.center);
		\draw [bend right=90, looseness=1.75] (13.center) to (14.center);
	\end{pgfonlayer}
\end{tikzpicture}$
\end{proposition}

In a linear bicategory, the structure described in Proposition \ref{Lemma: alternate presentation of linear monoids} 
is called a linear monad \cite{CKS00}: here, as we are in a 
simpler context, we use linear monoid. A linear monoid $A \linmonw B$, 
in a monoidal category gives a Frobenius algebra when it is a self-linear monoid that is $A = B$ and the duality coincide with the self-dual cup and cap. 
Note that 
while a Frobenius algebra is always on a self-dual object, a linear monoid allows Frobenius interaction 
between distinct objects which are duals of one another.

\begin{definition}
A {\bf morphism of linear monoids} is a pair of maps, $(f,g): (A \linmonw B) \to (A' \linmonw B')$, such that 
$f: A \to A'$ is a monoid morphism (or equivalently $g: B' \to B$ is a comonoid morphism),
 and $(f, g)$ and $(g,f)$ preserve the left and the right duals respectively. 
 \end{definition}

Note that a morphism of Frobenius algebras is usually given by a single monoid morphism 
which is an isomorphism. However, in the case of a morphism of linear monoids, 
the comonoid morphism, $g: B' \to B$, is the cyclic mate of the monoid morphism, $f: A \to A'$. 
This means that a linear monoid morphism is not restricted to being an isomorphism.  

Given an idempotent $e_A: A \to A$, and a monoid $(A, m, u)$ in a monoidal category, $e_A$ is
{\bf retractional} on the  monoid if $e_A m = e_A m (e_A \ox e_A)$. $e_A$ is 
{\bf sectional} on the monoid if $m (e_A \ox e_A) =  e_A m (e_A \ox e_A)$ and $u e_A = u$.

\begin{lemma} 
\label{Lemma: sectional monoid morphism}
	In a monoidal category, a split idempotent $e: A \to A$ on a monoid $(A,m,u)$, 
with splitting $A \to^r E \to^s$ A, is sectional (respectively retractional) if and only if the 
section $s$ (respectively the retraction $r$) is a monoid morphism for $(E, (s \ox s)m r, u r)$.
\end{lemma}

A binary idempotent $(\u,\v)$ is {\bf sectional} (respectively {\bf retractional}) on a linear monoid when 
$e_A=\u\v$ and $e_B=\v\u$ satisfies the conditions in the following table:

 \[	  \begin{tabular}{l|| l} 
		\hline
		(\u, \v)  \textbf{sectional} on $A \linmonw B$  & (\u, \v)  \textbf{retractional} on $A \linmonw B$\\
		\hline
		 $e_A$ preserves $(A,m,u)$ sectionally &  $e_A$ preserves $(A,m,u)$ retractionally \\
		\hline
		 $(e_A, e_B)$ preserves $(\eta_L, \epsilon_L): A \dashvv B$ sectionally &$(e_A, e_B)$ preserves $(\eta_L, \epsilon_L): A \dashvv B$ retractionally\\
		\hline
		 $(e_B, e_A)$ preserves $(\eta_R, \epsilon_R): B \dashvv A$ retractionally & $(e_B, e_A)$ preserves $(\eta_R, \epsilon_R): B \dashvv A$ sectionally\\
		\hline
	  \end{tabular} \]

Splitting a sectional/retractional binary idempotent on a linear monoid gives a self-linear monoid on the splitting. 

\begin{definition}
A  {\bf $\dagger$-linear monoid}, $(A,  \mulmap{1.2}{white}, \unitmap{1.2}{white}) 
\dagmonw (A^\dagger, \comulmap{1.2}{white}, \counitmap{1.2}{white})$, 
in a $\dagger$-LDC is a linear monoid such that $(\eta_L, \epsilon_L): A \dashvv A^\dagger$ 
and $(\eta_R, \epsilon_R): A^\dagger \dashvv A$ are $\dagger$-duals and:
\begin{equation}
	\label{MainEqn: rightdagmon}
	\begin{tikzpicture}
		\begin{pgfonlayer}{nodelayer}
			\node [style=circle] (17) at (1, 0.25) {};
			\node [style=none] (18) at (1, 2.5) {};
			\node [style=none] (21) at (1.25, 2.25) {$A^\dagger$};
		\end{pgfonlayer}
		\begin{pgfonlayer}{edgelayer}
			\draw (17) to (18.center);
		\end{pgfonlayer}
	\end{tikzpicture} := 	
	\begin{tikzpicture}
		\begin{pgfonlayer}{nodelayer}
			\node [style=circle] (0) at (-0.75, 3.5) {};
			\node [style=none] (1) at (-0.75, 2.5) {};
			\node [style=none] (2) at (-1.75, 2.5) {};
			\node [style=none] (3) at (-1.75, 4) {};
			\node [style=none] (4) at (-1.75, 2.5) {};
			\node [style=none] (5) at (-2.2, 3.75) {$A^\dagger$};
			\node [style=none] (6) at (-1.25, 1.75) {$\epsilon_L$};
		\end{pgfonlayer}
		\begin{pgfonlayer}{edgelayer}
			\draw (0) to (1.center);
			\draw [bend left=75, looseness=1.75] (1.center) to (2.center);
			\draw (2.center) to (3.center);
		\end{pgfonlayer}
	\end{tikzpicture} = \begin{tikzpicture}
		\begin{pgfonlayer}{nodelayer}
			\node [style=circle] (0) at (-1.75, 1.5) {};
			\node [style=none] (1) at (-1.75, 2.5) {};
			\node [style=none] (2) at (-1.4, 0.25) {$A^\dagger$};
			\node [style=none] (3) at (-2.5, 0.75) {};
			\node [style=none] (4) at (-1, 0.75) {};
			\node [style=none] (5) at (-1, 2.5) {};
			\node [style=none] (6) at (-2.5, 2.5) {};
			\node [style=none] (7) at (-1.75, 0.75) {};
			\node [style=none] (8) at (-1.75, -0) {};
		\end{pgfonlayer}
		\begin{pgfonlayer}{edgelayer}
			\draw (0) to (1.center);
			\draw (3.center) to (4.center);
			\draw (4.center) to (5.center);
			\draw (5.center) to (6.center);
			\draw (6.center) to (3.center);
			\draw (8.center) to (7.center);
		\end{pgfonlayer}
	\end{tikzpicture}   ~~~~~~~~
	\begin{tikzpicture}
		\begin{pgfonlayer}{nodelayer}
			\node [style=none] (12) at (2, 0.25) {$A^\dagger$};
			\node [style=none] (13) at (0, 0.25) {$A^\dagger$};
			\node [style=circle] (17) at (1, 1.25) {};
			\node [style=none] (18) at (1, 2.5) {};
			\node [style=none] (19) at (0.5, 0) {};
			\node [style=none] (20) at (1.5, 0) {};
			\node [style=none] (21) at (1.5, 2.25) {$A^\dagger$};
		\end{pgfonlayer}
		\begin{pgfonlayer}{edgelayer}
			\draw [in=-150, out=90, looseness=1.25] (19.center) to (17);
			\draw [in=90, out=-30, looseness=1.25] (17) to (20.center);
			\draw (17) to (18.center);
		\end{pgfonlayer}
	\end{tikzpicture} := \begin{tikzpicture}
		\begin{pgfonlayer}{nodelayer}
			\node [style=circle] (0) at (-2.75, 0.75) {};
			\node [style=none] (1) at (-3.25, 1.25) {};
			\node [style=none] (2) at (-2.75, 0.5) {};
			\node [style=none] (3) at (-2.25, 1.25) {};
			\node [style=none] (4) at (-1.5, 1.25) {};
			\node [style=none] (5) at (-1.5, -0) {};
			\node [style=none] (6) at (-3.25, 1.5) {};
			\node [style=none] (7) at (-1, 1.5) {};
			\node [style=none] (8) at (-1, -0) {};
			\node [style=none] (9) at (-3.75, 0.5) {};
			\node [style=none] (10) at (-3.75, 2.25) {};
			\node [style=none] (11) at (-4.1, 2) {$A^\dagger$};
			\node [style=none] (12) at (-1.85, 0.25) {$A^\dagger$};
			\node [style=none] (13) at (-0.65, 0.25) {$A^\dagger$};
			\node [style=none] (14) at (-2, 2.5) {$\eta_L$};
			\node [style=none] (15) at (-3.25, -0.25) {$\epsilon_L$};
		\end{pgfonlayer}
		\begin{pgfonlayer}{edgelayer}
			\draw [in=150, out=-90, looseness=1.00] (1.center) to (0);
			\draw [in=-90, out=30, looseness=1.00] (0) to (3.center);
			\draw (0) to (2.center);
			\draw [bend left=90, looseness=1.50] (2.center) to (9.center);
			\draw (9.center) to (10.center);
			\draw (6.center) to (1.center);
			\draw [bend left=90, looseness=1.00] (6.center) to (7.center);
			\draw (7.center) to (8.center);
			\draw (4.center) to (5.center);
			\draw [bend right=90, looseness=2.00] (4.center) to (3.center);
		\end{pgfonlayer}
	\end{tikzpicture} = \begin{tikzpicture}
		\begin{pgfonlayer}{nodelayer}
			\node [style=circle] (0) at (-2.75, 1.5) {};
			\node [style=none] (1) at (-3.25, 2.25) {};
			\node [style=none] (2) at (-2.75, 1) {};
			\node [style=none] (3) at (-2.25, 2.25) {};
			\node [style=none] (4) at (-4, 1) {};
			\node [style=none] (5) at (-1.5, 1) {};
			\node [style=none] (6) at (-1.5, 2.25) {};
			\node [style=none] (7) at (-4, 2.25) {};
			\node [style=none] (8) at (-2.75, 2.75) {};
			\node [style=none] (9) at (-2.75, 2.25) {};
			\node [style=none] (10) at (-3.5, 1) {};
			\node [style=none] (11) at (-3.5, 0.25) {};
			\node [style=none] (12) at (-2, 1) {};
			\node [style=none] (13) at (-2, 0.25) {};
			\node [style=none] (14) at (-2.35, 2.75) {$A^\dagger$};
			\node [style=none] (15) at (-1.4, 0.5) {$A^\dagger$};
			\node [style=none] (16) at (-4, 0.5) {$A^\dagger$};
		\end{pgfonlayer}
		\begin{pgfonlayer}{edgelayer}
			\draw [in=150, out=-90] (1.center) to (0);
			\draw [in=-90, out=30] (0) to (3.center);
			\draw (0) to (2.center);
			\draw (7.center) to (6.center);
			\draw (5.center) to (6.center);
			\draw (5.center) to (4.center);
			\draw (4.center) to (7.center);
			\draw (8.center) to (9.center);
			\draw (10.center) to (11.center);
			\draw (12.center) to (13.center);
		\end{pgfonlayer}
	\end{tikzpicture}	
\end{equation}
\end{definition}

A {\bf morphism} of $\dagger$-linear monoids is a pair of maps $(f, f^\dagger)$ which are 
morphisms of underlying linear monoids. 
 Similar to duals, splitting sectional/retractional binary idempotents on a linear monoid induces a 
self-linear monoid. In the presence of dagger, one gets a $\dagger$-self-linear monoid 
\cite[Lemma 5.8]{CoS21}. 

A $\dagger$-linear monoid in a unitary category is equivalent to a $\dagger$-Frobenius algebra under certain conditions: 
\begin{lemma}
	\label{Lemma: dagmondagFrob}
	In a compact LDC, a self-linear monoid $A \linmonw A'$ with an isomorphism 
	$\alpha: A \to A'$ precisely corresponds to a Frobenius algebra under the linear equivalence, 
	${\sf Mx}_\downarrow$ if and only if the linear monoid satisfies the equation  below. 
	In a unitary category, any $\dagger$-linear monoid $A \dagmonwtik A^\dag$ 
	precisely corresponds to a $\dagger$-Frobenius algebra under the same equivalence
	 if and only if  the $\dagger$-linear monoid satisfies the equation below for the 
	 unitary structure isomorphism $\varphi_A: A \to A^\dagger$.
		\begin{equation} 
			\label{eqnn: unitary coincidence}
				\begin{tikzpicture}
				\begin{pgfonlayer}{nodelayer}
					\node [style=none] (0) at (0.75, 1) {};
					\node [style=none] (1) at (0.75, -1.5) {};
					\node [style=none] (2) at (1.25, -1) {$A'$};
					\node [style=circle, scale=2] (3) at (0.75, -0.25) {};
					\node [style=none] (4) at (0.75, -0.25) {$\alpha$};
					\node [style=none] (5) at (1.25, 0.5) {$A$};
				\end{pgfonlayer}
				\begin{pgfonlayer}{edgelayer}
					\draw (1.center) to (3);
					\draw (0.center) to (3);
				\end{pgfonlayer}
			\end{tikzpicture} =  \begin{tikzpicture}
				\begin{pgfonlayer}{nodelayer}
					\node [style=none] (0) at (0, 0.5) {};
					\node [style=none] (1) at (-1.5, 1) {};
					\node [style=none] (2) at (1, 0.5) {};
					\node [style=none] (3) at (1, -1.75) {};
					\node [style=circle, scale=2] (4) at (-0.75, -1) {};
					\node [style=none] (5) at (-0.75, -1) {$\alpha$};
					\node [style=none] (6) at (-1.5, -1.5) {$A'$};
					\node [style=none] (7) at (-2, 0.75) {$A$};
					\node [style=none] (8) at (1.5, -1.25) {$A'$};
					\node [style=circle] (9) at (-0.75, -0.25) {};
					\node [style=circle] (10) at (-0.75, -1.75) {};
					\node [style=none] (11) at (0.5, 1.25) {$\eta_L$};
				\end{pgfonlayer}
				\begin{pgfonlayer}{edgelayer}
					\draw [bend right=90, looseness=1.50] (2.center) to (0.center);
					\draw (3.center) to (2.center);
					\draw [in=15, out=-90, looseness=1.25] (0.center) to (9);
					\draw [in=-90, out=165, looseness=1.25] (9) to (1.center);
					\draw (10) to (4);
					\draw (9) to (4);
				\end{pgfonlayer}
			\end{tikzpicture} = \begin{tikzpicture}
			\begin{pgfonlayer}{nodelayer}
				\node [style=none] (0) at (-0.5, 0.25) {};
				\node [style=none] (1) at (1, 1) {};
				\node [style=none] (2) at (-1.5, 0.25) {};
				\node [style=none] (3) at (-1.5, -1.75) {};
				\node [style=circle, scale=2] (4) at (0.25, -1) {};
				\node [style=none] (5) at (0.25, -1) {$\alpha$};
				\node [style=none] (6) at (1, -1.5) {$A'$};
				\node [style=none] (7) at (1.25, 0.75) {$A$};
				\node [style=none] (8) at (-1.75, -1.25) {$A'$};
				\node [style=circle] (9) at (0.25, -0.25) {};
				\node [style=circle] (10) at (0.25, -1.75) {};
				\node [style=none] (11) at (-1, 1) {$\eta_R$};
			\end{pgfonlayer}
			\begin{pgfonlayer}{edgelayer}
				\draw (3.center) to (2.center);
				\draw [in=180, out=-90, looseness=1.00] (0.center) to (9);
				\draw [in=-90, out=15, looseness=1.25] (9) to (1.center);
				\draw (10) to (4);
				\draw (9) to (4);
				\draw [bend left=90, looseness=1.75] (2.center) to (0.center);
			\end{pgfonlayer}
		\end{tikzpicture}
		\end{equation}
	\end{lemma}
The equation in the previous Lemma reminds us of involutive 
monoids \cite[Theorem 5.28]{HeV19} in $\dagger$-monoidal categories. 

One can get Frobenius algebras by splitting binary idempotents on linear monoids:
\begin{lemma}
	\label{Lemma: Frobenius splitting}
	In an isomix category $\X$, let  $E \linmonb E'$ be a self-linear monoid in $\Core(\X)$ given by splitting a coring sectional or retractional 
	binary idempotent $(\u, \v)$ on linear monoid $A \linmonw B$. Let $\alpha: E \to E'$ be 
	the isomorphism. Then, $E$ is a Frobenius Algebra under the linear equivalence ${\sf Mx}_\downarrow$
	 if and only if the binary idempotent satisfies the following equation: 
	\begin{equation}
		\label{eqn: dag Frob split}
		\begin{tikzpicture}
			\begin{pgfonlayer}{nodelayer}
				\node [style=onehalfcircle] (44) at (6.25, 3) {};
				\node [style=none] (45) at (6.25, 0.75) {};
				\node [style=none] (46) at (6.25, 3) {$\u$};
				\node [style=none] (47) at (6.25, 4.75) {};
				\node [style=none] (48) at (6.5, 4.5) {$A$};
				\node [style=none] (49) at (6.5, 1) {$B$};
			\end{pgfonlayer}
			\begin{pgfonlayer}{edgelayer}
				\draw (45.center) to (44);
				\draw (47.center) to (44);
			\end{pgfonlayer}
		\end{tikzpicture} = \begin{tikzpicture}
			\begin{pgfonlayer}{nodelayer}
				\node [style=onehalfcircle] (37) at (3, 3.5) {};
				\node [style=none] (20) at (4, 4) {};
				\node [style=none] (22) at (5, 4) {};
				\node [style=none] (23) at (5, 0.75) {};
				\node [style=onehalfcircle] (24) at (3.5, 1.75) {};
				\node [style=none] (25) at (3.5, 1.75) {$\u$};
				\node [style=none] (26) at (3.25, 1.25) {$B$};
				\node [style=none] (27) at (2.75, 4.25) {$E$};
				\node [style=none] (28) at (5.25, 1.25) {$B$};
				\node [style=circle] (29) at (3.5, 2.5) {};
				\node [style=circle] (30) at (3.5, 0.75) {};
				\node [style=none] (31) at (4.5, 4.75) {$\eta_L$};
				\node [style=none] (36) at (3, 4.5) {};
				\node [style=none] (33) at (3, 3.5) {$e_A$};
				\node [style=onehalfcircle] (38) at (5, 3) {};
				\node [style=none] (39) at (5, 3) {$e_B$};
				\node [style=onehalfcircle] (40) at (4, 3.5) {};
				\node [style=none] (41) at (4, 3.5) {$e_A$};
				\node [style=none] (42) at (4, 3) {};
				\node [style=none] (43) at (3, 3) {};
			\end{pgfonlayer}
			\begin{pgfonlayer}{edgelayer}
				\draw [bend right=90, looseness=1.50] (22.center) to (20.center);
				\draw (37) to (36.center);
				\draw (38) to (22.center);
				\draw (20.center) to (40);
				\draw [in=-90, out=15, looseness=1.25] (29) to (42.center);
				\draw [in=165, out=-90, looseness=1.25] (43.center) to (29);
				\draw (38) to (23.center);
				\draw (40) to (42.center);
				\draw (37) to (43.center);
				\draw (30) to (24);
				\draw (24) to (29);
			\end{pgfonlayer}
		\end{tikzpicture} = \begin{tikzpicture}
		\begin{pgfonlayer}{nodelayer}
			\node [style=onehalfcircle] (37) at (5, 3.5) {};
			\node [style=none] (20) at (4, 4) {};
			\node [style=none] (22) at (3, 4) {};
			\node [style=none] (23) at (3, 0.75) {};
			\node [style=onehalfcircle] (24) at (4.5, 1.75) {};
			\node [style=none] (25) at (4.5, 1.75) {$\u$};
			\node [style=none] (26) at (4.75, 1.25) {$B$};
			\node [style=none] (27) at (5.25, 4.25) {$E$};
			\node [style=none] (28) at (2.75, 1.25) {$B$};
			\node [style=circle] (29) at (4.5, 2.5) {};
			\node [style=circle] (30) at (4.5, 0.75) {};
			\node [style=none] (31) at (3.5, 4.75) {$\eta_R$};
			\node [style=none] (36) at (5, 4.5) {};
			\node [style=none] (33) at (5, 3.5) {$e_A$};
			\node [style=onehalfcircle] (38) at (3, 3) {};
			\node [style=none] (39) at (3, 3) {$e_B$};
			\node [style=onehalfcircle] (40) at (4, 3.5) {};
			\node [style=none] (41) at (4, 3.5) {$e_A$};
			\node [style=none] (42) at (4, 3) {};
			\node [style=none] (43) at (5, 3) {};
		\end{pgfonlayer}
		\begin{pgfonlayer}{edgelayer}
			\draw [bend left=90, looseness=1.50] (22.center) to (20.center);
			\draw (37) to (36.center);
			\draw (38) to (22.center);
			\draw (20.center) to (40);
			\draw [in=-90, out=165, looseness=1.25] (29) to (42.center);
			\draw [in=15, out=-90, looseness=1.25] (43.center) to (29);
			\draw (38) to (23.center);
			\draw (40) to (42.center);
			\draw (37) to (43.center);
			\draw (30) to (24);
			\draw (24) to (29);
		\end{pgfonlayer}
	\end{tikzpicture}		
	\end{equation}
	where $e_A = \u \v$ and $e_B = \v \u$.
	\end{lemma}

 In a $\dagger$-isomix category, splitting a sectional or retractional $\dagger$-coring binary idempotent 
 on a $\dagger$-linear monoid  gives a $\dagger$-self-linear monoid on a pre-unitary object. 
 If the binary idempotent satisfies equation \ref{eqn: dag Frob split}, then, by using 
 Lemmas \ref{Lemma: dagmondagFrob} and  \ref{Lemma: Frobenius splitting}, one 
 gets a $\dagger$-Frobenius algebra on the splitting. 

\subsection{Linear comonoid}

The bialgebra law is a central ingredient of a complimentary system. The directionality of the linear distributors in an LDC  
forbids a bialgebraic interaction between two linear monoids.  A linear monoid, however, can interact bialgebraically 
with a linear comonoid.

\begin{definition}
	A {\bf linear comonoid}, $A \lincomonw B$, in an LDC consists of a $\ox$-comonoid, 
	$(A, \trianglecomult{0.55}, \trianglecounit{0.55})$, and a left and a right dual, $(\tau_L, \gamma_L):A \dashvv B$, 
	and $(\tau_R, \gamma_R): B \dashvv A$, such that:
	\begin{equation}
		\label{eqn: lin comon}
	(a)~~~ \begin{tikzpicture}
		\begin{pgfonlayer}{nodelayer}
			\node [style=none] (19) at (4.5, 1.25) {};
			\node [style=none] (20) at (5, 1.5) {$B$};
			\node [style=none] (21) at (4.5, 2) {};
			\node [style=none] (22) at (4.25, 2.25) {};
			\node [style=none] (23) at (4.75, 2.25) {};
			\node [style=oa] (26) at (4.4, 3.5) {};
			\node [style=none] (27) at (4.4, 4.25) {};
			\node [style=none] (28) at (3.75, 4.1) {$B \oa B$};
		\end{pgfonlayer}
		\begin{pgfonlayer}{edgelayer}
			\draw (19.center) to (21.center);
			\draw (21.center) to (22.center);
			\draw (22.center) to (23.center);
			\draw (23.center) to (21.center);
			\draw (26) to (27.center);
			\draw [bend right=45] (26) to (22.center);
			\draw [bend left=60] (26) to (23.center);
		\end{pgfonlayer}
	\end{tikzpicture} := \begin{tikzpicture}
			\begin{pgfonlayer}{nodelayer}
				\node [style=none] (0) at (1, 2) {};
				\node [style=none] (1) at (2, 2) {};
				\node [style=none] (2) at (1.5, 3.5) {};
				\node [style=none] (3) at (1.75, 3.25) {$A$};
				\node [style=none] (4) at (0.5, 2) {};
				\node [style=none] (5) at (-0.5, 2) {};
				\node [style=none] (6) at (2.5, 3.5) {};
				\node [style=none] (7) at (2.5, 1.25) {};
				\node [style=none] (8) at (2.75, 1.75) {$B$};
				\node [style=oa] (9) at (0, 3) {};
				\node [style=none] (10) at (0, 3.75) {};
				\node [style=none] (11) at (-0.65, 3.35) {$B \oa B$};
				\node [style=none] (12) at (1.5, 2.75) {};
				\node [style=none] (13) at (1.25, 2.5) {};
				\node [style=none] (14) at (1.75, 2.5) {};
				\node [style=none] (15) at (2, 4.25) {$\tau_L$};
				\node [style=none] (16) at (0.5, 1) {$\gamma_L$};
			\end{pgfonlayer}
			\begin{pgfonlayer}{edgelayer}
				\draw [bend left=90, looseness=1.75] (0.center) to (4.center);
				\draw [bend right=90] (5.center) to (1.center);
				\draw [bend left=270, looseness=1.50] (6.center) to (2.center);
				\draw (6.center) to (7.center);
				\draw [in=90, out=-45] (9) to (4.center);
				\draw [in=90, out=-135, looseness=1.25] (9) to (5.center);
				\draw (9) to (10.center);
				\draw [in=75, out=-165] (13.center) to (0.center);
				\draw [in=90, out=-30] (14.center) to (1.center);
				\draw (2.center) to (12.center);
				\draw (12.center) to (13.center);
				\draw (13.center) to (14.center);
				\draw (14.center) to (12.center);
			\end{pgfonlayer}
		\end{tikzpicture} = \begin{tikzpicture}
			\begin{pgfonlayer}{nodelayer}
				\node [style=none] (0) at (1.1, 2) {};
				\node [style=none] (1) at (0.0999999, 2) {};
				\node [style=none] (2) at (0.6, 3.5) {};
				\node [style=none] (3) at (0.35, 3.25) {$A$};
				\node [style=none] (4) at (1.6, 2) {};
				\node [style=none] (5) at (2.6, 2) {};
				\node [style=none] (6) at (-0.4, 3.5) {};
				\node [style=none] (7) at (-0.4, 1.25) {};
				\node [style=none] (8) at (-0.65, 1.75) {$B$};
				\node [style=oa] (9) at (2.1, 3) {};
				\node [style=none] (10) at (2.1, 3.75) {};
				\node [style=none] (11) at (2.75, 3.35) {$B \oa B$};
				\node [style=none] (12) at (0.6, 2.75) {};
				\node [style=none] (13) at (0.85, 2.5) {};
				\node [style=none] (14) at (0.35, 2.5) {};
				\node [style=none] (15) at (0.0999999, 4.25) {$\tau_R$};
				\node [style=none] (16) at (1.6, 1) {$\gamma_R$};
			\end{pgfonlayer}
			\begin{pgfonlayer}{edgelayer}
				\draw [bend right=90, looseness=1.75] (0.center) to (4.center);
				\draw [bend left=90] (5.center) to (1.center);
				\draw [bend right=270, looseness=1.50] (6.center) to (2.center);
				\draw (6.center) to (7.center);
				\draw [in=90, out=-135] (9) to (4.center);
				\draw [in=90, out=-45, looseness=1.25] (9) to (5.center);
				\draw (9) to (10.center);
				\draw [in=105, out=-15] (13.center) to (0.center);
				\draw [in=90, out=-150] (14.center) to (1.center);
				\draw (2.center) to (12.center);
				\draw (12.center) to (13.center);
				\draw (13.center) to (14.center);
				\draw (14.center) to (12.center);
			\end{pgfonlayer}
		\end{tikzpicture}
~~~~~~~~
(b)~~~ \begin{tikzpicture}
	\begin{pgfonlayer}{nodelayer}
		\node [style=none] (19) at (4.5, 1.25) {};
		\node [style=none] (20) at (5, 1.75) {$B$};
		\node [style=none] (21) at (4.5, 3.5) {};
		\node [style=none] (22) at (4.25, 3.75) {};
		\node [style=none] (23) at (4.75, 3.75) {};
		\node [style=none] (24) at (4.5, 4) {};
		\node [style=none] (25) at (4.5, 3.75) {};
		\node [style=none] (26) at (5, 4) {$\bot$};
	\end{pgfonlayer}
	\begin{pgfonlayer}{edgelayer}
		\draw (19.center) to (21.center);
		\draw (21.center) to (22.center);
		\draw (22.center) to (23.center);
		\draw (23.center) to (21.center);
		\draw (24.center) to (25.center);
	\end{pgfonlayer}
\end{tikzpicture} := 
\begin{tikzpicture}
	\begin{pgfonlayer}{nodelayer}
		\node [style=none] (0) at (1.5, 3.5) {};
		\node [style=none] (1) at (1.25, 3.5) {$A$};
		\node [style=none] (2) at (2.5, 3.5) {};
		\node [style=none] (3) at (2.5, 1) {};
		\node [style=none] (4) at (2.75, 1.75) {$B$};
		\node [style=circle, scale=1.5] (5) at (0.75, 1.25) {};
		\node [style=none] (6) at (0.75, 1.25) {$\bot$};
		\node [style=none] (7) at (0.75, 4.25) {};
		\node [style=circle] (8) at (0.75, 2) {};
		\node [style=none] (9) at (1.5, 2.75) {};
		\node [style=none] (10) at (1.5, 3) {};
		\node [style=none] (11) at (1.25, 2.75) {};
		\node [style=none] (12) at (1.75, 2.75) {};
		\node [style=none] (13) at (2, 4.25) {$\tau_L$};
	\end{pgfonlayer}
	\begin{pgfonlayer}{edgelayer}
		\draw [bend right=90, looseness=1.50] (2.center) to (0.center);
		\draw (2.center) to (3.center);
		\draw (5) to (7.center);
		\draw [dotted, in=-90, out=0, looseness=1.25] (8) to (9.center);
		\draw (10.center) to (11.center);
		\draw (11.center) to (12.center);
		\draw (12.center) to (10.center);
		\draw (10.center) to (0.center);
	\end{pgfonlayer}
\end{tikzpicture}
 = \begin{tikzpicture}
	\begin{pgfonlayer}{nodelayer}
		\node [style=none] (0) at (2, 3.5) {};
		\node [style=none] (1) at (2.25, 3.5) {$A$};
		\node [style=none] (2) at (1, 3.5) {};
		\node [style=none] (3) at (1, 1) {};
		\node [style=none] (4) at (0.75, 1.75) {$B$};
		\node [style=circle, scale=1.5] (5) at (2.75, 1.25) {};
		\node [style=none] (6) at (2.75, 1.25) {$\bot$};
		\node [style=none] (7) at (2.75, 4.25) {};
		\node [style=circle] (8) at (2.75, 2) {};
		\node [style=none] (9) at (2, 2.75) {};
		\node [style=none] (10) at (2, 3) {};
		\node [style=none] (11) at (2.25, 2.75) {};
		\node [style=none] (12) at (1.75, 2.75) {};
		\node [style=none] (13) at (1.5, 4.25) {$\tau_R$};
	\end{pgfonlayer}
	\begin{pgfonlayer}{edgelayer}
		\draw [bend left=90, looseness=1.50] (2.center) to (0.center);
		\draw (2.center) to (3.center);
		\draw (5) to (7.center);
		\draw [dotted, in=-90, out=180, looseness=1.25] (8) to (9.center);
		\draw (10.center) to (11.center);
		\draw (11.center) to (12.center);
		\draw (12.center) to (10.center);
		\draw (10.center) to (0.center);
	\end{pgfonlayer}
\end{tikzpicture}
\end{equation}
\end{definition}

Note that while a linear monoid has a $\ox$-monoid and a $\oa$-comonoid, a linear comonoid has a $\ox$-comonoid and an $\oa$-monoid. 

A {\bf morphism} of linear comonoids, $(f,g) : (A \lincomonw B) \to (A' \lincomonw B')$, consists of a pair of maps,
$f: A \to A'$ and $g: B' \to B$, such that $f$ is a comonoid morphism, and $(f,g)$ and $(g,f)$ are morphisms of 
the left and the right duals respectively.  

In a monoidal category, an idempotent $e: A \to A$ 
is {\bf sectional} (respectively {\bf retractional}) on a comonoid $(A, d, k)$  
if $e d = e d (e \ox e)$ (respectively if $d (e \ox e) = e d (e \ox e)$ and $e k = k$). 
In an LDC, a binary idempotent $(\u,\v)$ is {\bf sectional} (respectively {\bf retractional}) on a linear monoid when 
$e_A=\u\v$ and $e_B=\v\u$ satisfies the conditions in  the table below. 
 \[ \begin{tabular}{l|| l} 
	\hline
		(\u, \v) \textbf{sectional} on $A \lincomonw B$  & $(\u, \v)$ \textbf{retractional} on $A \lincomonw B$\\
	\hline
		$e_A$ preserves $(A,d,k)$ sectionally &  $e_A$ preserves $(A,d,k)$ retractionally \\
	\hline
		$(e_A, e_B)$ preserves $(\eta_L, \epsilon_L): A \dashvv B$ sectionally &$(e_A, e_B)$ preserves $(\eta_L, \epsilon_L): A \dashvv B$ retractionally\\
	\hline
		$(e_B, e_A)$ preserves $(\eta_R, \epsilon_R): B \dashvv A$ retractionally & $(e_B, e_A)$ preserves $(\eta_R, \epsilon_R): B \dashvv A$ sectionally\\
	\hline
	\end{tabular} \] 
Splitting a sectional or retractional binary idempotent on a linear comonoid gives a self-linear comonoid.
\begin{definition}
A {\bf $\dagger$-linear comonoid} $A \dagcomonwtik A^\dagger$ in a $\dagger$-LDC is a linear comonoid 
$A \lincomonw A^\dagger$ such that $(\tau_L, \gamma_L): A \dashvv A^\dagger$ and $(\tau_R, \gamma_R): 
A^\dagger \dashvv A$ are $\dagger$-duals, and:
\begin{equation}
    \label{Eqn: leftdagcomon}
	\begin{tikzpicture}
		\begin{pgfonlayer}{nodelayer}
			\node [style=none] (19) at (4.5, 1.25) {};
			\node [style=none] (20) at (4.85, 1.75) {$A^\dagger$};
			\node [style=none] (21) at (4.5, 3.5) {};
			\node [style=none] (22) at (4.25, 3.75) {};
			\node [style=none] (23) at (4.75, 3.75) {};
			\node [style=none] (24) at (4.5, 4) {};
			\node [style=none] (25) at (4.5, 3.75) {};
			\node [style=none] (26) at (5, 4) {$\bot$};
		\end{pgfonlayer}
		\begin{pgfonlayer}{edgelayer}
			\draw (19.center) to (21.center);
			\draw (21.center) to (22.center);
			\draw (22.center) to (23.center);
			\draw (23.center) to (21.center);
			\draw (24.center) to (25.center);
		\end{pgfonlayer}
	\end{tikzpicture}  := 
    \begin{tikzpicture}
	\begin{pgfonlayer}{nodelayer}
		\node [style=none] (0) at (1.5, 3.5) {};
		\node [style=none] (1) at (1.25, 3.5) {$A$};
		\node [style=none] (2) at (2.5, 3.5) {};
		\node [style=none] (3) at (2.5, 1) {};
		\node [style=none] (4) at (2.75, 1.75) {$A^\dagger$};
		\node [style=circle, scale=1.5] (5) at (0.75, 1.25) {};
		\node [style=none] (6) at (0.75, 1.25) {$\bot$};
		\node [style=none] (7) at (0.75, 4.25) {};
		\node [style=circle] (8) at (0.75, 2) {};
		\node [style=none] (9) at (1.5, 2.75) {};
		\node [style=none] (10) at (1.5, 3) {};
		\node [style=none] (11) at (1.25, 2.75) {};
		\node [style=none] (12) at (1.75, 2.75) {};
	\end{pgfonlayer}
	\begin{pgfonlayer}{edgelayer}
		\draw [bend right=90, looseness=1.25] (2.center) to (0.center);
		\draw (2.center) to (3.center);
		\draw (5) to (7.center);
		\draw [dotted, in=-90, out=0, looseness=1.25] (8) to (9.center);
		\draw (10.center) to (11.center);
		\draw (11.center) to (12.center);
		\draw (12.center) to (10.center);
		\draw (10.center) to (0.center);
	\end{pgfonlayer}
\end{tikzpicture} = \begin{tikzpicture}
	\begin{pgfonlayer}{nodelayer}
		\node [style=none] (0) at (0, 2) {};
		\node [style=none] (1) at (0, 2.5) {};
		\node [style=none] (2) at (-0.75, 2.5) {};
		\node [style=none] (3) at (-0.75, 1.25) {};
		\node [style=none] (4) at (0.75, 1.25) {};
		\node [style=none] (5) at (0.75, 2.5) {};
		\node [style=none] (6) at (0, 0) {};
		\node [style=none] (7) at (0, 1.25) {};
		\node [style=none] (8) at (0.5, 0.25) {$A^\dagger$};
		\node [style=none] (9) at (-0.25, 1.75) {};
		\node [style=none] (10) at (0.25, 1.75) {};
		\node [style=none] (11) at (0, 1.75) {};
		\node [style=none] (12) at (0.5, 1.25) {};
		\node [style=circle, scale=1.5] (13) at (1.5, 0.25) {};
		\node [style=none] (14) at (1.5, 3.75) {};
		\node [style=circle] (15) at (1.5, 3) {};
		\node [style=none] (16) at (0.5, 2.5) {};
		\node [style=none] (17) at (1.5, 0.25) {$\bot$};
	\end{pgfonlayer}
	\begin{pgfonlayer}{edgelayer}
		\draw (1.center) to (0.center);
		\draw (2.center) to (5.center);
		\draw (5.center) to (4.center);
		\draw (4.center) to (3.center);
		\draw (3.center) to (2.center);
		\draw (6.center) to (7.center);
		\draw (9.center) to (10.center);
		\draw (10.center) to (0.center);
		\draw (0.center) to (9.center);
		\draw [dotted, bend right=15, looseness=1.00] (11.center) to (12.center);
		\draw (14.center) to (13);
		\draw [dotted, in=180, out=90, looseness=1.25] (16.center) to (15);
	\end{pgfonlayer}
\end{tikzpicture} ~~~~~~~~ 	\begin{tikzpicture}
	\begin{pgfonlayer}{nodelayer}
		\node [style=none] (19) at (4.5, 1.25) {};
		\node [style=none] (20) at (5, 1.5) {$A^\dagger$};
		\node [style=none] (21) at (4.5, 2) {};
		\node [style=none] (22) at (4.25, 2.25) {};
		\node [style=none] (23) at (4.75, 2.25) {};
		\node [style=oa] (26) at (4.4, 3.5) {};
		\node [style=none] (27) at (4.4, 4.25) {};
		\node [style=none] (28) at (3.6, 4.1) {$A^\dagger \oa A^\dagger$};
	\end{pgfonlayer} 
	\begin{pgfonlayer}{edgelayer}
		\draw (19.center) to (21.center);
		\draw (21.center) to (22.center);
		\draw (22.center) to (23.center);
		\draw (23.center) to (21.center);
		\draw (26) to (27.center);
		\draw [bend right=45] (26) to (22.center);
		\draw [bend left=60] (26) to (23.center);
	\end{pgfonlayer}
\end{tikzpicture} := \begin{tikzpicture}
	\begin{pgfonlayer}{nodelayer}
		\node [style=none] (0) at (1, 2) {};
		\node [style=none] (1) at (2, 2) {};
		\node [style=none] (2) at (1.5, 3.5) {};
		\node [style=none] (3) at (1.75, 3.25) {$A$};
		\node [style=none] (4) at (0.5, 2) {};
		\node [style=none] (5) at (-0.5, 2) {};
		\node [style=none] (6) at (2.5, 3.5) {};
		\node [style=none] (7) at (2.5, 1.25) {};
		\node [style=none] (8) at (2.75, 1.75) {$A^\dagger$};
		\node [style=oa] (9) at (0, 3) {};
		\node [style=none] (10) at (0, 3.75) {};
		\node [style=none] (11) at (-0.75, 3.35) {$A^\dagger \oa A^\dagger$};
		\node [style=none] (12) at (1.5, 2.75) {};
		\node [style=none] (13) at (1.25, 2.5) {};
		\node [style=none] (14) at (1.75, 2.5) {};
	\end{pgfonlayer}
	\begin{pgfonlayer}{edgelayer}
		\draw [bend left=90, looseness=1.75] (0.center) to (4.center);
		\draw [bend right=90, looseness=1.00] (5.center) to (1.center);
		\draw [bend right=90, looseness=1.25] (6.center) to (2.center);
		\draw (6.center) to (7.center);
		\draw [in=90, out=-45, looseness=1.00] (9) to (4.center);
		\draw [in=90, out=-135, looseness=1.25] (9) to (5.center);
		\draw (9) to (10.center);
		\draw [in=75, out=-165, looseness=1.00] (13.center) to (0.center);
		\draw [in=90, out=-30, looseness=1.00] (14.center) to (1.center);
		\draw (2.center) to (12.center);
		\draw (12.center) to (13.center);
		\draw (13.center) to (14.center);
		\draw (14.center) to (12.center);
	\end{pgfonlayer}
\end{tikzpicture}
	= \begin{tikzpicture}
	\begin{pgfonlayer}{nodelayer}
		\node [style=none] (0) at (0.75, -0.25) {};
		\node [style=none] (1) at (0.25, -1) {$A^\dagger$};
		\node [style=none] (2) at (0, 1.5) {};
		\node [style=none] (3) at (0.25, 0.75) {};
		\node [style=none] (4) at (-1, -0.25) {};
		\node [style=none] (5) at (-1, 1.5) {};
		\node [style=none] (6) at (1, 1.5) {};
		\node [style=none] (7) at (1, -0.25) {};
		\node [style=none] (8) at (0.5, 1.5) {};
		\node [style=none] (9) at (-0.5, 1.5) {};
		\node [style=none] (10) at (0, -0.25) {};
		\node [style=none] (11) at (0, -1) {};
		\node [style=oa] (12) at (0, 2.25) {};
		\node [style=none] (13) at (0, 3) {};
		\node [style=none] (14) at (0.65, 2.75) {$A^\dagger \oa A^\dagger$};
		\node [style=none] (15) at (0, 1) {};
		\node [style=none] (16) at (-0.25, 0.75) {};
		\node [style=none] (17) at (-0.25, 0.75) {};
		\node [style=none] (18) at (-0.75, -0.25) {};
	\end{pgfonlayer}
	\begin{pgfonlayer}{edgelayer}
		\draw [in=90, out=-30, looseness=1.25] (3.center) to (0.center);
		\draw (5.center) to (4.center);
		\draw (4.center) to (7.center);
		\draw (6.center) to (7.center);
		\draw (6.center) to (5.center);
		\draw (11.center) to (10.center);
		\draw [bend left, looseness=1.00] (9.center) to (12);
		\draw [bend left, looseness=1.00] (12) to (8.center);
		\draw (12) to (13.center);
		\draw [in=90, out=-150, looseness=1.25] (17.center) to (18.center);
		\draw (2.center) to (15.center);
		\draw (15.center) to (16.center);
		\draw (16.center) to (3.center);
		\draw (3.center) to (15.center);
	\end{pgfonlayer}
\end{tikzpicture}
\end{equation}
A {\bf $\dagger$-self-linear comonoid} consists of an 
isomorphism $\alpha: A \to A^\dagger$ such that $\alpha \alpha^{-1 \dagger} = \iota$. 
A {\bf morphism of $\dagger$-linear comonoids} is a 
pair $(f, f^\dagger)$ such that $(f, f^\dagger)$ is a morphism of the 
underlying linear comonoids. 
\end{definition}

In a $\dagger$-LDC, splitting a $\dagger$-binary idempotent on a $\dagger$-linear comonoid 
gives a $\dagger$-self-linear comonoid when the binary idempotent is either sectional or retractional. 
In the next section, we discuss linear bialgebras which are given by an interacting linear monoid and linear comonoid.

\subsection{Linear bialgebras}
All the results concerning bialgebras are necessarily set in symmetric LDCs and we shall assume that linear monoids
 and the linear comonoids are symmetric.
\begin{definition}
	\label{Defn: linear bialg}
	A {\bf linear bialgebra}, $\frac{(a,b)}{(a',b')}: A \linbialgwtik B$, in an LDC consists of  
	a linear monoid, $(a,b):\! A \linmonw B$ and  a linear comonoid, $(a',b') : A \lincomonwtritik B$ such that 
$(A, \mulmap{1.2}{white}, \unitmap{1.2}{white}, \trianglecomult{0.55} , \trianglecounit{0.55})$ and 
and  $(B, \trianglemult{0.55}, \triangleunit{0.55}, \comulmap{1.2}{white} , \counitmap{1.2}{white})$ are $\ox$- and $\oa$-bialgebras respectively. 
A {\bf morphism} of linear bialgebras is a morphism both of the linear monoids and linear comonoids.
\end{definition}	

A linear bialgebra is {\bf commutative} if the $\oa$-monoid and $\ox$-monoid are commutative.
A {\bf self-linear bialgebra} is a linear bialgebra in which there is an isomorphism $A \to^{\alpha} B$ 
(so essentially the algebra is on one object).

A binary idempotent on a linear bialgebra is {\bf sectional} (respectively {\bf retractional}) if it is sectional (respectively retractional) 
on the linear monoid, and the linear comonoid. In an LDC, splitting a sectional or retractional binary idempotent on a linear bialgebra induces a 
self-linear bialgebra on the splitting. 

\begin{definition}
A {\bf $\dagger$-linear bialgebra}, $\frac{(a,b)}{(a',b')} : A \dagbialgwtik A^\dagger$, is a linear bialgebra 
with a $\dagger$-linear monoid and a $\dagger$-linear comonoid. 
A {\bf $\dagger$-self-linear bialgebra} is $\dagger$-linear bialgebra which is also a self-linear bialgebra  such 
that the isomorphism, $\alpha: A \to A^\dagger$, satisfies $\alpha \alpha^{-1 \dagger} = \iota$. 
\end{definition}
Note that $A$ is a weak preunitary object: if it  was in the core as well, it would be a preunitary object. 
In a $\dagger$-LDC, splitting a $\dagger$-binary idempotent on a $\dagger$-linear bialgebra gives a 
$\dagger$-self-linear bialgebra if the idempotent is either a sectional or retractional. 

\subsection{Complementary systems}

In quantum mechanics, two quantum observables are complementary \cite{Gri18} if measuring one observable increases 
the uncertainty regarding the value of the other. Complementarity is a key feature distinguishing classical from quantum mechanics.
In CQM, the complimentarity principle is described using interacting commutative $\dagger$-Frobenius algebras. This section describes complementarity in isomix categories:
\begin{definition}
A {\bf complementary system} in an isomix category, $\X$, is a commutative and cocommutative 
self-linear bialgebra, $\frac{(a,b)}{(a',b')}: A \linbialgwtik A$ such that the following equations (with their `op' symmetries) hold:
\[ \mbox{ \small \bf [comp.1]} ~~~ \begin{tikzpicture}
	\begin{pgfonlayer}{nodelayer}
		\node [style=none] (20) at (6, 4.25) {};
		\node [style=none] (21) at (6.25, 4.5) {};
		\node [style=none] (22) at (5.75, 4.5) {};
		\node [style=none] (23) at (6, 3.5) {};
		\node [style=none] (24) at (6, 3.5) {};
		\node [style=circle] (26) at (6, 2.5) {};
		\node [style=none] (27) at (5.75, 3.5) {};
		\node [style=none] (28) at (5, 4.5) {};
		\node [style=none] (30) at (6, 3.75) {};
		\node [style=none] (31) at (6, 3.25) {};
	\end{pgfonlayer}
	\begin{pgfonlayer}{edgelayer}
		\draw (20.center) to (21.center);
		\draw (20.center) to (22.center);
		\draw (22.center) to (21.center);
		\draw (20.center) to (23.center);
		\draw (24.center) to (26);
		\draw [in=-90, out=-180, looseness=1.25] (27.center) to (28.center);
		\draw [bend right=90, looseness=1.75] (30.center) to (31.center);
	\end{pgfonlayer}
\end{tikzpicture}  = \begin{tikzpicture}
	\begin{pgfonlayer}{nodelayer}
		\node [style=none] (0) at (0.75, 3) {};
		\node [style=none] (1) at (1, 2.75) {};
		\node [style=none] (2) at (0.5, 2.75) {};
		\node [style=none] (3) at (0.75, 4.5) {};
	\end{pgfonlayer}
	\begin{pgfonlayer}{edgelayer}
		\draw (0.center) to (1.center);
		\draw (0.center) to (2.center);
		\draw (2.center) to (1.center);
		\draw (0.center) to (3.center);
	\end{pgfonlayer}
\end{tikzpicture} 
~~~~~~~~	
\mbox{ \small \bf [comp.2]} ~~~ \begin{tikzpicture}
		\begin{pgfonlayer}{nodelayer}
			\node [style=none] (0) at (-1.75, 4.5) {};
			\node [style=none] (1) at (-1.25, 4.5) {};
			\node [style=none] (2) at (-1.5, 4.25) {};
			\node [style=none] (3) at (-1.5, 3.75) {};
			\node [style=none] (4) at (-1.5, 3.75) {};
			\node [style=none] (5) at (-1.75, 3.5) {};
			\node [style=none] (6) at (-1.25, 3.5) {};
			\node [style=circle] (7) at (-0.75, 2.5) {};
			\node [style=none] (8) at (-2.25, 2.5) {};
		\end{pgfonlayer}
		\begin{pgfonlayer}{edgelayer}
			\draw (2.center) to (0.center);
			\draw (2.center) to (1.center);
			\draw (1.center) to (0.center);
			\draw (3.center) to (2.center);
			\draw (4.center) to (5.center);
			\draw (4.center) to (6.center);
			\draw (6.center) to (5.center);
			\draw [in=180, out=90, looseness=1.00] (8.center) to (5.center);
			\draw [in=90, out=-15, looseness=1.25] (6.center) to (7);
		\end{pgfonlayer}
	\end{tikzpicture} =  \begin{tikzpicture}
		\begin{pgfonlayer}{nodelayer}
			\node [style=none] (0) at (-0.75, 2.5) {};
			\node [style=circle] (1) at (-0.75, 4.5) {};
		\end{pgfonlayer}
		\begin{pgfonlayer}{edgelayer}
			\draw (0.center) to (1);
		\end{pgfonlayer}
	\end{tikzpicture} 
~~~~~~~~ 	\mbox{ \small \bf [comp.3]} ~~~ 
\begin{tikzpicture}
	\begin{pgfonlayer}{nodelayer}
		\node [style=none] (0) at (2.75, 0.5) {};
		\node [style=none] (1) at (3, 0.75) {};
		\node [style=none] (2) at (2.5, 0.75) {};
		\node [style=none] (3) at (2.75, -0.25) {};
		\node [style=none] (4) at (2.75, 0) {};
		\node [style=none] (5) at (2.5, -0.25) {};
		\node [style=none] (6) at (2.75, -1.25) {};
		\node [style=none] (7) at (1.75, -1.25) {};
		\node [style=none] (8) at (2.75, 0) {};
		\node [style=none] (9) at (2.75, -0.5) {};
	\end{pgfonlayer}
	\begin{pgfonlayer}{edgelayer}
		\draw (0.center) to (1.center);
		\draw (0.center) to (2.center);
		\draw (2.center) to (1.center);
		\draw (0.center) to (3.center);
		\draw (6.center) to (3.center);
		\draw [in=90, out=-180, looseness=1.50] (5.center) to (7.center);
		\draw [bend right=90, looseness=1.75] (8.center) to (9.center);
	\end{pgfonlayer}
\end{tikzpicture}
 = \begin{tikzpicture}
		\begin{pgfonlayer}{nodelayer}
			\node [style=none] (0) at (1.75, 2.5) {};
			\node [style=none] (1) at (2.25, 2.5) {};
			\node [style=none] (2) at (2, 2.25) {};
			\node [style=none] (3) at (2, 0.5) {};
			\node [style=none] (4) at (2.75, 2.25) {};
			\node [style=none] (5) at (2.75, 0.5) {};
			\node [style=none] (6) at (2.5, 2.5) {};
			\node [style=none] (7) at (3, 2.5) {};
		\end{pgfonlayer}
		\begin{pgfonlayer}{edgelayer}
			\draw[fill=white] (1.center) -- (0.center) -- (2.center) -- (1.center);
			\draw (3.center) to (2.center);
			\draw (4.center) to (6.center);
			\draw (4.center) to (7.center);
			\draw (7.center) to (6.center);
			\draw (5.center) to (4.center);
		\end{pgfonlayer}
	\end{tikzpicture} \]
	A {\bf $\dagger$-complementary system} in a $\dagger$-isomix category is a $\dagger$-self-linear bialgebra which 
	is also a complementary system. 
\end{definition}

Notice that we are using the alternative presentation of linear monoids by actions and coactions (see 
Proposition \ref{Lemma: alternate presentation of linear monoids}).  Thus, {\bf [comp.1]} requires that the counit of the  linear comonoid to be
 dual to the counit via the linear monoid dual, while {\bf [comp.2]} requires that the unit of the linear monoid to be dual to the 
counit via dual of the linear comonoid.  Finally, {\bf [comp.3] } requires that the coaction map of the linear monoid 
 duplicates the unit of $\dagger$-linear comonoid.  The `co' symmetry of the equations are immediate 
from the commutativity and cocommutativity of the linear bialgebra. The `op' symmetry of equations 
 holds automatically for a $\dagger$-complimentary system. 

\begin{lemma}
	If $A \linbialgwtik A$ is a complementary system in an isomix category,  then
	A  is a $\ox$-bialgebra with antipode given by $(a)$ and a $\oa$-bialgebra with antipode 
	given by  $(b)$:
	\begin{center}
	$(a)$ ~~~~ \begin{tikzpicture} [scale=0.75]
		\begin{pgfonlayer}{nodelayer}
			\node [style=none] (0) at (1, 4.75) {};
			\node [style=none] (1) at (1.5, 4.75) {};
			\node [style=circle] (2) at (1.25, 2.75) {};
			\node [style=none] (3) at (1.25, 5) {};
			\node [style=none] (4) at (1.5, 5.75) {};
			\node [style=none] (5) at (1, 5.75) {};
			\node [style=none] (6) at (1.25, 5.5) {};
			\node [style=circle] (7) at (1.25, 1.75) {};
			\node [style=none] (8) at (0.25, 1.75) {};
			\node [style=none] (9) at (0.25, 5.75) {};
		\end{pgfonlayer}
		\begin{pgfonlayer}{edgelayer}
			\draw (1.center) to (0.center);
			\draw (0.center) to (3.center);
			\draw (3.center) to (1.center);
			\draw [bend left=45, looseness=1.00] (1.center) to (2);
			\draw (4.center) to (5.center);
			\draw (5.center) to (6.center);
			\draw (6.center) to (4.center);
			\draw (3.center) to (6.center);
			\draw (7) to (2);
			\draw [in=90, out=-135, looseness=1.00] (0.center) to (8.center);
			\draw [in=-90, out=135, looseness=1.00] (2) to (9.center);
		\end{pgfonlayer}
	\end{tikzpicture}
	~~~~~~~~
	$(b)$~~~ \begin{tikzpicture}  [scale=0.75]
		\begin{pgfonlayer}{nodelayer}
			\node [style=none] (0) at (1, 2.75) {};
			\node [style=none] (1) at (1.5, 2.75) {};
			\node [style=circle] (2) at (1.25, 4.75) {};
			\node [style=none] (3) at (1.25, 2.5) {};
			\node [style=none] (4) at (1.5, 1.75) {};
			\node [style=none] (5) at (1, 1.75) {};
			\node [style=none] (6) at (1.25, 2) {};
			\node [style=circle] (7) at (1.25, 5.75) {};
			\node [style=none] (8) at (0.25, 5.75) {};
			\node [style=none] (9) at (0.25, 1.75) {};
		\end{pgfonlayer}
		\begin{pgfonlayer}{edgelayer}
			\draw (1.center) to (0.center);
			\draw (0.center) to (3.center);
			\draw (3.center) to (1.center);
			\draw [bend right=45, looseness=1.00] (1.center) to (2);
			\draw (4.center) to (5.center);
			\draw (5.center) to (6.center);
			\draw (6.center) to (4.center);
			\draw (3.center) to (6.center);
			\draw (7) to (2);
			\draw [in=-90, out=135, looseness=0.75] (0.center) to (8.center);
			\draw [in=90, out=-150, looseness=1.00] (2) to (9.center);
		\end{pgfonlayer}
	\end{tikzpicture} 
\end{center}
	\end{lemma}
	\begin{proof}
	Given a complementary system we show the $\ox$-bialgebra has an antipode so is a $\ox$-Hopf algebra:
	\[ \begin{tikzpicture}
		\begin{pgfonlayer}{nodelayer}
			\node [style=none] (0) at (0, 5) {};
			\node [style=none] (1) at (-0.25, 4.75) {};
			\node [style=none] (2) at (0.25, 4.75) {};
			\node [style=circle] (3) at (0, 2.75) {};
			\node [style=circle, scale=1.5] (4) at (0.5, 3.75) {};
			\node [style=none] (5) at (0.5, 3.75) {$s$};
			\node [style=none] (6) at (0, 5.75) {};
			\node [style=none] (7) at (0, 5.75) {};
			\node [style=none] (8) at (0, 1.75) {};
		\end{pgfonlayer}
		\begin{pgfonlayer}{edgelayer}
			\draw (1.center) to (2.center);
			\draw (2.center) to (0.center);
			\draw (0.center) to (1.center);
			\draw [bend right=45, looseness=1.00] (1.center) to (3);
			\draw [in=-90, out=45, looseness=1.00] (3) to (4);
			\draw [bend left=15, looseness=1.00] (2.center) to (4);
			\draw (6.center) to (0.center);
			\draw (3) to (8.center);
		\end{pgfonlayer}
	\end{tikzpicture} := \begin{tikzpicture}
		\begin{pgfonlayer}{nodelayer}
			\node [style=none] (0) at (0, 5) {};
			\node [style=none] (1) at (-0.25, 4.75) {};
			\node [style=none] (2) at (0.25, 4.75) {};
			\node [style=circle] (3) at (0, 2.75) {};
			\node [style=none] (4) at (0, 5.75) {};
			\node [style=none] (5) at (0, 5.75) {};
			\node [style=none] (6) at (0, 1.75) {};
			\node [style=none] (7) at (1, 4.75) {};
			\node [style=none] (8) at (1.5, 4.75) {};
			\node [style=circle] (9) at (1.25, 2.75) {};
			\node [style=none] (10) at (1.25, 5) {};
			\node [style=none] (11) at (1.5, 5.75) {};
			\node [style=none] (12) at (1, 5.75) {};
			\node [style=none] (13) at (1.25, 5.5) {};
			\node [style=circle] (14) at (1.25, 1.75) {};
		\end{pgfonlayer}
		\begin{pgfonlayer}{edgelayer}
			\draw (1.center) to (2.center);
			\draw (2.center) to (0.center);
			\draw (0.center) to (1.center);
			\draw [bend right=45, looseness=1.00] (1.center) to (3);
			\draw (4.center) to (0.center);
			\draw (3) to (6.center);
			\draw (8.center) to (7.center);
			\draw (7.center) to (10.center);
			\draw (10.center) to (8.center);
			\draw [bend left=45, looseness=1.00] (8.center) to (9);
			\draw (11.center) to (12.center);
			\draw (12.center) to (13.center);
			\draw (13.center) to (11.center);
			\draw (10.center) to (13.center);
			\draw (14) to (9);
			\draw (7.center) to (3);
			\draw (2.center) to (9);
		\end{pgfonlayer}
	\end{tikzpicture} = \begin{tikzpicture}
		\begin{pgfonlayer}{nodelayer}
			\node [style=none] (0) at (0, 3.25) {};
			\node [style=none] (1) at (-0.25, 3) {};
			\node [style=none] (2) at (0.25, 3) {};
			\node [style=circle] (3) at (0, 4) {};
			\node [style=none] (4) at (0, 3.25) {};
			\node [style=none] (5) at (-0.75, 5) {};
			\node [style=circle] (6) at (0.75, 2) {};
			\node [style=none] (7) at (-0.75, 2) {};
			\node [style=none] (8) at (0.5, 5) {};
			\node [style=none] (9) at (1, 5) {};
			\node [style=none] (10) at (0.75, 4.75) {};
		\end{pgfonlayer}
		\begin{pgfonlayer}{edgelayer}
			\draw (1.center) to (2.center);
			\draw (2.center) to (0.center);
			\draw (0.center) to (1.center);
			\draw (3) to (4.center);
			\draw (10.center) to (8.center);
			\draw (8.center) to (9.center);
			\draw (9.center) to (10.center);
			\draw [in=150, out=-90, looseness=1.00] (5.center) to (3);
			\draw [in=-90, out=30, looseness=1.00] (3) to (10.center);
			\draw [in=90, out=-150, looseness=0.75] (1.center) to (7.center);
			\draw [in=90, out=-30, looseness=1.25] (2.center) to (6);
		\end{pgfonlayer}
	\end{tikzpicture} =\begin{tikzpicture}
		\begin{pgfonlayer}{nodelayer}
			\node [style=none] (0) at (0, 3.25) {};
			\node [style=none] (1) at (-0.25, 3) {};
			\node [style=none] (2) at (0.25, 3) {};
			\node [style=circle] (3) at (0, 4.5) {};
			\node [style=none] (4) at (0, 3.25) {};
			\node [style=none] (5) at (-0.75, 5.5) {};
			\node [style=circle] (6) at (0.75, 2) {};
			\node [style=none] (7) at (-0.75, 2) {};
			\node [style=none] (8) at (0.5, 5.5) {};
			\node [style=none] (9) at (1, 5.5) {};
			\node [style=none] (10) at (0.75, 5.25) {};
			\node [style=none] (11) at (0, 4) {};
			\node [style=none] (12) at (0, 3.75) {};
			\node [style=none] (13) at (-0.25, 3.75) {};
			\node [style=circle] (14) at (-1.5, 2) {};
			\node [style=none] (15) at (0, 4) {};
			\node [style=none] (16) at (0, 3.5) {};
		\end{pgfonlayer}
		\begin{pgfonlayer}{edgelayer}
			\draw (1.center) to (2.center);
			\draw (2.center) to (0.center);
			\draw (0.center) to (1.center);
			\draw (3) to (4.center);
			\draw (10.center) to (8.center);
			\draw (8.center) to (9.center);
			\draw (9.center) to (10.center);
			\draw [in=150, out=-90, looseness=1.25] (5.center) to (3);
			\draw [in=-90, out=30] (3) to (10.center);
			\draw [in=90, out=-150, looseness=0.75] (1.center) to (7.center);
			\draw [in=90, out=-30, looseness=1.25] (2.center) to (6);
			\draw [in=90, out=-165] (13.center) to (14);
			\draw [bend right=90, looseness=1.75] (15.center) to (16.center);
		\end{pgfonlayer}
	\end{tikzpicture}	
	 = \begin{tikzpicture}
		\begin{pgfonlayer}{nodelayer}
			\node [style=circle] (0) at (3.25, 3) {};
			\node [style=none] (1) at (4.25, 4.75) {};
			\node [style=none] (2) at (4.25, 5.25) {};
			\node [style=none] (3) at (4.25, 4.5) {};
			\node [style=none] (4) at (4, 4.5) {};
			\node [style=none] (5) at (4.25, 3.25) {};
			\node [style=none] (6) at (4, 5.5) {};
			\node [style=none] (7) at (4.5, 5.5) {};
			\node [style=none] (10) at (3, 3.75) {};
			\node [style=none] (11) at (2.5, 5.5) {};
			\node [style=none] (12) at (4.25, 3.25) {};
			\node [style=none] (13) at (4.25, 3.25) {};
			\node [style=circle] (14) at (5, 2) {};
			\node [style=none] (15) at (4, 3) {};
			\node [style=none] (16) at (3.5, 2) {};
			\node [style=none] (17) at (4.5, 3) {};
			\node [style=none] (18) at (4.25, 4.75) {};
			\node [style=none] (19) at (4.25, 4.25) {};
			\node [style=none] (20) at (3.25, 4) {};
			\node [style=none] (21) at (3.25, 3.5) {};
		\end{pgfonlayer}
		\begin{pgfonlayer}{edgelayer}
			\draw [in=90, out=-180, looseness=1.50] (4.center) to (0);
			\draw (2.center) to (5.center);
			\draw (6.center) to (7.center);
			\draw (7.center) to (2.center);
			\draw (2.center) to (6.center);
			\draw [in=-90, out=165] (10.center) to (11.center);
			\draw (15.center) to (17.center);
			\draw (17.center) to (12.center);
			\draw (12.center) to (15.center);
			\draw [in=90, out=-150, looseness=0.75] (15.center) to (16.center);
			\draw [in=90, out=-30, looseness=1.25] (17.center) to (14);
			\draw [bend right=90, looseness=1.75] (18.center) to (19.center);
			\draw [bend right=90, looseness=1.75] (20.center) to (21.center);
		\end{pgfonlayer}
	\end{tikzpicture} \stackrel{\bf [comp.3]}{=}
	\begin{tikzpicture}
		\begin{pgfonlayer}{nodelayer}
			\node [style=circle] (0) at (3, 3) {};
			\node [style=none] (1) at (4, 4) {};
			\node [style=none] (2) at (4, 3.5) {};
			\node [style=none] (3) at (3.75, 4.25) {};
			\node [style=none] (4) at (4.25, 4.25) {};
			\node [style=none] (5) at (3, 3.75) {};
			\node [style=none] (6) at (3, 3.5) {};
			\node [style=none] (7) at (2.75, 3.75) {};
			\node [style=none] (8) at (2, 5.25) {};
			\node [style=none] (9) at (4, 3.5) {};
			\node [style=none] (10) at (4, 3.5) {};
			\node [style=circle] (11) at (4.5, 2.25) {};
			\node [style=none] (12) at (3.75, 3.25) {};
			\node [style=none] (13) at (3.5, 2) {};
			\node [style=none] (14) at (4.25, 3.25) {};
			\node [style=none] (15) at (3, 4.75) {};
			\node [style=none] (16) at (3.25, 5) {};
			\node [style=none] (17) at (2.75, 5) {};
			\node [style=none] (18) at (3, 4) {};
			\node [style=none] (19) at (3, 3.5) {};
		\end{pgfonlayer}
		\begin{pgfonlayer}{edgelayer}
			\draw (1.center) to (2.center);
			\draw (3.center) to (4.center);
			\draw (4.center) to (1.center);
			\draw (1.center) to (3.center);
			\draw [in=-90, out=-180, looseness=1.25] (7.center) to (8.center);
			\draw (12.center) to (14.center);
			\draw (14.center) to (9.center);
			\draw (9.center) to (12.center);
			\draw [in=90, out=-150, looseness=0.75] (12.center) to (13.center);
			\draw [in=90, out=-45] (14.center) to (11);
			\draw (17.center) to (16.center);
			\draw (16.center) to (15.center);
			\draw (15.center) to (17.center);
			\draw (15.center) to (0);
			\draw [bend right=90, looseness=1.75] (18.center) to (19.center);
		\end{pgfonlayer}
	\end{tikzpicture}  \stackrel{\bf [comp.2]}{=}
	\begin{tikzpicture}
		\begin{pgfonlayer}{nodelayer}
			\node [style=circle] (0) at (3, 3) {};
			\node [style=none] (5) at (3, 3.75) {};
			\node [style=none] (6) at (3, 3.5) {};
			\node [style=none] (7) at (2.75, 3.75) {};
			\node [style=none] (8) at (2, 5.25) {};
			\node [style=circle] (11) at (3.5, 4.25) {};
			\node [style=none] (13) at (3.5, 2) {};
			\node [style=none] (15) at (3, 4.75) {};
			\node [style=none] (16) at (3.25, 5) {};
			\node [style=none] (17) at (2.75, 5) {};
			\node [style=none] (18) at (3, 4) {};
			\node [style=none] (19) at (3, 3.5) {};
		\end{pgfonlayer}
		\begin{pgfonlayer}{edgelayer}
			\draw [in=-90, out=-180, looseness=1.25] (7.center) to (8.center);
			\draw (17.center) to (16.center);
			\draw (16.center) to (15.center);
			\draw (15.center) to (17.center);
			\draw (15.center) to (0);
			\draw [bend right=90, looseness=1.75] (18.center) to (19.center);
			\draw (11) to (13.center);
		\end{pgfonlayer}
	\end{tikzpicture}	
	\stackrel{ \bf [comp.1]}{=} 
	\begin{tikzpicture}
		\begin{pgfonlayer}{nodelayer}
			\node [style=circle] (0) at (3, 3) {};
			\node [style=none] (1) at (3, 4) {};
			\node [style=none] (2) at (3, 5) {};
			\node [style=none] (3) at (2.75, 3.75) {};
			\node [style=none] (4) at (3.25, 3.75) {};
			\node [style=none] (5) at (3, 2) {};
		\end{pgfonlayer}
		\begin{pgfonlayer}{edgelayer}
			\draw (1.center) to (2.center);
			\draw (3.center) to (4.center);
			\draw (4.center) to (1.center);
			\draw (1.center) to (3.center);
			\draw (5.center) to (0);
		\end{pgfonlayer}
	\end{tikzpicture} \]
	
	Similarly, the $\oa$-bialgebra has an antipode using the `op' versions of {\bf[comp.1]} to {\bf[comp.3]}. 
	\end{proof}

	A $\dagger$-complimentary system in a unitary category corresponds to the usual notion 
	of interacting commutative $\dagger$-Frobenius algebras \cite{CoD11} when its linear monoid and linear 
	comonoid satisfy condition \ref{eqnn: unitary coincidence}. Splitting binary idempotents on 
	a linear bialgebra produces a complimentary system under the following conditions:

	\begin{lemma} 
	\label{Lemma: complementary idempotent}
	In an isomix category, a self linear bialgebra given by  splitting a coring binary idempotent $(\u, \v)$ on a commutative 
and cocommutative linear bialgebra $A \linbialgwtik B$ is a complimentary system if and only if the binary idempotent 
satisfies the following conditions (and their `op' symmetric forms):
\begin{equation}
	\label{eqn: idemcomp}
(a) ~ \begin{tikzpicture}
	\begin{pgfonlayer}{nodelayer}
		\node [style=none] (0) at (2, 5) {};
		\node [style=none] (1) at (2.25, 5.25) {};
		\node [style=none] (2) at (1.75, 5.25) {};
		\node [style=none] (3) at (2, 3.25) {};
		\node [style=none] (4) at (2, 3) {};
		\node [style=none] (5) at (2, 3.25) {};
		\node [style=circle] (6) at (2, 1.5) {};
		\node [style=none] (7) at (1.75, 3.25) {};
		\node [style=none] (9) at (0.75, 5) {$A$};
		\node [style=none] (10) at (2.5, 3.25) {$B$};
		\node [style=none] (11) at (1, 5.25) {};
		\node [style=circle, scale=1.8] (12) at (2, 4.25) {};
		\node [style=none] (13) at (2, 4.25) {$e_B$};
		\node [style=none] (14) at (2.5, 4.75) {$B$};
		\node [style=circle, scale=1.8] (15) at (2, 2.25) {};
		\node [style=none] (16) at (2, 2.25) {$e_B$};
		\node [style=circle, scale=1.8] (17) at (1, 4.25) {};
		\node [style=none] (18) at (1, 4.25) {$e_A$};
		\node [style=none] (19) at (2, 3.5) {};
		\node [style=none] (20) at (2, 3) {};
	\end{pgfonlayer}
	\begin{pgfonlayer}{edgelayer}
		\draw (0.center) to (1.center);
		\draw (0.center) to (2.center);
		\draw (2.center) to (1.center);
		\draw (4.center) to (5.center);
		\draw (0.center) to (12);
		\draw (4.center) to (15);
		\draw (6) to (15);
		\draw (12) to (5.center);
		\draw [in=180, out=-90, looseness=1.25] (17) to (7.center);
		\draw (11.center) to (17);
		\draw [bend right=90, looseness=1.75] (19.center) to (20.center);
	\end{pgfonlayer}
\end{tikzpicture} = \begin{tikzpicture}
	\begin{pgfonlayer}{nodelayer}
		\node [style=none] (19) at (4, 2) {};
		\node [style=none] (20) at (4.25, 1.75) {};
		\node [style=none] (21) at (3.75, 1.75) {};
		\node [style=none] (22) at (4, 2) {};
		\node [style=none] (23) at (3.5, 5) {$A$};
		\node [style=none] (24) at (4, 5.5) {};
		\node [style=circle, scale=1.8] (25) at (4, 3.75) {};
		\node [style=none] (26) at (4, 3.75) {$e_A$};
	\end{pgfonlayer}
	\begin{pgfonlayer}{edgelayer}
		\draw (19.center) to (20.center);
		\draw (19.center) to (21.center);
		\draw (21.center) to (20.center);
		\draw (24.center) to (25);
		\draw (25) to (22.center);
	\end{pgfonlayer}
\end{tikzpicture}
~~~~~~~~ (b) ~
\begin{tikzpicture}
	\begin{pgfonlayer}{nodelayer}
		\node [style=none] (33) at (5.5, 3.75) {};
		\node [style=none] (34) at (6, 3.75) {};
		\node [style=none] (36) at (5.75, 4) {};
		\node [style=none] (37) at (5.75, 5.25) {};
		\node [style=circle] (38) at (5, 1.75) {};
		\node [style=circle, scale=1.5] (39) at (5.75, 4.75) {};
		\node [style=none] (40) at (5.75, 4.75) {$\v$};
		\node [style=none] (41) at (6.25, 5) {$B$};
		\node [style=circle, scale=1.5] (42) at (5, 2.75) {};
		\node [style=none] (43) at (5, 2.75) {$\u$};
		\node [style=none] (44) at (5.5, 5.5) {};
		\node [style=none] (45) at (6, 5.5) {};
		\node [style=none] (47) at (6.25, 4.25) {$A$};
		\node [style=none] (48) at (6.5, 1. 5) {};
		\node [style=circle, scale=1.5] (50) at (6.5, 2.75) {};
		\node [style=none] (51) at (6.5, 2.75) {$e_A$};
		\node [style=none] (52) at (7, 2.25) {$A$};
	\end{pgfonlayer}
	\begin{pgfonlayer}{edgelayer}
		\draw (34.center) to (33.center);
		\draw (33.center) to (36.center);
		\draw (36.center) to (34.center);
		\draw (37.center) to (39);
		\draw (36.center) to (39);
		\draw (38) to (42);
		\draw (44.center) to (45.center);
		\draw (45.center) to (37.center);
		\draw (37.center) to (44.center);
		\draw (48.center) to (50);
		\draw [in=90, out=-150, looseness=1.25] (33.center) to (42);
		\draw [in=90, out=-30, looseness=1.25] (34.center) to (50);
	\end{pgfonlayer}
\end{tikzpicture} = \begin{tikzpicture}
	\begin{pgfonlayer}{nodelayer}
		\node [style=none] (39) at (5.5, 2.5) {$A$};
		\node [style=none] (45) at (5, 1.75) {};
		\node [style=circle] (46) at (5, 5.5) {};
		\node [style=circle, scale=1.5] (47) at (5, 3.75) {};
		\node [style=none] (48) at (5, 3.75) {$e_A$};
	\end{pgfonlayer}
	\begin{pgfonlayer}{edgelayer}
		\draw (45.center) to (47);
		\draw (47) to (46);
	\end{pgfonlayer}
\end{tikzpicture}
~~~~~~~~ 	(c) ~
\begin{tikzpicture}
	\begin{pgfonlayer}{nodelayer}
		\node [style=none] (49) at (7.25, 3.5) {};
		\node [style=none] (50) at (7, 3.25) {};
		\node [style=none] (51) at (7, 1.75) {};
		\node [style=none] (52) at (7, 3.75) {};
		\node [style=none] (53) at (7, 5.25) {};
		\node [style=none] (54) at (8, 1.75) {};
		\node [style=circle, scale=1.5] (55) at (7, 4.5) {};
		\node [style=none] (56) at (7, 4.5) {$\v$};
		\node [style=none] (57) at (8.5, 2) {$B$};
		\node [style=none] (58) at (7.25, 5.5) {};
		\node [style=none] (59) at (6.75, 5.5) {};
		\node [style=none] (60) at (6.5, 2) {$B$};
		\node [style=none] (61) at (6.5, 5) {$B$};
		\node [style=circle, scale=1.5] (62) at (7, 2.5) {};
		\node [style=none] (63) at (7, 2.5) {$\u$};
		\node [style=none] (64) at (8, 3) {};
		\node [style=none] (65) at (6.5, 3.5) {$A$};
		\node [style=circle, scale=1.5] (66) at (8, 2.5) {};
		\node [style=none] (67) at (8, 2.5) {$e_B$};
		\node [style=none] (68) at (7, 3.75) {};
		\node [style=none] (69) at (7, 3.25) {};
	\end{pgfonlayer}
	\begin{pgfonlayer}{edgelayer}
		\draw (52.center) to (50.center);
		\draw (53.center) to (55);
		\draw (52.center) to (55);
		\draw (58.center) to (59.center);
		\draw (59.center) to (53.center);
		\draw (53.center) to (58.center);
		\draw (51.center) to (62);
		\draw (50.center) to (62);
		\draw [in=0, out=105] (64.center) to (49.center);
		\draw (54.center) to (66);
		\draw (64.center) to (66);
		\draw [bend left=90, looseness=1.75] (68.center) to (69.center);
	\end{pgfonlayer}
\end{tikzpicture}  =  \begin{tikzpicture}
	\begin{pgfonlayer}{nodelayer}
		\node [style=none] (68) at (9.75, 2) {};
		\node [style=none] (69) at (9.75, 3.75) {};
		\node [style=none] (70) at (9.75, 5.25) {};
		\node [style=circle, scale=1.5] (71) at (9.75, 4.5) {};
		\node [style=none] (72) at (9.75, 4.5) {$\v$};
		\node [style=circle, scale=1.5] (73) at (9.75, 3.25) {};
		\node [style=none] (74) at (9.75, 3.25) {$\u$};
		\node [style=none] (75) at (9.5, 5.5) {};
		\node [style=none] (76) at (10, 5.5) {};
		\node [style=circle, scale=1.5] (77) at (10.75, 4.5) {};
		\node [style=none] (78) at (10.75, 4.5) {$\v$};
		\node [style=none] (79) at (11.25, 2.25) {$B$};
		\node [style=circle, scale=1.5] (80) at (10.75, 3.25) {};
		\node [style=none] (81) at (11, 5.5) {};
		\node [style=none] (82) at (10.5, 5.5) {};
		\node [style=none] (83) at (10.75, 5.25) {};
		\node [style=none] (84) at (10.75, 2) {};
		\node [style=none] (85) at (11.25, 5) {$B$};
		\node [style=none] (86) at (10.75, 3.75) {};
		\node [style=none] (87) at (11.25, 4) {$A$};
		\node [style=none] (88) at (10.75, 3.25) {$\u$};
	\end{pgfonlayer}
	\begin{pgfonlayer}{edgelayer}
		\draw (70.center) to (71);
		\draw (69.center) to (71);
		\draw (75.center) to (76.center);
		\draw (76.center) to (70.center);
		\draw (70.center) to (75.center);
		\draw (68.center) to (73);
		\draw (69.center) to (73);
		\draw (83.center) to (77);
		\draw (86.center) to (77);
		\draw (82.center) to (81.center);
		\draw (81.center) to (83.center);
		\draw (83.center) to (82.center);
		\draw (84.center) to (80);
		\draw (86.center) to (80);
	\end{pgfonlayer}
\end{tikzpicture}
\end{equation}
where $e_A = \u \v$, and $e_B = \v \u$.
\end{lemma}

In the next section, we provide an example of the previous Lemma using 
exponential modalities.


\section{Exponential modalities}
\label{Sec: exponential modalities}

An LDC has exponential modalities if it is equipped with a linear comonad 
$((!,?),(\epsilon,\eta),(\delta,\mu))$ \cite{BCS96}. The linearity of the functors in a 
$(!,?)$-LDC means that $(!,\delta,\epsilon)$ is a monoidal comonad while $(?,\mu,\eta)$ 
is a comonoidal monad, and $(!(A),\Delta_A,\trianglecounit{0.55}_A)$ is a natural 
cocommutative comonoid while $(?(A),\nabla_A,\triangleunit{0.55}_A)$ is a natural commutative monoid. 
A $\dagger$-$(!,?)$-LDC is a $(!,?)$-LDC in which all the functors and natural transformations are 
$\dagger$-linear (see \cite{CCS18}).

In a $(!,?)$-LDC, any dual, $(\alpha, \beta): A \dashvv ~B$, induces a dual, $(\alpha_!, \beta_?): !A \dashvv ~?B$ 
(see the below diagrams), on the exponential modalities using the 
linearity of $(!,?)$. This means that any dual induces a linear comonoid, $(\alpha_!, \beta_?): !A \lincomonw ?B$,
where the comonoid structure is given by the modalities.

{ \centering $ \alpha_! := \begin{tikzpicture}
	\begin{pgfonlayer}{nodelayer}
		\node [style=none] (23) at (3, 2) {};
		\node [style=none] (24) at (5, 2) {};
		\node [style=none] (25) at (3, 1) {};
		\node [style=none] (26) at (5, 1) {};
		\node [style=none] (27) at (3.25, 1) {};
		\node [style=none] (28) at (3.75, 1) {};
		\node [style=none] (29) at (4.5, 1) {};
		\node [style=none] (30) at (3.5, 1) {};
		\node [style=none] (31) at (3.5, 0.5) {};
		\node [style=none] (32) at (4.5, 0.5) {};
		\node [style=none] (33) at (4, 1.75) {$\alpha$};
		\node [style=none] (34) at (4.75, 1.25) {$!$};
		\node [style=none] (35) at (3, 0.5) {$!A$};
		\node [style=none] (36) at (5, 0.5) {$?B$};
	\end{pgfonlayer}
	\begin{pgfonlayer}{edgelayer}
		\draw [bend left=75, looseness=1.25] (27.center) to (28.center);
		\draw [bend left=90, looseness=2.00] (30.center) to (29.center);
		\draw (23.center) to (24.center);
		\draw (24.center) to (26.center);
		\draw (26.center) to (25.center);
		\draw (25.center) to (23.center);
		\draw (31.center) to (30.center);
		\draw (32.center) to (29.center);
	\end{pgfonlayer}
\end{tikzpicture} = m_\top (!\alpha)\nu_\ox  ~~~~~~~~ \beta_? := \begin{tikzpicture}
	\begin{pgfonlayer}{nodelayer}
		\node [style=none] (23) at (5, 0.5) {};
		\node [style=none] (24) at (3, 0.5) {};
		\node [style=none] (25) at (5, 1.5) {};
		\node [style=none] (26) at (3, 1.5) {};
		\node [style=none] (27) at (3.75, 1.5) {};
		\node [style=none] (28) at (3.25, 1.5) {};
		\node [style=none] (29) at (3.5, 1.5) {};
		\node [style=none] (30) at (4.5, 1.5) {};
		\node [style=none] (31) at (4.5, 2) {};
		\node [style=none] (32) at (3.5, 2) {};
		\node [style=none] (33) at (4, 0.75) {$\beta$};
		\node [style=none] (34) at (4.75, 0.75) {$?$};
		\node [style=none] (35) at (5, 2) {$!A$};
		\node [style=none] (36) at (3, 2) {$?B$};
	\end{pgfonlayer}
	\begin{pgfonlayer}{edgelayer}
		\draw [bend left=75, looseness=1.25] (27.center) to (28.center);
		\draw [bend left=90, looseness=1.75] (30.center) to (29.center);
		\draw (23.center) to (24.center);
		\draw (24.center) to (26.center);
		\draw (26.center) to (25.center);
		\draw (25.center) to (23.center);
		\draw (31.center) to (30.center);
		\draw (32.center) to (29.center);
	\end{pgfonlayer}
\end{tikzpicture} =  \nu_\oa (?\epsilon) n_\bot ~~~~~~~~ m_{F_\ox} := \begin{tikzpicture}
	\begin{pgfonlayer}{nodelayer}
		\node [style=none] (23) at (3, 2) {};
		\node [style=none] (24) at (5, 2) {};
		\node [style=none] (25) at (3, 1) {};
		\node [style=none] (26) at (5, 1) {};
		\node [style=none] (29) at (4.5, 2.5) {};
		\node [style=none] (30) at (3.5, 2.5) {};
		\node [style=none] (31) at (4.75, 1.25) {$F_\ox$};
		\node [style=circle] (32) at (4, 1.5) {};
		\node [style=none] (33) at (4, 0.5) {};
		\node [style=none] (34) at (3.75, 1) {};
		\node [style=none] (35) at (4.25, 1) {};
		\node [style=none] (36) at (2.8, 2.25) {$F_\ox(A)$};
		\node [style=none] (37) at (5.2, 2.25) {$F_\ox(A)$};
		\node [style=none] (38) at (3.3, 0.5) {$F_\ox(A)$};
	\end{pgfonlayer}
	\begin{pgfonlayer}{edgelayer}
		\draw (23.center) to (24.center);
		\draw (24.center) to (26.center);
		\draw (26.center) to (25.center);
		\draw (25.center) to (23.center);
		\draw [bend right, looseness=1.25] (30.center) to (32);
		\draw [bend right, looseness=1.25] (32) to (29.center);
		\draw (33.center) to (32);
		\draw [bend left=60, looseness=1.50] (34.center) to (35.center);
	\end{pgfonlayer}
\end{tikzpicture}  = m_\ox F_\ox(m) $ \par }

Any linear functor $(F_\ox,F_\oa)$ applied to a linear monoid $(\alpha,\beta):A \linmonw B$ always produces a linear monoid 
$(\alpha_F,\beta_F): F_\ox(A) \expmonwtik F_\oa(B)$ with multplication $m_F$ as shown in the right diagram above.  This simple observation when applied to the 
exponential modalities has a striking effect: 
\begin{lemma} 
	\label{Lemma: !? linear monoid}
	In any $(!,?)$-LDC any linear monoid $(a, b): A \linmonw B$ 
	and an arbitrary dual $(a', b') :  A \dashvv B$ give a linear bialgebra
	$\frac{(a_!, b_!)}{(a'_!, b'_?)} : !A \expbialgwtik ?B$ using the natural cocommutative comonoid 
	$(!A, \Delta_A, \tricounit{0.65})$. 
\end{lemma}
	
The bialgebra structure results from the naturality of $\Delta$ and $\trianglecounit{0.55}$ over the 
functorially induced monoid structure.

A (!,?)-LDC has {\bf free exponential modalities} if, for any object $A$, $(!A, \Delta_A, \trianglecounit{0.55}_A)$ is a cofree 
cocommutative comonoid, and $(?(A), \nabla_A, \triangleunit{0.55}_A)$ is a free commutative monoid \cite{Laf88}.
An example of a $\dagger$-LDC with free ($\dagger$-)exponential modalites is finiteness matrices over the complex numbers, $\FMat(\C)$. 
Moreover, $\FMat(\C)$, is a $\dagger$-isomix category and gives a key example of a MUC as 
 discussed in \cite{CCS18} (although exponentials are not discussed).  
 The universal property of free exponential modalities in a $(!,?)$-LDC implies the following:
\begin{lemma}
	\label{Lemma: free !? linear comonoid}
If $(f,g)$ is a morphism of duals, then the unique map $(f^\flat, g^\sharp)$ induced by the 
universal property of the free exponential is a morphism of linear comonoids.
\end{lemma}
This is illustrated by the commuting diagram $(a)$, below.
\[ (a) \!\! \xymatrixrowsep{10mm} \xymatrix{
	(x,y)  :  X \! \lincomonbtik \! Y \ar@{..>}[d]_{(f^\flat,g^\sharp)}  \ar[drr]^{(f,g)} \\
	(a_!,b_?) \! : \! !A \! \lincomonwtik \! ?B  \ar[rr]_{(\epsilon,\eta)} & & (a,b) \! : \! A \! \dashvv \! B }
(b)  \!\! \xymatrixcolsep{6mm} \xymatrixrowsep{8mm} \xymatrix{
	\frac{(x,y)}{(x',y')}  : X \! \linbialgbtik \! Y \ar@{..>}[d]_{(f^\flat,g^\sharp)}  \ar[drr]^{(f,g)} \\
	\frac{(a_!,b_?)}{(a'_!, b'_?)}\! :\! !A \!\expbialgwtik \! ?B  \ar[rr]_{(\epsilon,\eta)~~~~~~} & & (a,b) \!: \!A  \!\linmonwtik  \! B; (a', b') \!: \!A \dashvv B } \]
The results discussed so far can be combined to give, as shown in diagram $(b)$ above, a more complicated observation:
\begin{proposition} 
	\label{Prop: free !? linear bialgebra}
In a $(!,?)$-LDC with free exponential modalities, let $\frac{(x,y)}{(x',y')}:X \linbialgbtik Y$ be a linear bialgebra, $(a,b): A \linmonwtik B$ a linear 
monoid, and $(a',b'): A \dashvv B$ a dual,  then
\[ (f^\flat,g^\sharp): \left(\frac{(x,y)}{(x',y')}:X \linbialgbtik Y \right) \to \left(\frac{(a_!,b_?)}{(a'_!,b'_?)}:!A \expbialgwtik ?B \right)\] 
is a morphism of bialgebras, whenever $f: (X, \mulmap{1.5}{black}, \unitmap{1.5}{black}) \to 
(A, \mulmap{1.5}{white}, \unitmap{1.5}{white})$ 
is a morphism of monoids, and $(f,g)$ is a morphism of both duals:
 \[ (f,g): ((x,y): X \linmonw Y) \to ((a,b): A \linmonw B) ~~~ \text{ and }~~~ (f,g): ((x',y'): X \dashvv Y) \to ((a',b') : A \dashvv B) \] 
\end{proposition}

\begin{corollary} \label{Corr: splitting expmod}
	In a $(!,?)$-LDC with free exponential modalities, 
	if $A \linbialgbtik B$ is a linear bialgebra then $(1^\flat,1^\sharp): (A \linbialgbtik B) \to (!A\expbialgbtik ?A)$ 
	is a morphism of bialgebras, making $A \linbialgbtik B$ a retract of $!A\expbialgbtik ?A$.
\end{corollary}

 The corollary shows that every self-linear bialgebra in a (!,?)-LDC, with free exponential modalities, induces a 
 sectional binary idempotent on the induced linear bialgebra on the exponential modalities:
 
{ \centering $ \xymatrix{  !A \ar@<0.5ex>[r]^{\epsilon}  
	 & A \simeq B \ar@<0.5ex>[r]^{\eta}  \ar@<0.5ex>[l]^{1^\flat}
	 & ?B \ar@<0.5ex>[l]^{1^\sharp} } $ \par }

\medskip

Combining Corollary \ref{Corr: splitting expmod} and Lemma \ref{Lemma: complementary idempotent}, we get: 

\begin{theorem}
In a $(!,?)$-isomix category with free exponential modalities, every complimentary system
arises as a  splitting of a sectional binary idempotent on the free exponential modalities. 
\end{theorem}

The above results extend directly to $\dagger$-linear bialgebras in $\dagger$-LDCs with free exponential modalities due to the 
$\dagger$-linearity of $(!,?)$, $(\eta, \epsilon)$, $(\Delta, \nabla)$, and $(\tricounit{0.65}, \triunit{0.65})$. 

\section{Conclusion}

Bohr's principle of complementarity \cite{Gri18} states that, due to the wave and particle nature of matter, physical properties occur in 
complimentary pairs.  In the formulation of measurements in a MUC, a measurement on $A$
induces a measurement on $A^\dagger$, and vice versa.  A measurement transfers the structures of $A$ and $A^\dagger$ --- 
and the interactions between these --- onto a single compact object.  Our main result displays a complementary system 
as the result of a measurement of a $\dagger$-linear bialgebra in which two distinct dual structures  have been ``compacted'' into one structure. 
This provides an interesting perspective on Bohr's principle.

\bibliographystyle{eptcs}
\bibliography{measurement}



\end{document}